\documentclass[smallcondensed]{svjour3}       % onecolumn (second format)
\smartqed  % flush right qed marks, e.g. at end of proof

\usepackage{geometry}
\usepackage{xcolor}

\definecolor{matlab_blue}{rgb}{        0 ,   0.4470 ,   0.7410}
\definecolor{matlab_red}{rgb}{1.0,    0.0,    0.0}
\definecolor{matlab_yellow}{rgb}{0.9290,    0.6940 ,   0.1250}
\definecolor{matlab_green}{rgb}{0.4660,    0.6740,    0.1880}
\usepackage{graphicx}
\usepackage{subfigure}

\usepackage{amssymb,amsmath,sansmath,lmodern}
\usepackage{color}

\usepackage{accents}
\usepackage{stmaryrd}
\usepackage{tabu}

\usepackage[numbers]{natbib}

\newcommand{\aver}[1]{\ensuremath{\{\!\{#1\}\!\}}}
\newcommand{\jump}[1]{\ensuremath{\left\llbracket #1 \right\rrbracket}}

\newcommand\R{\mathbb R}

\newcommand\Rey{\mathrm{Re}}

\newcommand*\diff{\mathop{}\!\mathrm{d}}

\DeclareMathAccent{\svec}{\mathord}{letters}{126}
\newcommand\acclrvec[1]{\accentset{\,\leftrightarrow}{#1}}
\newcommand\stvec[1]{\mathbf #1}				
\newcommand\stvecg[1]{\boldsymbol #1}				
\newcommand\ssvec[1]{\acclrvec{\stvec{#1}}}	
\newcommand\ssvecg[1]{\acclrvec{\stvecg{#1}}}	
\newcommand\cssvec[1]{\acclrvec{\tilde{\stvec{#1}}}} 
\newcommand\cssvecg[1]{\acclrvec{\tilde{\stvecg{ #1}}}} 
\newcommand{\smat}[1]{\underline{\stvec{#1}}}

\newcommand{\snabla}{\svec{\nabla}}
\newcommand\interiorfaces{{\mathrm{interior}\atop\mathrm{faces}}}
\newcommand\boundaryfaces{{\mathrm{boundary}\atop\mathrm{faces}}}

\usepackage{xargs}    
\usepackage[colorinlistoftodos,prependcaption,textsize=tiny]{todonotes}
\newcommandx{\unsure}[2][1=]{\todo[linecolor=blue,backgroundcolor=blue!25,bordercolor=blue,#1]{#2}}
\newcommandx{\change}[2][1=]{\todo[linecolor=red,backgroundcolor=red!25,bordercolor=red,#1]{#2}}
\newcommandx{\info}[2][1=]{\todo[linecolor=OliveGreen,backgroundcolor=OliveGreen!25,bordercolor=OliveGreen,#1]{#2}}
\newcommandx{\improvement}[2][1=]{\todo[linecolor=Plum,backgroundcolor=Plum!25,bordercolor=Plum,#1]{#2}}
\newcommandx{\thiswillnotshow}[2][1=]{\todo[disable,#1]{#2}}

\usepackage{scalerel,stackengine}
\stackMath
\newcommand\reallywidehat[1]{%
\savestack{\tmpbox}{\stretchto{%
  \scaleto{%
    \scalerel*[\widthof{\ensuremath{#1}}]{\kern.1pt\mathchar"0362\kern.1pt}%
    {\rule{0ex}{\textheight}}%WIDTH-LIMITED CIRCUMFLEX
  }{\textheight}% 
}{2.4ex}}%
\stackon[-6.9pt]{#1}{\tmpbox}%
}
\usepackage{tikz}
\usetikzlibrary{plotmarks}

\newtheorem{prop}{Property}
\newtheorem{rem}{Remark}

\begin{document}

\title{Entropy--stable discontinuous Galerkin approximation with summation--by--parts property for the incompressible Navier--Stokes/Cahn--Hilliard system}%\thanks{Grants or other notes}
%\subtitle{Do you have a subtitle?\\ If so, write it here}

\titlerunning{ES--DGSEM with SBP for the incompressible NS/CH system}%\thanks{Grants or other notes}        % if too long for running head

\author{Juan Manzanero        \and Gonzalo Rubio \and \\
David A. Kopriva \and Esteban Ferrer \and 
Eusebio Valero 
}

%\authorrunning{Short form of author list} % if too long for running head

\institute{Juan Manzanero (\email{juan.manzanero@upm.es}), Gonzalo Rubio, Esteban Ferrer, Eusebio Valero  \at
	ETSIAE-UPM - School of Aeronautics, Universidad Polit\'ecnica de Madrid. Plaza Cardenal Cisneros 3, E-28040 Madrid, Spain. //
	Center for Computational Simulation, Universidad Polit\'ecnica de Madrid, Campus de Montegancedo, Boadilla del Monte, 28660, Madrid, Spain. \\
              David A. Kopriva \at Department of Mathematics, Florida State University and Computational Science Research Center, San Diego State University. 
}
\date{Received: date / Accepted: date}
% The correct dates will be entered by the editor

\maketitle

\begin{abstract}
We develop an entropy stable two--phase 
incompressible Navier--Stokes/Cahn--Hilliard discontinuous Galerkin (DG) flow solver method. The model poses
 the Cahn--Hilliard equation as the phase field method, a skew--symmetric 
form of the momentum equation, and an artificial compressibility method to compute the 
pressure. We design the model so that it satisfies an entropy law, including 
free-- and no--slip wall boundary conditions with non--zero 
wall contact angle.
 We then construct a high--order DG approximation of 
the model that satisfies the SBP--SAT property. 
With the help of a discrete 
stability analysis, the scheme has two modes: an entropy 
conserving approximation with central advective fluxes and the Bassi--Rebay 1 (BR1) method for diffusion, and an entropy 
stable approximation with an exact Riemann solver for advection and interface stabilization added to the BR1 
method. The scheme is applicable to, and the stability proofs hold for, three--dimensional unstructured meshes 
with curvilinear hexahedral elements.
%Two choices are proposed to integrate in time: an explicit third--order 
%Runge--Kutta scheme, and an implicit--explicit BDF method with first or second 
%order of accuracy. 
We test the convergence 
of the schemes on a manufactured solution, and their robustness by solving a flow 
initialized from random numbers. In the latter, we find that a similar scheme that does 
not satisfy an entropy inequality had 30$\%$ probability to fail, while 
the entropy stable scheme never does. 
We also solve the static and rising bubble test problems, and to challenge 
the solver capabilities we compute a three--dimensional pipe flow in the annular 
regime.

% \PACS{PACS code1 \and PACS code2 \and more}
% \subclass{MSC code1 \and MSC code2 \and more}
\end{abstract}

\section{Introduction}\label{sec:Introduction}

The study of multiphase flows is of broad interest, from both scientific and industrial perspectives. In particular, the oil industry investigates two--phase flows of immiscible fluids (e.g. oil and water). 
Without mixing, the typical flow configuration is the segregation of the two fluids, separated by a thin interfacial region that behaves like a permeable membrane.
At the interface of two dissimilar fluids, the forces acting on the molecules of each of these fluids 
are not the same as within each phase and generate interfacial tension.

Among the various techniques to approach the solution of immiscible two--phase flows, one finds two broad categories: 
sharp and diffuse interface models. 
In the former, the interface is considered infinitely thin, and acts as a physical boundary condition where the two fluids 
are coupled and balanced by interfacial forces \cite{2007:Sussman}. 
An example of a sharp interface model is the level--set method \cite{2005:Olsson}, where 
the interface is tracked by an additional variable advected by the flow velocity field. 

Diffuse interface approaches, however, regularize the problem with the introduction of
 an interface with non--zero thickness and a smooth transition of thermodynamic variables between both fluids. 
This time, interfacial forces (which in the sharp interface approach are delta functions 
 applied at the interface) are transformed to body forces whose effect is concentrated 
 at the interface \cite{1998:Lowengrub}. 
 
 In this work we study the diffuse interface model 
 of Cahn--Hilliard \cite{1958:Cahn,1959:Cahn} combined with the incompressible 
 Navier--Stokes equations with variable density and artificial 
 (also called \textit{pseudo}) compressibility 
 \cite{1996:Shen}. A review of alternate 
 Navier--Stokes/Cahn--Hilliard models can be found in \cite{2017:Hosseini}.

We discretize the system of partial differential equations that represents the model with a high--order Discontinuous Galerkin Spectral Element Method (DGSEM). 
The DGSEM allows for arbitrary order of accuracy and can represent complex geometries through the use of
 unstructured meshes with curvilinear elements. 
Here, we develop an entropy--stable discretization, which enhances its robustness and enables the industrialization of the method. 

This work is the third on a roadmap to obtain an entropy stable multiphase solver.
In \cite{2019:Manzanero-CH}, we developed a free--energy stable discretization of the 
Cahn--Hilliard equation, while in \cite{2019:Manzanero-iNS} we derived an 
entropy--stable approximation of the incompressible Navier--Stokes equations with 
artificial compressibility. We combine the developments in those previous papers to construct an 
entropy--stable DG approximation of the incompressible Navier--Stokes/Cahn--Hilliard 
system. 

As in \cite{2019:Manzanero-CH,2019:Manzanero-iNS}, the method is a nodal DG approximation with Gauss--Lobatto points, which satisfies the summation--by--parts simultaneous--aproximation--term (SBP--SAT) property \cite{2013:Fisher,2014:Carpenter}. The SBP--SAT property satisfied by the Gauss--Lobatto variant of the DGSEM allows us to mimic the continuous entropy analysis semi--discretely (i.e. discrete in space, continuous in time).
The SBP-SAT property originated with finite difference methods 
\cite{2013:Fisher,2014:Carpenter,2019:Chan},
and has been exploited in recent years to obtain entropy--stable DG schemes 
for the linear advection equation \cite{Kopriva2,2017:Manzanero}, 
Burgers equation \cite{2013:Gassner,2017:Gassner}, shallow water equations \cite{2016:Gassner-shallow}, Euler and Navier--Stokes equations \cite{2016:Gassner,2017:Gassner,2019:Manzanero-iNS}, the magneto--hydrodynamics equations \cite{2016:Winters}, 
and the Cahn--Hilliard equation \cite{2019:Manzanero-CH} among others. Although we develop a DG approximation, the approach taken here can be applied to any discretization that satisfies the SBP--SAT property (e.g. finite difference).

An entropy--stable approximation requires that the continuous system of partial differential equations satisfies an 
entropy law. When the density is allowed to 
vary between the two fluids in a multiphase flow, as opposed to the more traditional Model H \cite{1977:Hohenberg} with constant density, the design of an entropy--stable scheme requires care.
The variable density compromises the formulation used for the momentum equation \cite{2010:Shen}. Two main approaches to the design of a variable density entropy--stable approximation are commonly adopted: to use a skew--symmetric version of the momentum equation \cite{2000:Guermond}, 
or to augment the momentum equation with an additional diffusive flux that depends on the density difference  \cite{2007:Ding,2018:Dong}. 

The addition of the artificial compressibility model, where a divergence--free velocity is not strictly enforced, completes the formulation used 
for the body interfacial forces.
Without proper choices, terms proportional to $\svec{\nabla}\cdot\svec{u}$,
which can not be neglected nor bounded, might appear in the entopy equation.
With the proper choices, we confirm that the entropy of the Navier--Stokes/Cahn--Hilliard system obeys 
the Onsager principle; the total entropy is the sum of the entropy of the incompressible Navier--Stokes equations, plus the free--energy of the Cahn--Hilliard equation \cite{1958:Cahn}. 
In \cite{2019:Manzanero-iNS} we derived a mathematical entropy of the incompressible 
Navier--Stokes equations with artificial compressibility as the sum of the traditional kinetic energy, 
plus an additional energy term that accounts for compressibility effects. As a result, the total entropy becomes the sum of the kinetic energy, the artificial compressibility energy, and the free--energy.

In this paper we construct a multiphase approximation that is appropriate to solve typical applications 
in the oil industry. An example is the transport of crude oil in pipes, where the flow of oil and 
water is solved under high Reynolds number conditions. 
Since one desires to keep the mesh and degrees of freedom as low as possible, the flow configuration 
is often under--resolved. Under these circumstances, a robust method that is provably stable 
is attractive, as it avoids aliasing driven numerical instabilities that might lead to numerical 
divergence \cite{2017:Manzanero}.

The rest of this work is organized as follows: in Sec.~\ref{sec:Governing} we describe the incompressible Navier--Stokes/Cahn--Hilliard model, with the continuous entropy analysis in Sec.~\ref{subsec:Governing:iNS/CHEntropy}. The construction of the discrete DG approximation is described in Sec.\ref{sec:DGSEM}. Then, the semi--discrete stability analysis of the approximation is performed in Sec.~\ref{sec:Stability}, for which a summary can be found in Sec.~\ref{subsec:stability:Summary}. Lastly, we provide numerical experiments in Sec.~\ref{sec:NumericalExperiments} that assess both the accuracy and robustness of the method. We solve a manufactured solution convergence analysis in Sec.~\ref{subsec:num:conv}, and an assessment of the robustness solving a random initial condition in Sec.~\ref{subsec:num:Random}. Then, we solve classic static and rising bubble test problems in Sec.~\ref{subsec:num:StaticBubble} and Sec.~\ref{subsec:num:RisingBubble} respectively. Lastly, in Sec.~\ref{subsec:num:Annular} we challenge the technique by solving a multiphase pipe in the annular flow regime. Conclusions and discussions can be found in Sec.~\ref{sec:Conclusions}.

\section{Governing equations. Continuous entropy analysis}\label{sec:Governing}

In this section we describe a model that combines the incompressible 
Navier--Stokes and the Cahn--Hilliard equations, which we will refer to as iNS/CH. In the single phase variable density
incompressible Navier--Stokes equations, the density $\rho(\svec{x},t)$ is an independent variable that satisfies 
the continuity equation (see \cite{2000:Guermond,2017:Bassi,2019:Manzanero-iNS}),
\begin{equation}
  \rho_t + \svec{\nabla}\cdot\left(\rho\svec{u}\right) = 0,
  \label{eq:governing:continuity-iNS}
\end{equation}
where $\svec{u}(\svec{x},t)=\left(u,v,w\right)$ is the velocity 
field. In the two--phase iNS/CH system, however, the density is computed from the concentration $c(\svec{x},t)$ of 
the two fluids \cite{1998:Lowengrub},
\begin{equation}
  \rho = \rho(c) = \rho_1 c + \rho_2 (1-c),
\label{eq:governing:density-compatibility}
\end{equation}
where  $\rho_{1,2}$ are the 
densities of fluids 1 and 2, respectively, which are constant in space and time. 

In phase field methods, an advection--diffusion equation
drives the concentration. Among the 
model choices, we pick the Cahn--Hilliard equation \cite{1958:Cahn,1959:Cahn},
\begin{equation}
c_t + \svec{\nabla}\cdot\left(c\svec{u}\right) = M_0 \svec{\nabla}^2 \mu.
\label{eq:governing:cahn--hilliard}
\end{equation}
In \eqref{eq:governing:cahn--hilliard}, $\mu$ is the chemical potential,
\begin{equation}
  \mu = \frac{\diff f_0(c)}{\diff c} - \frac{3}{2}\sigma\varepsilon 
  \svec{\nabla}^2 c,
  \label{eq:governing:chemical-potential}
\end{equation}
where $f_0(c)$ is the chemical free--energy,
\begin{equation}
  f_0 = \frac{12\sigma}{\varepsilon}c^2(1-c)^2,
  \label{eq:governing:chemical-free-energy}
\end{equation}
$\sigma$ is the coefficient of interface tension between the fluids, $\varepsilon$ is the 
interface width, and $M_0$ is the mobility, computed in this work with the chemical characteristic time, $t_{\mathrm{CH}}$, as
\begin{equation}
  M_0 = \frac{\varepsilon}{\sigma t_{\mathrm {CH}}}.
\end{equation}
The parameters $\sigma$, $\varepsilon$, $M_0$, and $t_{\mathrm{CH}}$ 
are positive constants. The Cahn--Hilliard equation \eqref{eq:governing:cahn--hilliard}, with the chemical potential definition \eqref{eq:governing:chemical-potential}, has an associated free--energy
\begin{equation}
\mathcal F(c,\svec{\nabla}c) = f_0(c) + \frac{3}{4}\sigma\varepsilon|\svec{\nabla}c|^{2}.
\label{eq:governing:freeenergy-def}
\end{equation}

The velocity field is computed from the momentum equation. If one considers that the continuity 
equation \eqref{eq:governing:continuity-iNS} holds, the 
conservative,
\begin{equation}
  \left(\rho \svec{u}\right)_t +\svec{\nabla}\cdot\left(\rho \svec{u}\svec{u}\right) 
  = -\svec{\nabla}p_{s} -\svec{\nabla}\cdot\left(\svec{\nabla}c\otimes\frac{\partial \mathcal F}{\partial \svec{\nabla}c}\right)+ \svec{\nabla}\cdot\left(\eta\left(\svec{\nabla}\svec{u} + 
  \svec{\nabla}\svec{u}^{T}\right)\right)+\rho\svec{g},
  \label{eq:governing:momentum-conservative}
\end{equation}
non--conservative,
\begin{equation}
\rho  \left(\svec{u}\right)_t +\rho\svec{u}\cdot\svec{\nabla}\svec{u}
   = -\svec{\nabla}p_{s} -\svec{\nabla}\cdot\left(\svec{\nabla}c\otimes\frac{\partial \mathcal F}{\partial \svec{\nabla}c}\right)+ \svec{\nabla}\cdot\left(\eta\left(\svec{\nabla}\svec{u} + 
  \svec{\nabla}\svec{u}^{T}\right)\right)+\rho\svec{g},
  \label{eq:governing:momentum-nonconservative}
\end{equation}
and skew--symmetric
\begin{equation}
 \frac{1}{2}\left(\rho \svec{u}\right)_t+\frac{1}{2}\rho  \left(\svec{u}\right)_t  +\svec{\nabla}\cdot\left(\frac{1}{2}\rho \svec{u}\svec{u}\right) 
 +\frac{1}{2}\rho\svec{u}\cdot\svec{\nabla}\svec{u}
   = -\svec{\nabla}p_{s} -\svec{\nabla}\cdot\left(\svec{\nabla}c\otimes\frac{\partial \mathcal F}{\partial \svec{\nabla}c}\right)+ \svec{\nabla}\cdot\left(\eta\left(\svec{\nabla}\svec{u} + 
  \svec{\nabla}\svec{u}^{T}\right)\right)+\rho\svec{g},
  \label{eq:governing:momentum-skewsymmetric}
\end{equation}
forms of the momentum equation are identical in the continuous setting. However, since the density does not satisfy the continuity 
equation \eqref{eq:governing:continuity-iNS}, but rather the compatibility condition \eqref{eq:governing:density-compatibility}, the three 
forms \eqref{eq:governing:momentum-conservative}, \eqref{eq:governing:momentum-nonconservative} 
and \eqref{eq:governing:momentum-skewsymmetric} are no longer equivalent. In this situation, the momentum equation is \textit{frozen} as one of the three forms. Following \cite{2000:Guermond,2010:Shen}, 
 we choose the skew--symmetric version of the momentum \eqref{eq:governing:momentum-skewsymmetric}
 since is the only one that satisfies an entropy inequality for any positive density field. (An alternative that produces an entropy--stable scheme for the iNS/CH system is to include a relative flux in the momentum equation that models the diffusion of the components \cite{2012:Abels}.) We use the chain rule in time to perform an additional
manipulation on the time derivative of \eqref{eq:governing:momentum-skewsymmetric} and obtain the momentum form by Guermond et 
al. \cite{2000:Guermond},
\begin{equation}
\sqrt{\rho}\left(\sqrt{\rho}\svec{u}\right)_t+\svec{\nabla}\cdot\left(\frac{1}{2}\rho \svec{u}\svec{u}\right) 
 +\frac{1}{2}\rho\svec{u}\cdot\svec{\nabla}\svec{u}
   = -\svec{\nabla}p_s -\svec{\nabla}\cdot\left(\svec{\nabla}c\otimes\frac{\partial \mathcal F}{\partial \svec{\nabla}c}\right)+ \svec{\nabla}\cdot\left(\eta\left(\svec{\nabla}\svec{u} + 
  \svec{\nabla}\svec{u}^{T}\right)\right)+\rho\svec{g}.
  \label{eq:governing:momentum-skewsymmetric-sqrtRho-oldCapillary}
\end{equation}
In \eqref{eq:governing:momentum-skewsymmetric-sqrtRho-oldCapillary}, $p_s\left(\svec{x},t\right)$ is the static pressure. The second term is the body force approximation of the capillary 
pressure \cite{1998:Lowengrub},
\begin{equation}
\svec{p}_{c} = -\svec{\nabla}\cdot\left(\svec{\nabla}c\otimes\frac{\partial \mathcal F}{\partial \svec{\nabla}c}\right)=-\frac{3}{2}\sigma\varepsilon\svec{\nabla}\cdot\left(\svec{\nabla}c\otimes\svec{\nabla}c\right).
\label{eq:governing:capillary-pressure}
\end{equation}
Finally $\eta$ is the viscosity, computed from the concentration $c$ and the 
viscosities of fluids 1 and 2, $\eta_{1,2}$,
\begin{equation}
  \eta = \eta(c) = \eta_1 c + \eta_2 \left(1-c\right).
  \label{eq:governing:compatibility-viscosity}
\end{equation}

As in \cite{1998:Lowengrub}, we write the capillary pressure more conveniently. Applying the derivative product rule,
\begin{equation}
\svec{p}_{c} = -\frac{3}{2}\sigma\varepsilon\left(\svec{\nabla}^{2}c\svec{\nabla}c + \svec{\nabla}\left(\frac{1}{2}|\svec{\nabla}c|^2\right)\right),
\end{equation}
and adding and subtracting the term $\svec{\nabla} f_0=f'_0(c)\svec{\nabla}c$, we can rewrite the expression for the body force as
\begin{equation}
\svec{p}_{c}=\left(f'_0(c)-\frac{3}{2}\sigma\varepsilon\svec{\nabla}^{2}c\right)\svec{\nabla}c-\svec{\nabla}\left(f_0(c) + \frac{3}{4}\sigma\varepsilon|\svec{\nabla}c|^2\right)= \mu\svec{\nabla}c-\svec{\nabla}\mathcal F = -c\svec{\nabla}\mu-\svec{\nabla}\left(\mathcal F-\mu c\right).
\label{eq:governing:capillary-pressure-equivalency}
\end{equation}
In the last form, an additional application of the product derivative rule for the first term has been performed. Eq.~\eqref{eq:governing:capillary-pressure-equivalency} represents equivalent expressions of the capillary pressure \eqref{eq:governing:capillary-pressure}, the last form being the one adopted here. The first term of the last form, $-c\svec{\nabla}\mu$, is a non--conservative term, and the second, $-\svec{\nabla}\left(\mathcal F-\mu c\right)$, is combined with the static pressure gradient to define an auxiliary pressure
\begin{equation}
p = p_s + \mathcal F-\mu c.
\label{eq:governing:pressure}
\end{equation}

With all of these manipulations, we get the final expression of the momentum equation,
\begin{equation}
\sqrt{\rho}\left(\sqrt{\rho}\svec{u}\right)_t+\svec{\nabla}\cdot\left(\frac{1}{2}\rho \svec{u}\svec{u}\right) 
+\frac{1}{2}\rho\svec{u}\cdot\svec{\nabla}\svec{u}+c\svec{\nabla}\mu
= -\svec{\nabla}p + \svec{\nabla}\cdot\left(\eta\left(\svec{\nabla}\svec{u} + 
\svec{\nabla}\svec{u}^{T}\right)\right)+\rho\svec{g}.
\label{eq:governing:momentum-skewsymmetric-sqrtRho}
\end{equation}

Many authors have explored how to relax the incompressibility constraint  $\svec{\nabla}\cdot\svec{u}=0$ in Navier--Stokes and Navier--Stokes/Cahn--Hilliard systems by using the artificial (or \textit{pseudo})
compressibility method \cite{1996:Shen,1997:Shen,2018:Feng,2019:Zhu}. We use the time derivative 
of the pressure \cite{1997:Shen,2017:Bassi,2019:Manzanero-iNS}, 
\begin{equation}
  p_t + \rho_0 c_0^2 \svec{\nabla}\cdot\svec{u} = 0,
  \label{eq:governing:ACM}
\end{equation}
where $c_0$ is the artificial sound speed, and  $\rho_0 = 
\max\left(\rho_1,\rho_2\right)$.%
\begin{rem}
Alternatively, one can take the time derivative of the pressure Laplacian \cite{1996:Shen,2010:Shen,2015:Shen},
\begin{equation}
  -\svec{\nabla}^{2}p_t + \rho_0 c_0^2 \svec{\nabla}\cdot\svec{u} = 0.
  \label{eq:governing:ACM-2}  
\end{equation}
Since both models share the same velocity 
divergence term, we can modify the proofs of stability for the approximation of \eqref{eq:governing:ACM} and study the approximation of the second model 
\eqref{eq:governing:ACM-2}. We do so in Appendix~\ref{app:second-ACM}. The use of the Laplacian of the pressure leads to a more complicated implicit implementation, so  for that scheme we include only the stability proof in this paper.
\end{rem}

The combination of \eqref{eq:governing:cahn--hilliard}, \eqref{eq:governing:momentum-skewsymmetric-sqrtRho} 
and \eqref{eq:governing:ACM} is the iNS/CH system studied in this paper. The coupling between the 
Cahn--Hilliard equation and the incompressible Navier--Stokes equations is 
\textit{two--way}, through the compatibility conditions 
\eqref{eq:governing:density-compatibility}, \eqref{eq:governing:compatibility-viscosity}, 
and the capillary pressure. We note that even in the continuous setting, the momentum is not 
conserved because of the choice of the skew--symmetric form \eqref{eq:governing:momentum-skewsymmetric-sqrtRho}. 
Although one can argue about the physical implications, we note that the capillary force is a non--conservative 
term as well, so momentum conservation is not guaranteed in any of the forms 
\eqref{eq:governing:momentum-conservative}, 
\eqref{eq:governing:momentum-nonconservative}, and 
\eqref{eq:governing:momentum-skewsymmetric-sqrtRho}. Alternative incompressible Navier--Stokes/Cahn--Hilliard discretizations using the skew--symmetric form of the momentum equation can be 
found in \cite{2010:Boyer,2010:Shen}. 

To simplify the notation, we write the iNS/CH system as a general advection--diffusion equation,
\begin{equation}
\begin{split}
&\left(\begin{array}{ccccc}1 & 0 & 0 & 0 & 0 \\
0 & \sqrt{\rho} & 0 & 0 & 0 \\
0 & 0 & \sqrt{\rho} & 0 & 0 \\
0 & 0 & 0 & \sqrt{\rho} & 0 \\
0 & 0 & 0 & 0 & 1\end{array}\right)\left(\begin{array}{c}c \\ \sqrt{\rho} u \\
\sqrt{\rho}v \\ \sqrt{\rho} w \\ p \end{array}\right)_{t} + \svec{\nabla}\cdot\left(\begin{array}{c}c\svec{u} \\
\frac{1}{2}\rho \svec{u}u + p\svec{e}_{1} \\
\frac{1}{2}\rho \svec{u}v + p\svec{e}_{2} \\
\frac{1}{2}\rho \svec{u}w + p\svec{e}_{3} \\
0
\end{array}\right) \\
&+ \left(\begin{array}{c}\svec{0} \\
\frac{1}{2}\rho \svec{u}\cdot\svec{\nabla}u+c\svec{e}_{1}\cdot\svec{\nabla}\mu\\
\frac{1}{2}\rho \svec{u}\cdot\svec{\nabla}v+c\svec{e}_{2}\cdot\svec{\nabla}\mu\\
\frac{1}{2}\rho \svec{u}\cdot\svec{\nabla}w+c\svec{e}_{3}\cdot\svec{\nabla}\mu\\
\rho_0c_0^2\left(\svec{e}_{1}\cdot\svec{\nabla}u+\svec{e}_{2}\cdot\svec{\nabla}v+\svec{e}_{3}\cdot\svec{\nabla}w\right)
\end{array}\right)
=\svec{\nabla}\cdot\left(\begin{array}{c}M_0\svec{\nabla}\mu\\
2\eta\tens{S}\cdot\svec{e}_{1} \\
2\eta\tens{S}\cdot\svec{e}_{2} \\
2\eta\tens{S}\cdot\svec{e}_{3} \\
\svec{0} 
\end{array}\right)+\left(\begin{array}{c} 
0 \\
\rho \svec{g}\cdot\svec{e}_{1} \\
\rho \svec{g}\cdot\svec{e}_{2} \\
\rho \svec{g}\cdot\svec{e}_{3} \\
0
\end{array}\right),
\end{split}
\label{eq:governing:iNS/CH-system}
\end{equation}
where,
\begin{equation}
\tens{S}=\mathrm{sym}\left(\svec{\nabla}\svec{u}\right)=\frac{1}{2}\left(\svec{\nabla}\svec{u}+\svec{\nabla}\svec{u}^{T}\right),
\label{eq:governing:strain-tensor}
\end{equation}
is the strain tensor, and $\svec{e}_{i}$ are the space unit vectors. Note that we have grouped the velocity divergence in the artificial compressibility 
equation \eqref{eq:governing:ACM} into the non--conservative terms. 
Although treating the velocity divergence as conservative or non--conservative 
is equivalent, the latter makes it easier to show stability.

We adopt the notation in \cite{2017:Gassner} to work with vectors of different nature.  We define space vectors (e.g. $\svec{x}=\left(x,y,z\right)\in\R^{3}$) with
an arrow on top, and state vectors (e.g. $\stvec{q}=\left(c,\sqrt{\rho}\svec{u},p\right)\in\R^{5}$) in bold. 
Moreover, we define block vectors as the result of stacking three state vectors (e.g. fluxes),
\begin{equation}
\ssvec{f}_{e} = \left(\begin{array}{ccc}\stvec{f}_{e,1} \\ \stvec{f}_{e,2} \\ \stvec{f}_{e,3}\end{array}\right)= \left(\begin{array}{ccc}\stvec{f}_e \\ \stvec{g}_e \\ \stvec{h}_e\end{array}\right),~~
\ssvec{f}_{v} = \left(\begin{array}{ccc}\stvec{f}_{v,1} \\ \stvec{f}_{v,2} \\ \stvec{f}_{v,3}\end{array}\right)= \left(\begin{array}{ccc}\stvec{f}_v \\ \stvec{g}_v \\ 
\stvec{h}_v\end{array}\right),
\end{equation}
and define the operator $\Upsilon$ that transforms a $5\times 3$ (state--space) 
matrix into a $15\times 1$ block vector,
\begin{equation}
  \ssvec{f} = \Upsilon\left(\begin{array}{ccc}f_1 & g_1 & h_1 \\
  f_2 & g_2 & h_2 \\
  f_3& g_3& h_3\\
  f_4& g_4& h_4\\
  f_5& g_5& h_5\end{array}\right) = \left(\begin{array}{c}\stvec{f} \\ \stvec{g} \\ 
  \stvec{h}\end{array}\right).
\end{equation}
This notation allows us to define products of state, space, and block vectors,
\begin{equation}
\ssvec{f}\cdot\ssvec{g} = \sum_{i=1}^{3}\stvec{f}_{i}^T\stvec{g}_{i},~~\svec{g}\cdot\ssvec{f} = \sum_{i=1}^{3}g_i 
\stvec{f}_{i},~~\svec{g}\stvec{f} = \left(\begin{array}{c}g_1\stvec{f} \\ g_2\stvec{f} \\ 
g_3\stvec{f}\end{array}\right).
\label{eq:governing:space-state-block-ops}
\end{equation}
We can then define the divergence and gradient operators using \eqref{eq:governing:space-state-block-ops},
\begin{equation}
\svec{\nabla}\cdot\ssvec{f} = \sum_{i=1}^{3}\frac{\partial\stvec{f}_{i}}{\partial x_i},~~\svec{\nabla}\stvec{q} = \left(\begin{array}{c}\stvec{q}_{x} \\ \stvec{q}_{y} \\ 
\stvec{q}_z\end{array}\right).
\end{equation}
Lastly, we refer to state matrices (i.e. $5\times 5$ matrices) with an underline,  e.g. $\smat{B}$, which can be combined to construct a block matrix,
\begin{equation}
\mathcal B = \left(\begin{array}{ccc}\smat{B}_{11} & \smat{B}_{12} & \smat{B}_{13} \\
\smat{B}_{21} & \smat{B}_{22} & \smat{B}_{23}\\
\smat{B}_{31} & \smat{B}_{32} & \smat{B}_{33}\end{array}\right).
\label{eq:notation:block-matrix}
\end{equation}
Block matrices can be directly multiplied to a block vector to obtain another block vector. For instance, to perform a matrix multiplication in space (e.g. a rotation),
\begin{equation}
\svec{g} = \tens{M}\svec{f},
\end{equation}
for each of the variables in the state vector, we construct the block matrix version of $\tens{M}$ ($\smat{I}_{5}$ is the $5\times 5$ identity matrix),
\begin{equation}
\mathcal M = \left(\begin{array}{ccc}M_{11}\smat{I}_{5} & M_{12}\smat{I}_{5} & M_{13}\smat{I}_{5}\\
M_{21}\smat{I}_{5} & M_{22}\smat{I}_{5} & M_{23}\smat{I}_{5}\\
M_{31}\smat{I}_{5} & M_{32}\smat{I}_{5} & M_{33}\smat{I}_{5}\end{array}\right),
\label{eq:notation:block-matrix-from-space-matrix}
\end{equation}
so that we can compactly write
\begin{equation}
\ssvec{g} = \mathcal M \ssvec{f}.
\end{equation}
For more details, see \cite{2017:Gassner}.

The notation introduced above makes it possible to write the iNS/CH system \eqref{eq:governing:iNS/CH-system} compactly
in the form of a general advection--diffusion equation,
\begin{equation}
  \smat{m}\stvec{q}_{t} + \svec{\nabla}\cdot\ssvec{f}_{e}\left(\stvec{q}\right) + 
  \sum_{m=1}^{5}\ssvecg{\phi}_{m}\left(\stvec{q}\right)\cdot\svec{\nabla}w_{m} = 
  \svec{\nabla}\cdot\ssvec{f}_{v}\left(\stvec{q},\svec{\nabla}\stvec{w}\right) + \stvec{s}\left(\stvec{q}\right),
  \label{eq:governing:general-advdiff}
\end{equation}
with state vector $\stvec{q} = \left(c,\sqrt{\rho}\svec{u},p\right)$, gradient variables vector $\stvec{w}=\left(w_{1},w_{2},...,w_{5}\right)=\left(\mu,\svec{u},p\right)$, mass 
matrix $\smat{m}$,
\begin{equation}
  \smat{m} = \left(\begin{array}{ccc}1 & 0 & 0 \\
  0& \sqrt{\rho}\tens{I}_{3} & 0 \\ 0 & 0 & 1\end{array}\right),
\end{equation}
inviscid fluxes $\ssvec{f}_{e}\left(\stvec{q}\right)$,
\begin{equation}
\stvec{f}_{e,1} = \stvec{f}_{e} = \left(\begin{array}{c}
  c u \\
  \rho u^2 + p \\
  \rho u v \\
  \rho u w \\
  0 
  \end{array}\right),~~
  \stvec{f}_{e,2} = \stvec{g}_{e} = \left(\begin{array}{c}
  c v \\
  \rho uv \\
  \rho v^2 + p \\
  \rho v w \\
  0 
  \end{array}\right),~~
  \stvec{f}_{e,3} = \stvec{h}_{e} = \left(\begin{array}{c}
  c w \\
  \rho uw \\
  \rho vw \\
  \rho w^2+p \\
  0 
  \end{array}\right),
\end{equation}
non--conservative term coefficients $\ssvecg{\phi}_{m}\left(\stvec{q}\right)$,
\begin{equation}
\begin{split}
&  \ssvecg{\phi}_{1} = \Upsilon\left(\begin{array}{c} 0 \\ c \svec{e}_{1} \\
c\svec{e}_{2} \\ c\svec{e}_{3} \\ 0 
  \end{array}\right),~~\ssvecg{\phi}_{2} = \Upsilon\left(\begin{array}{c}
  0 \\
  \frac{1}{2}\rho\svec{u} \\
  0 \\
  0 \\
  \rho_0c_0^2 \svec{e}_{1}
\end{array}  \right),~~\ssvecg{\phi}_{3} = \Upsilon\left(\begin{array}{c}
  0 \\
  0 \\
  \frac{1}{2}\rho\svec{u} \\
  0 \\
  \rho_0c_0^2 \svec{e}_{2}
\end{array}  \right), \\
&\ssvecg{\phi}_{4} = \Upsilon\left(\begin{array}{c}
  0 \\
  0 \\
  0 \\
  \frac{1}{2}\rho\svec{u} \\
  \rho_0c_0^2 \svec{e}_{3}
\end{array}  \right),~~\ssvecg{\phi}_{5} = \Upsilon\left(\begin{array}{c}0 \\ 0 \\ 0 \\ 0 \\ 0 
\end{array}\right),
\end{split}
\end{equation}
viscous fluxes $\ssvec{f}_{v}\left(\svec{\nabla}\stvec{w}\right)$,
\begin{equation}
  \stvec{f}_{v,1} = \stvec{f}_{v} = \left(\begin{array}{c}M_0 \mu_{x}\\
  2\eta \tens{S}_{11}\\
  2\eta \tens{S}_{21}\\
  2\eta \tens{S}_{31}\\
  0 \end{array}\right),~~  \stvec{f}_{v,2} = \stvec{g}_{v} = \left(\begin{array}{c}M_0 \mu_{y}\\
  2\eta \tens{S}_{12}\\
  2\eta \tens{S}_{22}\\
  2\eta \tens{S}_{32}\\
  0 \end{array}\right),~~  \stvec{f}_{v,3} = \stvec{h}_{v} = \left(\begin{array}{c}M_0 \mu_{z}\\
  2\eta \tens{S}_{13}\\
  2\eta \tens{S}_{23}\\
  2\eta \tens{S}_{33}\\
  0 \end{array}\right),
\end{equation}
and source term $\stvec{s}\left(\stvec{q}\right)=\left(0,\rho\svec{g},0\right)$.

An attractive property of the iNS/CH system from the point of view of computational efficiency is that the same gradient 
variables $\stvec{w}$ are used in both 
the non--conservative terms and the viscous fluxes. Moreover, these gradient variables will be 
shown to be the entropy variables associated with the 
mathematical entropy  in Section \ref{subsec:Governing:iNS/CHEntropy}. Additionally, recall that the velocity 
divergence from the artificial compressibility equation \eqref{eq:governing:ACM} was grouped 
into the non--conservative terms since it is beneficial when proving stability.

\subsection{Entropy analysis of the iNS/CH system}\label{subsec:Governing:iNS/CHEntropy}

The entropy analysis rests on the existence of a pair $\mathcal E\left(\stvec{q}\right)$ 
(mathematical entropy) and $\stvec{w}$ (entropy variables) that 
contract the system of equations \eqref{eq:governing:general-advdiff} into a 
conservation law \cite{2003:Tadmor},
\begin{equation}
\begin{split}
  \stvec{w}^{T}&\left(\smat{m}\stvec{q}_{t} + \svec{\nabla}_{}\cdot\ssvec{f}_{e}\left(\stvec{q}\right) + 
  \sum_{m=1}^{5}\ssvecg{\phi}_{m}\left(\stvec{q}\right)\cdot\svec{\nabla}w_{m} - 
  \svec{\nabla}\cdot\ssvec{f}_{v}\left(\stvec{q},\svec{\nabla}\stvec{w}\right)\right) 
  \\
&  = \mathcal E_{t} + \stvec{\nabla}\cdot\svec{f}^{\mathcal E} + 
  \stvec{\nabla}\stvec{w}^{T}\cdot\ssvec{f}_{v} = 0,
  \end{split}
  \label{eq:continuous:entropy-pde}
\end{equation}
with an entropy flux $\svec{f}^{\mathcal E}$, and a viscous dissipative contribution,
\begin{equation}
  \svec{\nabla}\stvec{w}^{T}\cdot\ssvec{f}_{v} \geqslant 0. 
  \label{eq:continuous:viscous-def-positive}
\end{equation}

The difference between the usual entropy analysis \cite{2003:Tadmor,2014:Carpenter,2014:Carpenter} and this work 
is that we will incorporate the Cahn--Hilliard free--energy into the 
mathematical entropy, which makes it depend not only on the solution, but also its 
gradient. To obtain the entropy of the iNS/CH system, we 
add the Cahn--Hilliard free--energy $\mathcal F$ to the incompressible NSE entropy defined in 
\cite{2019:Manzanero-iNS}. The latter combines the traditional kinetic energy $\mathcal K$ 
with an extra energy term due to the artificial compressibility $\mathcal 
E_{\mathrm{AC}}$. For the iNS/CH system, the total entropy is,
\begin{equation}
\begin{split}
  \mathcal E &= \mathcal F(c,\svec{\nabla} c) + \mathcal K\left(\sqrt{\rho}\svec{u}\right) + \mathcal E_{\mathrm{AC}}(p) 
  \\
  &= f_0(c) + \frac{3}{4}\sigma\varepsilon|\svec{\nabla}c|^2 + \frac{1}{2}\rho 
  v_{tot}^2 + \frac{p^2}{2\rho_0 c_0^2}\geqslant 0 ,
  \end{split}
  \label{eq:continuous:entropy-def}
\end{equation}
where $v_{tot}^2 = u^2 + v^2 + w^2$ is the square of the total speed. Eq.~\eqref{eq:continuous:entropy-def} 
assumes positivity on the density $\rho$, which is ensured by limiting the 
maximum and minimum density values used in the momentum equation with a simple cutoff,
\begin{equation}
  \rho\left(\hat{c}\right) = \rho_1 \hat{c} + 
  \rho_{2}\left(1-\hat{c}\right),~~\hat{c}=\min\left(\max\left(c,0\right),1\right).
\end{equation}

We find a set of entropy variables $\stvec{w}$ that contracts the time derivative of the state vector, $\stvec{q}_{t}$, into the time derivative of the entropy $\mathcal E_{t}$, plus an additional divergence term of a time derivative flux $\svec{f}_{t}^{\mathcal E}$,
\begin{equation}
\stvec{w}^{T}\smat{m}\stvec{q}_{t} = \mathcal E_{t} + \svec{\nabla}\cdot\svec{f}_{t}^{\mathcal E}.
\label{eq:continuous:entropy-time-contraction}
\end{equation}
To do so, we compute the time derivative of the entropy \eqref{eq:continuous:entropy-def},
\begin{equation}
\begin{split}
\mathcal E_{t} &= f_0'(c)c_{t} + \frac{3}{2}\sigma\varepsilon \nabla c\cdot\nabla c_{t} + \sqrt{\rho}\svec{u}\cdot\left(\sqrt{\rho}\svec{u}\right)_{t}+\frac{pp_{t}}{\rho_0 c_0^2} \\
&=\left(f_0'(c) - \frac{3}{2}\sigma\varepsilon\svec{\nabla}^{2}c\right)c_{t}+\svec{\nabla}\cdot\left(\frac{3}{2}\sigma\varepsilon c_{t}\svec{\nabla}c \right) + \sqrt{\rho}\svec{u}\cdot\left(\sqrt{\rho}\svec{u}\right)_{t}+\frac{pp_{t}}{\rho_0 c_0^2} \\
&=\mu c_{t} + \svec{\nabla}\cdot\left(\frac{3}{2}\sigma\varepsilon c_{t}\svec{\nabla}c \right) + \sqrt{\rho}u\left(\sqrt{\rho} u\right)_{t}+ \sqrt{\rho}v\left(\sqrt{\rho} v\right)_{t} + \sqrt{\rho}w\left(\sqrt{\rho} w\right)_{t}+\frac{pp_{t}}{\rho_0 c_0^2},
\end{split}
\label{eq:continuous:entropy-time-contraction:entropy-time}
\end{equation}
and the product of the entropy variables with the time derivative of the state vector,
\begin{equation}
w_{1}c_{t} + \sqrt{\rho}w_{2}\left(\sqrt{\rho}u\right)_{t}+ \sqrt{\rho}w_{3}\left(\sqrt{\rho}v\right)_{t}+ \sqrt{\rho}w_{4}\left(\sqrt{\rho}w\right)_{t}+w_{5}p_{t}.
\label{eq:continuous:entropy-time-contraction:state-time}
\end{equation}
Thus, we replace \eqref{eq:continuous:entropy-time-contraction:entropy-time} and \eqref{eq:continuous:entropy-time-contraction:state-time} in \eqref{eq:continuous:entropy-time-contraction} and rearrange,
\begin{equation}
\begin{split}
\left(w_{1}-\mu\right)c_{t} &+ \sqrt{\rho}\left(w_{2}-u\right)\left(\sqrt{\rho}u\right)_{t}+ \sqrt{\rho}\left(w_{3}-v\right)\left(\sqrt{\rho}v\right)_{t}+ \sqrt{\rho}\left(w_{4}-w\right)\left(\sqrt{\rho}w\right)_{t}\\
&+\left(w_{5}-\frac{p}{\rho_0c_0^2}\right)p_{t} = \svec{\nabla}\cdot\left(\svec{f}_{t}^{\mathcal E} + \frac{3}{2}\sigma\varepsilon c_{t}\svec{\nabla}c \right),
\end{split}
\end{equation}
to extract the entropy variables $\stvec{w}$ and time entropy flux $\svec{f}_{t}^{\mathcal E}$
\begin{equation}
\stvec{w}=\left(\mu,u,v,w,\frac{p}{\rho_0 c_0^2}\right),~~ \svec{f}_{t}^{\mathcal E}=-\frac{3}{2}\sigma\varepsilon c_{t}\svec{\nabla}c.
\end{equation}
Note that the minimum of the entropy is found when all the entropy variables are zero.

The contraction of the inviscid fluxes rests on one property to be satisfied by 
the inviscid fluxes and the non--conservative terms. 

\begin{prop}
If inviscid fluxes, $\ssvec{f}_{e}$, and non--conservative terms coefficients, $\ssvecg{\phi}_{m}$, satisfy
\begin{equation}
  \stvec{e}_{m}^{T}\ssvec{f}_{e} = \stvec{w}^{T}\ssvec{\phi}_{m},
  \label{eq:continuous:entropy-condition}
\end{equation}
where $\stvec{e}_{m}$ is the state unit vector along the $m$--th state variable, then the entropy variables automatically contract the inviscid fluxes, and the 
entropy flux is computed as
\begin{equation}
\svec{f}_{e}^{\mathcal E} = \stvec{w}^{T}\ssvec{f}_{e}.
\label{eq:continuous:entropy-flux}
\end{equation}
\end{prop}
To show the contraction property \eqref{eq:continuous:entropy-condition}, we multiply the divergence of the inviscid fluxes by the entropy variables, apply the derivative product rule, 
and write the second scalar product as the sum of the product of the five state 
components. Lastly we use \eqref{eq:continuous:entropy-condition} and \eqref{eq:continuous:entropy-flux} to see that
\begin{equation}
\begin{split}
  \stvec{w}^{T}\svec{\nabla}\cdot\ssvec{f}_{e} &=  
  \svec{\nabla}\cdot\left(\stvec{w}^{T}\ssvec{f}_{e}\right)-\nabla\stvec{w}^{T}\cdot\ssvec{f}_{e} 
  = \svec{\nabla}\cdot\left(\stvec{w}^{T}\ssvec{f}_{e}\right)- 
  \sum_{m=1}^{5}\stvec{e}^{T}_{m}\ssvec{f}_{e}\cdot\svec{\nabla}w_{m}\\
  &= \svec{\nabla}\cdot\svec{f}_{e}^{\mathcal E} - 
  \sum_{m=1}^{5}\stvec{w}^{T}\ssvec{\phi}_{m}\cdot\svec{\nabla}w_{m}.
  \end{split}
\end{equation}
Therefore,
\begin{equation}
\stvec{w}^{T}\left(\svec{\nabla}_{}\cdot\ssvec{f}_{e}\left(\stvec{q}\right) + 
  \sum_{m=1}^{5}\ssvecg{\phi}_{m}\left(\stvec{q}\right)\cdot\svec{\nabla}w_{m}\right) 
  =\svec{\nabla}\cdot\svec{f}_{e}^{\mathcal E}.
  \label{eq:continuous:inviscid-fluxes-contraction}
\end{equation}

The iNS/CH system does satisfy the property 
\eqref{eq:continuous:entropy-condition}, for
\begin{equation}
\begin{split}
  \stvec{e}_{1}^{T}\ssvec{f}_{e} &= c\svec{u} = \stvec{w}^{T}\ssvecg{\phi}_{1},\\
    \stvec{e}_{2}^{T}\ssvec{f}_{e} &= \frac{1}{2}\rho u\svec{u}+p\svec{e}_{1} = 
    \stvec{w}^{T}\ssvecg{\phi}_{2},\\
    \stvec{e}_{3}^{T}\ssvec{f}_{e} &= \frac{1}{2}\rho v\svec{u}+p\svec{e}_{2} = 
    \stvec{w}^{T}\ssvecg{\phi}_{3},\\
    \stvec{e}_{4}^{T}\ssvec{f}_{e} &= \frac{1}{2}\rho w\svec{u}+p\svec{e}_{3} = 
    \stvec{w}^{T}\ssvecg{\phi}_{4},\\
        \stvec{e}_{5}^{T}\ssvec{f}_{e} &= 0 = 
    \stvec{w}^{T}\ssvecg{\phi}_{5}.\\
  \end{split}
  \label{eq:continuous:cancellation-condition}
\end{equation}
So, the inviscid entropy flux for the iNS/CH system follows from \eqref{eq:continuous:entropy-flux} and is
\begin{equation}
  \svec{f}_{e}^{\mathcal E} = \left(\frac{1}{2}\rho v_{tot}^{2}+p+\mu
  c\right)\svec{u}.
\end{equation}
For completeness, insight into how each of the different terms are combined into an entropy flux is 
provided in Appendix~\ref{app:inviscid-contraction}.

Lastly, the entropy variables contract viscous fluxes. The entropy flux follows \eqref{eq:continuous:entropy-flux},
\begin{equation}
\svec{f}_{v}^{\mathcal E} = \stvec{w}^{T}\ssvec{f}_{v} = M_0 \mu 
\svec{\nabla}\mu + 2\eta \tens{S}\cdot\svec{u}.
\end{equation}
Also, \eqref{eq:continuous:viscous-def-positive} holds,
\begin{equation}
\begin{split}
\svec{\nabla} \stvec{w}^{T}\cdot\ssvec{f}_{v} &= M_0 |\svec{\nabla}\mu|^2 + 
\eta\left(\svec{\nabla}\svec{u} + 
\svec{\nabla}\svec{u}^{T}\right):\svec{\nabla}\svec{u}\\
&=M_0|\svec{\nabla}\mu|^2 + 2\eta\frac{1}{2}\left(\svec{\nabla}\svec{u} + 
\svec{\nabla}\svec{u}^{T}\right):\frac{1}{2}\left(\svec{\nabla}\svec{u} + 
\svec{\nabla}\svec{u}^{T}\right) \\
&=M_0|\svec{\nabla}\mu|^2 + 2\eta \tens{S}:\tens{S}\geqslant 0.
\end{split}
\label{eq:continuous:viscous-volume-terms}
\end{equation}
Therefore,
\begin{equation}
\stvec{w}^{T}\svec{\nabla}\cdot\ssvec{f}_{v} = \svec{\nabla}\cdot 
\left(\stvec{w}^{T}\ssvec{f}_{v}\right) - 
\svec{\nabla}\stvec{w}^{T}\cdot\ssvec{f}_{v} \leqslant \svec{\nabla}\cdot\svec{f}^{\mathcal E}_{v}.
\end{equation}
In \eqref{eq:continuous:viscous-volume-terms}, we replaced $\svec{\nabla}\svec{u}$ by its symmetric part \eqref{eq:governing:strain-tensor} since it multiplies a symmetric tensor. 

With these results, we confirm that the iNS/CH system, \eqref{eq:governing:iNS/CH-system}, with entropy, \eqref{eq:continuous:entropy-def}, satisfies the conservation 
law
\begin{equation}
  \mathcal E_{t} + \svec{\nabla}\cdot\svec{f}^{\mathcal E}  = 
  -\svec{\nabla}\stvec{w}^{T}\cdot\ssvec{f}_{v} = -M_0|\svec{\nabla}\mu|^2 - 2\eta\tens{S}:\tens{S}\leqslant 0,
  \label{eq:continuous:entropy-pde-final}
\end{equation}
with entropy flux $\svec{f}^{\mathcal E}$,
\begin{equation}
  \svec{f}^{\mathcal E} = \svec{f}^{\mathcal E}_{t}+ \svec{f}^{\mathcal E}_{e}+ \svec{f}^{\mathcal E}_{v} 
  = -\frac{3}{2}\sigma\varepsilon c_{t}\svec{\nabla} c + \left(\frac{1}{2}\rho v_{tot}^{2}+p+\mu 
  c \right)\svec{u}+ 2\eta\tens{S}\cdot\svec{u} + M_0\mu\svec{\nabla}\mu.
\end{equation}

Eq.~\eqref{eq:continuous:entropy-pde-final} 
shows that the entropy \eqref{eq:continuous:entropy-def} is always dissipated in the interior of the domain, 
and it can only increase due to boundary exchanges. The entropy equation can be written as a global equation by integrating 
over the whole domain $\Omega$,
\begin{equation}
\frac{\diff \bar{\mathcal E}}{\diff t} + \int_{\partial \Omega}\svec{f}^{\mathcal E}\cdot\svec{n}\diff S = -\int_{\Omega}\left(M_0|\svec{\nabla}\mu|^2 + 2\eta\tens{S}:\tens{S}\right)\diff\svec{x} \leqslant 0
\label{eq:continuous:entropy-balance}
\end{equation}
where $\bar{\mathcal E}$ is the total entropy,
\begin{equation}
\bar{\mathcal E}=\int_{\Omega}\mathcal E\diff\svec{x}.
\end{equation}

The discrete approximation will be constructed so that it mimics \eqref{eq:continuous:entropy-balance}, which will guarantee that the discrete entropy of the approximation will remain bounded by the boundary and initial data.

\subsubsection{Boundary conditions}

Boundedness of the total entropy depends on proper specification of boundary conditions.
We examine the effect of free-- and no--slip boundary conditions. 
For the Cahn--Hilliard equation, we use a non--homogeneous Neumann condition for the concentration \cite{2012:Dong}, and a homogeneous Neumann condition for the chemical potential,
\begin{equation}
-\frac{3}{2}\sigma\varepsilon\frac{\partial c}{\partial \svec{n}}\biggr|_{\partial \Omega} = f_{w}'(c),~~\frac{\partial \mu}{\partial \svec{n}}\biggr|_{\partial \Omega} = 0,
\label{eq:continuous:CahnHilliardBCs}
\end{equation}
where $f_{w}(c)$ is the boundary free--energy function that controls the wall contact angle. One choice is to use the function \cite{2012:Dong}
\begin{equation}
f_{w}(c)=\frac{1}{2}\sigma\cos\theta_{w}(2c-1)(1+2c-2c^2),~~f_{w}'(c)=6\sigma\cos\theta_{w}c(1-c),
\end{equation}
where $\theta_{w}$ is the imposed contact angle with the wall. 
In most 
simulations performed in this work we use a $90^\circ$ angle, which simplifies $f_{w}(c)=0$.

For free--slip 
boundary conditions, we impose zero normal velocity and zero normal stress for momentum, while for no--slip 
walls, all velocity components are zero ($\svec{u}=0$), and we do not impose any 
conditions on the stress tensor. Either way, the boundary entropy flux is,
\begin{equation}
\begin{split}
\svec{f}^{\mathcal E}\cdot\svec{n} &= -\frac{3}{2}\sigma\varepsilon c_{t}\svec{\nabla} c\cdot\svec{n} + \left(\frac{1}{2}\rho v_{tot}^{2}+p+\mu 
c\right)\svec{u}\cdot\svec{n} + 2\eta\svec{n}\cdot\tens{S}\cdot\svec{u} + M_0\mu\svec{\nabla}\mu\cdot\svec{n}.\\
&=c_t f_{w}'(c) = \frac{\diff f_{w}}{\diff t},
\end{split}
\end{equation}
i.e., both entropy fluxes for the free-- and no--slip walls coincide. The entropy balance \eqref{eq:continuous:entropy-balance} with wall boundary conditions is therefore
\begin{equation}
\frac{\diff}{\diff t}\left(\bar{\mathcal E}+\int_{\partial\Omega}f_{w}(c)\diff S\right) = -\int_{\Omega}\left(M_0|\svec{\nabla}\mu|^2 + 2\eta\tens{S}:\tens{S}\right)\diff\svec{x} \leqslant 0,
\label{eq:continuous:entropy-balance-w-bcs}
\end{equation}
where the volume entropy is augmented with the surface free energy as in \cite{2012:Dong,2019:Manzanero-CH}.

\section{Space and time discretization}\label{sec:DGSEM}

We now construct the entropy--stable DGSEM approximation.
We restrict ourselves to the tensor product DGSEM with 
Gauss--Lobatto (GL) points, since it satisfies the Summation--By--Parts
Simultaneous--Approximation--Term (SBP--SAT) property \cite{2014:Carpenter}. The SBP--SAT property allows us to discretely follow the continuous stability steps to construct a discrete entropy law.

\subsection{Differential geometry and curvilinear elements}\label{sec:DGSEM:geometry-weak}

The physical domain $\Omega$ is tessellated with non--overlapping hexahedral elements, $e$, which are geometrically transformed from a reference element $E=[-1,1]^3$. This transformation is performed using a (polynomial) transfinite mapping $\svec{X}$ that relates physical coordinates ($\svec{x}=\left(x^1,x^2,x^3\right)=\left(x,y,z\right)=x\hat{x}+y\hat{y}+z\hat{z}$) to local reference coordinates ($\svec{\xi}=\left(\xi^1,\xi^2,\xi^3\right)=\left(\xi,\eta,\zeta\right)=\xi\hat{\xi}+\eta\hat{\eta}+\zeta\hat{\zeta}$) through
\begin{equation}
\svec{x}=\svec{X}\left(\svec{\xi}\right)=\svec{X}\left(\xi,\eta,\zeta\right).
\label{eq:Mapping}
\end{equation}
The space vectors $\hat{x}_{i}$ and $\hat{\xi}^{i}$ are unit vectors in the three Cartesian directions of physical and reference coordinates, respectively. 

From the transformation \eqref{eq:Mapping} one can define three covariant basis vectors,
\begin{equation}
\svec{a}_{i} = \frac{\partial\svec{X}}{\partial \xi^{i}},~~i=1,2,3,
\end{equation}
and three contravariant basis vectors,
\begin{equation}
\svec{a}^{i} = \svec{\nabla}\xi^{i} = \frac{1}{J}\left(\svec{a}_{j}\times\svec{a}_k\right),~~(i,j,k)\text{ cyclic},
\end{equation}
where 
\begin{equation}
J = \svec{a}_{1}\cdot\left(\svec{a}_{2}\times\svec{a}_{3}\right)
\end{equation}
is the Jacobian of the mapping $\svec{X}$. The contravariant coordinate vectors satisfy the metric identities \cite{2006:Kopriva},
\begin{equation}
\sum_{i=1}^{3}\frac{\partial \left(Ja_n^{i}\right)}{\partial\xi^{i}} = 0,~~n=1,2,3,
\label{eq:metrics:metric-id-cont}
\end{equation}
where $a_{n}^{i}$ is the $n$--th Cartesian component of the contravariant vector $\svec{a}^{i}$.

We use the volume weighted contravariant basis $J\svec{a}^{i}$ to transform differential operators from physical ($\svec{\nabla}$) to reference ($\svec{\nabla}_{\xi}$) space. The divergence of a vector is \cite{2017:Gassner}
\begin{equation}
\snabla\cdot\svec{f} = \frac{1}{J}\snabla_{\xi}\cdot\left(\tens{M}^T\svec{f}\right),
\label{eq:metrics:div-transf-1eqn}
\end{equation}
where $\tens{M}=\left(J\svec{a}^{\xi},J\svec{a}^{\eta},J\svec{a}^{\zeta}\right)$. We use \eqref{eq:notation:block-matrix-from-space-matrix} to write the divergence of an entire block vector compactly. Thus, we define the block matrix $\mathcal M$,
\begin{equation}
\mathcal M = \left(\begin{array}{ccc}Ja^{1}_{1}\smat{I}_{5} & Ja^{2}_{1}\smat{I}_{5} & Ja^{3}_{1}\smat{I}_{5}\\
Ja^{1}_{2}\smat{I}_{5} & Ja^{2}_{2}\smat{I}_{5} & Ja^{3}_{2}\smat{I}_{5}\\
Ja^{1}_{3}\smat{I}_{5} & Ja^{2}_{3}\smat{I}_{5} & Ja^{3}_{3}\smat{I}_{5}\end{array}\right),
\label{eq:metrics:metrics-matrix}
\end{equation}
which allows us to write \eqref{eq:metrics:div-transf-1eqn} for all the state variables simultaneously,
\begin{equation}
\snabla\cdot\ssvec{f} = \frac{1}{J}\snabla_{\xi}\cdot\left(\mathcal M^T\ssvec{f}\right)=\frac{1}{J}\snabla_{\xi}\cdot\cssvec{f}  ,
\label{eq:metrics:div-transf}
\end{equation}
with $\cssvec{f}$ being the block vector of the contravariant fluxes,
\begin{equation}
\cssvec{f}=\mathcal M^{T}\ssvec{f},~~\tilde{\stvec{f}}^{i} = J\svec{a}^{i}\cdot\ssvec{f}.
\label{eq:metrics:contravariant-flux}
\end{equation}
The gradient of a scalar is \cite{2017:Gassner}
\begin{equation}
\snabla w = \frac{1}{J}\tens{M}\snabla_{\xi}w,
\label{eq:metrics:gradient-transf-1eqn}
\end{equation}
which we can also extend to all entropy variables using \eqref{eq:metrics:metrics-matrix},
\begin{equation}
\ssvec{g}=\snabla \stvec{w} = \frac{1}{J}{\mathcal M}\snabla_{\xi}\stvec{w},
\label{eq:metrics:gradient-transf}
\end{equation}
and to non--conservative terms,
\begin{equation}
\ssvecg{\phi}_{m}\cdot\svec{\nabla}w_{m} = \frac{1}{\mathcal J}\cssvecg{\phi}_{m}\cdot\svec{\nabla}_{\xi}w_{m}.
\label{eq:metrics:ncterms-transf}
\end{equation}

To transform the iNS/CH system, \eqref{eq:governing:iNS/CH-system},  into reference space, we first write it as a first order system. To do so, we define the auxiliary variables $\ssvec{g}=\snabla\stvec{w}$ and $\svec{g}_{c}=\svec{\nabla}c$ so that
%DAK the letter g is already being used for the boundary values
%%
\begin{equation}
\begin{split}
&\smat{m}\stvec{q}_{t} +\snabla\cdot\ssvec{f}_{e}(\stvec{q})+\sum_{m=1}^{5}\left(\ssvecg{\phi}_{m}\left(\stvec{q}\right)\cdot\svec{\nabla}w_{m}\right) = \snabla\cdot\ssvec{f}_{v}\left(\ssvec{g}\right) + \stvec{s}(\stvec{q}),\\
&\ssvec{g}= \snabla\stvec{w},\\
&\mu=f_0'(c)-\frac{3}{2}\sigma\varepsilon\svec{\nabla}\cdot\svec{g}_{c}, \\
&\svec{g}_{c} = \svec{\nabla}c.
\end{split}
\end{equation}
Recall that the viscous fluxes of the incompressible NSE  depend only on the gradient of the entropy variables, $\ssvec{g}$, and not on the state vector $\stvec{q}$. Note that in the first equation the non--conservative terms are written in terms of $\svec{\nabla}\stvec{w}$ and not in terms of $\ssvec{g}$. 
This allows us to construct an entropy--stable scheme.

Next, we transform the operators to reference space using \eqref{eq:metrics:div-transf}, \eqref{eq:metrics:gradient-transf}, and \eqref{eq:metrics:ncterms-transf}
\begin{subequations}\label{eq:dg:system-1st-order-computational}
	\begin{align}
	&J\smat{m}\stvec{q}_{t} +\snabla_{\xi}\cdot\cssvec{f}_{e}(\stvec{q})+\sum_{m=1}^{5}\left(\cssvecg{\phi}_{m}\left(\stvec{q}\right)\cdot\svec{\nabla}_{\xi}w_{m}\right) = \snabla_{\xi}\cdot\cssvec{f}_{v}\left(\ssvec{g}\right) + J\stvec{s}(\stvec{q}),
	\label{eq:dg:system-1st-order-computational-q}\\
	&J\ssvec{g}= {\mathcal M}\snabla_{\xi}\stvec{w},
	\label{eq:dg:system-1st-order-computational-g}\\
	&J\mu=Jf_0'(c)-\frac{3}{2}\sigma\varepsilon\svec{\nabla}_{\xi}\cdot\svec{\tilde{g}}_{c}, \\
	&J\svec{g}_{c} = \tens{M}\svec{\nabla}_{\xi}c,
	\end{align}
\end{subequations}
to obtain the final form of the equations to be approximated.

The DG approximation is obtained from weak forms of the equations \eqref{eq:dg:system-1st-order-computational}. We first define the inner product in the reference element, $E$, for state and block vectors
\begin{equation}
\left\langle \stvec{f},\stvec{g}\right\rangle_{E} = \int_{E}\stvec{f}^{T}\stvec{g}\diff E,~~\left\langle\ssvec{f},\ssvec{g}\right\rangle_{E} = \int_{E}\ssvec{f}\cdot\ssvec{g}\diff E.
\label{eq:dg:weak-forms}
\end{equation}
We construct four weak forms by multiplying \eqref{eq:dg:system-1st-order-computational} by four test functions $\stvecg{\varphi}_{q}$, $\ssvecg{\varphi}_{g}$, ${\varphi}_{\mu}$, and $\svec{\varphi}_{c}$, 
then we integrate over the reference element $E$, and finally we integrate by parts to get
\begin{equation}
\begin{split}
\left\langle J\smat{m}\stvec{q}_{t},\stvecg{\varphi}_{q}\right\rangle_{E} &+\int_{\partial E}\stvecg{\varphi}_{q}^{T}\left(\cssvec{f}_{e}+\sum_{m=1}^{5}\cssvecg{\phi}_{m}w_m-\cssvec{f}_{v}\right)\cdot\hat{n}\diff S_{\xi}-\left\langle\cssvec{f}_{e},\snabla_{\xi}\stvecg{\varphi}_{q}\right\rangle_{E} \\
&-\sum_{m=1}^{5}\left\langle w_{m},\svec{\nabla}_{\xi}\cdot\left(\stvecg{\varphi}_{q}^{T}\cssvecg{\phi}_{m}\right)\right\rangle_{E}=- \left\langle\cssvec{f}_{v},\snabla_{\xi}\stvecg{\varphi}_{q}\right\rangle_{E} + \left\langle J\stvec{s},\stvecg{\varphi}_{q}\right\rangle_{E},\\
\left\langle J\ssvec{g},\ssvecg{\varphi}_{g}\right\rangle_{E}&= \int_{\partial E}\stvec{w}^{T}\left(\cssvecg{\varphi}_{g}\cdot\hat{n}\right)\diff S_{\xi} - \left\langle \stvec{w},\snabla_{\xi}\cdot\cssvecg{\varphi}_{g}\right\rangle_{E},\\
	\left\langle J\mu,\varphi_{\mu}\right\rangle_{E}&=\left\langle Jf_0',\varphi_{\mu}\right\rangle_{E}-\int_{\partial E}\frac{3}{2}\sigma\varepsilon\varphi_{\mu}\svec{\tilde{g}}_{c}\cdot\hat{n}\diff S_{\xi}+\left\langle\frac{3}{2}\sigma\varepsilon\svec{\tilde{g}}_{c},\svec{\nabla}_{\xi}\varphi_{\mu}\right\rangle_{E}, \\
\left\langle J\svec{g}_{c},\svec{\varphi}_{c}\right\rangle_{E}&=\int_{\partial E}c\svec{\tilde \varphi}_{c}\cdot\hat{n}\diff S_{\xi} - \left\langle c,\svec{\nabla}_{\xi}\cdot\svec{\tilde \varphi}_{c}\right\rangle_{E}.
\end{split}
\label{eq:dg:weak-form-1}
\end{equation}
The quantities $\hat{n}$ and $\diff S_{\xi}$ are the unit outward pointing normal and surface differential at the faces of $E$, respectively. The contravariant test functions $\cssvecg{\varphi}_{g}$ and $\svec{\tilde \varphi}_{c}$ follow the definition \eqref{eq:metrics:contravariant-flux}. Finally, surface integrals extend to all six faces of an element,
\begin{equation}
\int_{\partial E} \svec{\tilde{f}}\cdot\hat{n}\diff S_{\xi} = \int_{[-1,1]^{2}}\!\!\!\!\!\!\tilde{f}^{1}\diff \eta \diff \zeta\biggr|_{\xi=-1}^{\xi=1}\!\!+\int_{[-1,1]^{2}}\!\!\!\!\!\!\tilde{f}^{2}\diff \xi \diff \zeta\biggr|_{\eta=-1}^{\eta=1}\!\!+\int_{[-1,1]^{2}}\!\!\!\!\!\!\tilde{f}^{3}\diff \xi\diff \eta \biggr|_{\zeta=-1}^{\zeta=1}.
\label{eq:dg:surface-integral-def-cont}
\end{equation}

We can write surface integrals in either physical or reference space. The relation between physical and reference surface differentials is given by,
\begin{equation}
\diff S^i = \left| J\svec{a}^i\right|\diff\xi^{j}\diff\xi^{k} = \mathcal J_f^i \diff S_\xi^i,
\end{equation}
where we have defined the face Jacobian $\mathcal J_{f}^i = \left|\mathcal J \svec{a}^i\right|$. We can write 
the surface flux in either reference element, $\svec{\tilde f}\cdot\hat{{n}}$, or physical, $\svec{f}\cdot\svec{n}$, variables through
\begin{equation}
{\svec{\tilde f}}\cdot\hat{{n}}^{i}\diff S_\xi = \left(\boldsymbol{\mathcal M}^T\svec{f}\right)\cdot\hat{{n}}^{i}\diff S_\xi = \svec{f}\cdot\left(\boldsymbol{\mathcal M}\hat{{n}}^{i}\right)\diff S_\xi = \svec{f}\cdot\svec{n}\left|J\svec{a}^i\right|\diff S_\xi = \svec{f}\cdot\svec{n}^{i}\diff S.
\label{eq:dg:local-physical-fluxes-relationship}
\end{equation}
Therefore, the surface integrals can be written in both physical and reference 
spaces,%
\begin{equation}
\int_{\partial E}\tilde{\svec{f}}\cdot\hat{\stvec{n}}\diff S_\xi = \int_{\partial 
	e}\svec{f}\cdot\svec{n}\diff S,
\label{eq:dg:surface-integrals-relation}
\end{equation}
and we will use one or the other depending on whether we are studying an 
isolated element (reference space) or the entire mesh (physical space).

\subsection{Polynomial approximation and the DGSEM}\label{subsec:DGSEM:Approximation}

We now construct the discrete version of \eqref{eq:dg:weak-form-1}. The approximation of the state vector inside each element $E$ is an order $N$ polynomial,
\begin{equation}
\stvec{q}\approx \stvec{Q}\left(\svec{\xi}\right) = \sum_{i,j,k=0}^{N}\stvec{Q}_{ijk}(t)l_i(\xi)l_j(\eta)l_k(\zeta)\in\mathbb P^{N},
\label{eq:dg:interpolation}
\end{equation}
where $\mathbb P^{N}$ is the space of polynomials of degree less than or equal to $N$ on $[-1,1]^{3}$. The state values $\stvec{Q}_{ijk}(t)=\stvec{Q}(\xi_i,\eta_j,\zeta_k,t)$ are the nodal degrees of freedom (time dependent coefficients) at the tensor product of each of the Gauss--Lobatto (GL) points $\{\xi_i\}_{i=0}^{N}$.  Then, Lagrange polynomials $l_i(\xi)$ are
\begin{equation}
l_i(\xi) = \prod_{\substack{j=0 \\ j\neq i}}^{N}\frac{\xi-\xi_j}{\xi_i-\xi_j}.
\end{equation}

The geometry and metric terms are also approximated with order N polynomials. Let us  denote $\mathbb I^{N}$ as the polynomial interpolation operator \cite{2009:Kopriva}. The transfinite mapping is approximated using $\mathcal X=\mathbb I^{N}\left(X\right)$, but special attention must be paid to its derivatives (i.e. the contravariant basis) since $\mathcal J\svec{a}^{i}\neq\mathbb I^{N}\left(\svec{a}_{j}\times\svec{a}_k\right)$. For the metric identities \eqref{eq:metrics:metric-id-cont} to hold discretely,
\begin{equation}
\sum_{i=1}^{3}\frac{\partial \mathbb I^{N}\left(J a_n^i\right)}{\partial \xi^i} = 0,~~n=1,2,3,
\end{equation}
we approximate metric terms using the curl form \cite{2006:Kopriva},%
\begin{equation}
\mathcal J a_{n}^{i} = -\hat{x}_i \cdot\snabla_{\xi}\times\left(\mathbb I^{N}\left(\mathcal X_{l}\snabla_{\xi}\mathcal X_{m}\right)\right)\in\mathbb{P}^{N},~~~i=1,2,3,~~n=1,2,3,~~(n,m,l)\text{ cyclic}.
\label{eq:dg:contravariant-basis}
\end{equation}
If we compute $\mathcal J \svec{a}^{i}$ using \eqref{eq:dg:contravariant-basis}, we ensure discrete free--stream preservation, which is crucial to avoid grid induced solution changes \cite{2006:Kopriva}.

Next, we approximate integrals that arise in the weak formulation using Gauss quadratures. Let $\{w_i\}_{i=0}^{N}$ be the quadrature weights associated to Gauss--Lobatto nodes $\{\xi_{i}\}_{i=0}^{N}$. Then, in one dimension,
\begin{equation}
\int_{-1}^{1}f(\xi)\diff \xi \approx \int_{N}f(\xi)\diff \xi = \sum_{m=0}^{N}w_m f(\xi_{m})=\sum_{m=0}^{N}w_mF_{m}.
\end{equation}
For Gauss--Lobatto points, the approximation is exact if $f\left(\xi\right)\in\mathbb P^{2N-1}$. The extension to three dimensions has three nested quadratures, one for each of the three reference space dimensions. For example, the inner product is
\begin{equation}
\left\langle f,g\right\rangle_{E} \approx\left\langle f,g\right\rangle_{E,N} = 
\sum_{m,n,l=0}^{N}w_{mnl}F_{mnl}G_{mnl},~~w_{mnl}=w_mw_nw_l,
\label{eq:dg:discrete-inner}
\end{equation}
with a similar definition for block vectors.
The approximation of surface integrals is performed similarly, replacing exact integrals by Gauss quadratures in \eqref{eq:dg:surface-integral-def-cont},
\begin{equation}
\begin{split}
\int_{\partial E} \svec{\tilde{f}}\cdot\hat{n}\diff S_{\xi} &\approx\int_{\partial E,N} \svec{\tilde{f}}\cdot\hat{n}\diff S_{\xi} = \int_{N}\tilde{f}^{1}\diff \eta \diff \zeta\biggr|_{\xi=-1}^{\xi=1}\!\!+\int_{N}\tilde{f}^{2}\diff \xi \diff \zeta\biggr|_{\eta=-1}^{\eta=1}\!\!+\int_{N}\tilde{f}^{3}\diff \xi\diff \eta 
\biggr|_{\zeta=-1}^{\zeta=1}\\
&=\sum_{j,k=0}^{N}w_{jk}\left(\tilde{{F}}^{1}_{Njk}-\tilde{{F}}^{1}_{0jk}\right)+\sum_{i,k=0}^{N}w_{ik}\left(\tilde{{F}}^{2}_{iNk}-\tilde{{F}}^{2}_{i0k}\right)+\sum_{i,j=0}^{N}w_{ij}\left(\tilde{{F}}^{3}_{ijN}-\tilde{{F}}^{3}_{ij0}\right).
\end{split}
\label{eq:dg:surface-integral-def-disc}
\end{equation}

Gauss--Lobatto points are used to construct entropy--stable schemes using split--forms \cite{2014:Carpenter,2016:Gassner}. 
 There is no need to perform an interpolation from the 
volume polynomials in \eqref{eq:dg:surface-integral-def-disc} to the boundaries,
since boundary points are included, which is known as the Simultaneous--Approximation--Term (SAT) property. 
The exactness of the numerical quadrature and 
the SAT property yield the discrete Gauss law 
\cite{2014:Carpenter,2017:Kopriva}: for any polynomials $\ssvec{F}$ and $\stvec V$ in $\mathbb{P}^{N}$,
\begin{equation}
\left\langle \svec{\nabla}_{\xi}\cdot\cssvec{F},\stvec{V}\right\rangle_{E,N} = \int_{\partial E,N}\left(\cssvec{F}\cdot\hat{n}\right)\stvec{V}\diff S_{\xi} - \left\langle \cssvec{F},\snabla_\xi 
\stvec{V}\right\rangle_{E,N}.
\label{eq:dg:discreteGaussLaw}
\end{equation}

We start the discretization with the insertion of the discontinuous polynomial ansatz \eqref{eq:dg:interpolation}, \eqref{eq:dg:discrete-inner}, and \eqref{eq:dg:surface-integral-def-disc} into the continuous weak forms \eqref{eq:dg:weak-form-1},
\begin{equation}
\begin{split}
\left\langle \mathcal J\smat{M}\stvec{Q}_{t},\stvecg{\varphi}_{q}\right\rangle_{E,N} &+\int_{\partial E,N}\stvecg{\varphi}_{q}^{T}\left(\cssvec{F}_{e}+\sum_{m=1}^{5}\cssvecg{\Phi}_{m}W_m-\cssvec{F}_{v}\right)\cdot\hat{n}\diff S_{\xi}-\left\langle\cssvec{F}_{e},\snabla_{\xi}\stvecg{\varphi}_{q}\right\rangle_{E,N} \\
&-\sum_{m=1}^{5}\left\langle W_{m},\svec{\nabla}_{\xi}\cdot\left(\stvecg{\varphi}_{q}^{T}\cssvecg{\Phi}_{m}\right)\right\rangle_{E,N}=- \left\langle\cssvec{F}_{v},\snabla_{\xi}\stvecg{\varphi}_{q}\right\rangle_{E,N} + \left\langle \mathcal J\stvec{S},\stvecg{\varphi}_{q}\right\rangle_{E,N},\\
\left\langle \mathcal J\ssvec{G},\ssvecg{\varphi}_{g}\right\rangle_{E,N}&= \int_{\partial E,N}\stvec{W}^{T}\left(\cssvecg{\varphi}_{g}\cdot\hat{n}\right)\diff S_{\xi} - \left\langle \stvec{W},\snabla_{\xi}\cdot\cssvecg{\varphi}_{g}\right\rangle_{E,N},\\
	\left\langle \mathcal J\mu,\varphi_{\mu}\right\rangle_{E,N}&=\left\langle \mathcal JF_0',\varphi_{\mu}\right\rangle_{E,N}-\int_{\partial E,N}\frac{3}{2}\sigma\varepsilon\varphi_{\mu}\svec{\tilde{G}}_{c}\cdot\hat{n}\diff S_{\xi}+\left\langle\frac{3}{2}\sigma\varepsilon\svec{\tilde{G}}_{c},\svec{\nabla}_{\xi}\varphi_{\mu}\right\rangle_{E,N}, \\
\left\langle \mathcal J\svec{G}_{c},\svec{\varphi}_{c}\right\rangle_{E,N}&=\int_{\partial E,N}C\svec{\tilde \varphi}_{c}\cdot\hat{n}\diff S_{\xi} - \left\langle C,\svec{\nabla}_{\xi}\cdot\svec{\tilde 
\varphi}_{c}\right\rangle_{E,N}.
\end{split}
\label{eq:dg:discrete-weak-form}
\end{equation}
In \eqref{eq:dg:discrete-weak-form}, the test functions are restricted to polynomial spaces, $\mathbb P^{N}$. Taken in context, $\mu$  (which lacks an upper case symbol) 
refers to the polynomial approximation of the chemical potential. 

Euler conservative fluxes and viscous fluxes have inter--element coupling and physical boundary conditions enforced 
by \textit{numerical fluxes} in the element boundary quadratures in \eqref{eq:dg:discrete-weak-form},
\begin{equation}
\begin{split}
&\ssvec{F}_{e}\approx{\ssvec{F}}_{e}^{\star}\left(\stvec{Q}_L,\stvec{Q}_R\right),~~ \ssvec{F}_{v}\approx{\ssvec{F}}_{v}^{\star}\left(\ssvec{G}_L,\ssvec{G}_R\right),~~ \stvec{W}\approx 
\stvec{W}^{\star}\left(\stvec{Q}_L,\stvec{Q}_R\right),\\
&\svec{G}_{c}\approx\svec{G}_{c}^{\star}\left(\svec{G}_{cL},\svec{G}_{cR}\right),~~C\approx 
C^{\star}\left(C_L,C_R\right).
\label{eq:dg:numerical-fluxes}
\end{split}
\end{equation}
Whereas for non--conservative terms, we follow \cite{2018:Bohm} and use \textit{diamond 
fluxes} at the boundaries,
\begin{equation}
  \ssvecg{\Phi}_{m}W_m \approx 
  \left(\ssvec{\Phi}_{m}W_{m}\right)^{\diamondsuit}\left(\stvec{Q}_{L},\stvec{Q}_{R}\right).
  \label{eq:dg:diamond-fluxes}
\end{equation}
Both numerical and diamond flux functions are detailed below in Sec. \ref{subsec:NumericalFluxes}. 
In \eqref{eq:dg:numerical-fluxes} and \eqref{eq:dg:diamond-fluxes}, $L$ and $R$ represent the values from
the left and right adjacent elements.
Although numerical fluxes are single valued at each interface, diamond 
fluxes are not so constrained and their value can jump from one side to the other.

Inserting the numerical  \eqref{eq:dg:numerical-fluxes} and diamond \eqref{eq:dg:diamond-fluxes} fluxes into \eqref{eq:dg:discrete-weak-form} 
completes the semi--discretization,
\begin{subequations}\label{eq:dg:discrete-weak-form-numflux}
	\begin{align}
\left\langle \mathcal J\smat{M}\stvec{Q}_{t},\stvecg{\varphi}_{q}\right\rangle_{E,N} &+\int_{\partial E,N}\stvecg{\varphi}_{q}^{T}\left(\cssvec{F}_{e}^{\star}+\sum_{m=1}^{5}\left(\cssvecg{\Phi}_{m}W_m\right)^{\diamondsuit}-\cssvec{F}_{v}^{\star}\right)\cdot\hat{n}\diff S_{\xi}-\left\langle\cssvec{F}_{e},\snabla_{\xi}\stvecg{\varphi}_{q}\right\rangle_{E,N} \nonumber\\
&-\sum_{m=1}^{5}\left\langle W_{m},\svec{\nabla}_{\xi}\cdot\left(\stvecg{\varphi}_{q}^{T}\cssvecg{\Phi}_{m}\right)\right\rangle_{E,N}=- \left\langle\cssvec{F}_{v},\snabla_{\xi}\stvecg{\varphi}_{q}\right\rangle_{E,N} + \left\langle \mathcal J\stvec{S},\stvecg{\varphi}_{q}\right\rangle_{E,N},\label{eq:dg:discrete-weak-form-numflux:q}\\
\left\langle \mathcal J\ssvec{G},\ssvecg{\varphi}_{g}\right\rangle_{E,N}=& \int_{\partial E,N}\stvec{W}^{\star,T}\left(\cssvecg{\varphi}_{g}\cdot\hat{n}\right)\diff S_{\xi} - \left\langle \stvec{W},\snabla_{\xi}\cdot\cssvecg{\varphi}_{g}\right\rangle_{E,N},\label{eq:dg:discrete-weak-form-numflux:g}\\
	\left\langle \mathcal J\mu,\varphi_{\mu}\right\rangle_{E,N}=&\left\langle \mathcal JF_0',\varphi_{\mu}\right\rangle_{E,N}-\int_{\partial E,N}\frac{3}{2}\sigma\varepsilon\varphi_{\mu}\svec{\tilde{G}}_{c}^{\star}\cdot\hat{n}\diff S_{\xi}+\left\langle\frac{3}{2}\sigma\varepsilon\svec{\tilde{G}}_{c},\svec{\nabla}_{\xi}\varphi_{\mu}\right\rangle_{E,N} \label{eq:dg:discrete-weak-form-numflux:mu},\\
\left\langle \mathcal J\svec{G}_{c},\svec{\varphi}_{c}\right\rangle_{E,N}=&\int_{\partial E,N}C^{\star}\svec{\tilde \varphi}_{c}\cdot\hat{n}\diff S_{\xi} - \left\langle C,\svec{\nabla}_{\xi}\cdot\svec{\tilde 
\varphi}_{c}\right\rangle_{E,N}.\label{eq:dg:discrete-weak-form-numflux:c}
	\end{align}
\end{subequations}
Because it enhances the algorithm efficiency and makes it easier to show stability, we apply the discrete Gauss law \eqref{eq:dg:discreteGaussLaw}
to the non--conservative terms and inviscid fluxes of 
\eqref{eq:dg:discrete-weak-form-numflux:q}, in 
\eqref{eq:dg:discrete-weak-form-numflux:g}, and in 
\eqref{eq:dg:discrete-weak-form-numflux:c}, and use \eqref{eq:dg:surface-integrals-relation} to write surface integrals in 
physical variables,
\begin{subequations}\label{eq:dg:semi-discrete-approx}
	\begin{align}
\left\langle \mathcal J\smat{M}\stvec{Q}_{t},\stvecg{\varphi}_{q}\right\rangle_{E,N} &+\int_{\partial e,N}\stvecg{\varphi}_{q}^{T}\left(\ssvec{F}_{e}^{\star}-\ssvec{F}_{e}+\sum_{m=1}^{5}\left(\ssvecg{\Phi}_{m}W_m\right)^{\diamondsuit}-\ssvecg{\Phi}_{m}W_{m}\right)\cdot\svec{n}\diff S \nonumber\\
&+\left\langle\stvecg{\varphi}_{q},\svec{\nabla}_{\xi}\cdot\cssvec{F}_{e}\right\rangle_{E,N}+\sum_{m=1}^{5}\left\langle\stvecg{\varphi}_{q},\cssvecg{\Phi}_{m}\cdot\svec{\nabla}_{\xi}W_{m}\right\rangle_{E,N} \nonumber\\
=&\int_{\partial e,N}\stvecg{\varphi}^{T}_{q}\ssvec{F}_{v}^{\star}\cdot\svec{n}\diff S- \left\langle\cssvec{F}_{v},\snabla_{\xi}\stvecg{\varphi}_{q}\right\rangle_{E,N} + \left\langle \mathcal J\stvec{S},\stvecg{\varphi}_{q}\right\rangle_{E,N},\label{eq:dg:semi-discrete-approx:q}\\
\left\langle \mathcal J\ssvec{G},\ssvecg{\varphi}_{g}\right\rangle_{E,N}=& \int_{\partial e,N}\left(\stvec{W}^{\star,T}-\stvec{W}^{T}\right)\left(\ssvecg{\varphi}_{g}\cdot\svec{n}\right)\diff S + \left\langle \cssvecg{\varphi}_{g},\svec{\nabla}_{\xi}\stvec{W}\right\rangle_{E,N},\label{eq:dg:semi-discrete-approx:g}\\
	\left\langle \mathcal J\mu,\varphi_{\mu}\right\rangle_{E,N}=&\left\langle \mathcal JF_0',\varphi_{\mu}\right\rangle_{E,N}-\int_{\partial e,N}\frac{3}{2}\sigma\varepsilon\varphi_{\mu}\svec{{G}}_{c}^{\star}\cdot\svec{n}\diff S+\left\langle\frac{3}{2}\sigma\varepsilon\svec{\tilde{G}}_{c},\svec{\nabla}_{\xi}\varphi_{\mu}\right\rangle_{E,N} \label{eq:dg:semi-discrete-approx:mu},\\
\left\langle \mathcal J\svec{G}_{c},\svec{\varphi}_{c}\right\rangle_{E,N}=&\int_{\partial e,N}\left(C^{\star}-C\right)\svec{\varphi}_{c}\cdot\svec{n}\diff S + \left\langle \svec{\tilde 
\varphi}_{c},\svec{\nabla}_{\xi}C\right\rangle_{E,N}.\label{eq:dg:semi-discrete-approx:c}
	\end{align}
\end{subequations}
The algorithm efficiency is improved by the use of single calculation of the local gradient of the entropy variables, $\svec{\nabla}_{\xi}\stvec{W}$, for both non--conservative terms in \eqref{eq:dg:semi-discrete-approx:q} 
and gradients in \eqref{eq:dg:semi-discrete-approx:g}. As a result, 
non--conservative terms are as expensive as source terms (i.e. no additional matrix multiplications are 
required). Then for the gradients $\ssvec{G}$ we augment the
 pre--computed local gradients $\svec{\nabla}_{\xi}\stvec{W}$ with the surface integral in 
 \eqref{eq:dg:semi-discrete-approx:g}.
 
 The point--wise discretization is obtained by replacing the test function in \eqref{eq:dg:semi-discrete-approx} by the Lagrange polynomials. For completeness, we write the implementation of the discretization in Appendix \ref{app:point-wise}.

\subsection{Numerical fluxes}\label{subsec:NumericalFluxes}

The approximation \eqref{eq:dg:semi-discrete-approx} is completed with the specification of numerical fluxes ${\stvec{F}}_{e}^{\star}$, ${\stvec{F}}_{v}^{\star}$, $\stvec{W}^{\star}$, $\svec{G}_{c}^{\star}$, $C^{\star}$, and diamond fluxes 
$(\ssvec{\Phi}_{m}W_m)^{\diamondsuit}$.
For the inviscid fluxes we propose two options: an entropy conserving option using using central fluxes, and an entropy--stable approximation using the exact Riemann solver derived in \cite{2017:Bassi}. 
For the viscous fluxes and concentration gradient, we use 
the Bassi--Rebay 1 (BR1) scheme \cite{BR1}. 

Before we write the numerical fluxes, we first describe the rotational 
invariance property satisfied by inviscid fluxes and non--conservative terms.
The rotational invariance of the flux \cite{2009:Toro} allows us to write the 
normal flux $\ssvec{F}_{e}\cdot\svec{n}$ from a rotated version of the 
inviscid flux  $x$--component $\stvec{F}_{e}$,
\begin{equation}
\ssvec{F}_{e}\cdot\svec{n} = \smat{T}^{T}\stvec{F}_{e}\left(\smat{T}\stvec{Q}\right)=\smat{T}^{T}\stvec{F}_{e}\left(\stvec{Q}_{n}\right),~~\smat{T} = \left(\begin{array}{ccccc}1 & 0 & 0 & 0 & 0 \\
0 & n_{x} & n_{y} & n_{z} & 0 \\
0 & t_{1,x} & t_{1,y} & t_{1,z} & 0 \\
0 & t_{2,x} & t_{2,y} & t_{2,z} & 0 \\
0 & 0 & 0 & 0 & 1
\end{array}\right),
\label{eq:Riemann:rotational-invariance}
\end{equation}
where $\smat{T}$ is a rotation matrix that only affects velocities, $\svec{n}=\left(n_x,n_y,n_z\right)$ is the normal unit vector to the face, and $\svec{t}_{1}$ and $\svec{t}_{2}$ are two tangent unit vectors to the face. 
When the rotation matrix $\smat{T}$ multiplies the state vector $\stvec{Q}$, we obtain the face normal state vector $\stvec{Q}_{n}$,
\begin{equation}
\stvec{Q}_n = \smat{T}\stvec{Q}=\left(C, \sqrt{\rho} U_{n}, \sqrt{\rho} V_{t1}, \sqrt{\rho} V_{t2}, P\right),
\label{eq:Riemann:rotated-state}
\end{equation}
where $U_n = \svec{U}\cdot\svec{n}$ the normal velocity, and $V_{ti}=\svec{U}\cdot\svec{t}_{i}$ ($i=1,2$) are the two tangent velocities. Note that the reference system rotation does not affect the total speed 
\begin{equation}
V_{tot}^2=U^2+V^2+W^2=U_n^2+V_{t1}^2+V_{t2}^2.
\end{equation}

The non--conservative terms,
\begin{equation}
  \sum_{m=1}^{5} \ssvec{\Phi}_{m}W_{m} = \Upsilon\left(\begin{array}{ccc}0 & 0 & 0 
  \\
  \frac{1}{2}\rho U^2 + \mu C & \frac{1}{2}\rho UV & \frac{1}{2}\rho UW \\
  \frac{1}{2}\rho UV & \frac{1}{2}\rho V^2 + \mu C & \frac{1}{2}\rho VW \\
  \frac{1}{2}\rho UW & \frac{1}{2}\rho VW & \frac{1}{2}\rho W^2 + \mu C\\
  \rho_0c_0^2U & \rho_0c_0^2V & \rho_0c_0^2W
  \end{array}
  \right),
\end{equation}
are also rotationally invariant,
 \begin{equation}
     \sum_{m=1}^{5} \left(\ssvec{\Phi}_{m}W_{m}\right)\cdot\svec{n} = \smat{T}^{T}\left(\begin{array}{c} 
     0 \\
     \frac{1}{2}\rho U_n^2 + \mu C\\
     \frac{1}{2}\rho U_n V_{t1} \\
     \frac{1}{2}\rho U_n V_{t2} \\
     \rho_0c_0^2U_{n}
     \end{array}\right) = \smat{T}^{T}\stvec{F}_{\Phi 
     W}\left(\stvec{Q}_{n}\right),
 \end{equation}
 with its equivalent $x$--component,
 
 \begin{equation}
   \stvec{f}_{\phi w}\left(\stvec{q}\right) = \left(\begin{array}{c} 0 \\ \frac{1}{2}\rho u^2 + \mu c \\
   \frac{1}{2}\rho uv \\
   \frac{1}{2}\rho uw \\
   \rho_0c_0^2u\end{array}\right).
 \end{equation}
 
 \subsubsection{Inviscid fluxes: entropy conserving central fluxes}
 
 The first choice for the inviscid numerical and diamond fluxes is to use central 
 fluxes, which will lead to an entropy conserving approximation. We adapt the approach used for the resistive MHD in \cite{2018:Bohm}, so that
 \begin{equation}
   \ssvec{F}_{e}^{\star}\cdot\svec{n} = \aver{\ssvec{F}_{e}}\cdot\svec{n}=\smat{T}^{T}\left(\begin{array}{c}
   \aver{CU_n} \\
   \aver{\frac{1}{2}\rho U_n^2 + p}\\
   \aver{\frac{1}{2}\rho U_nV_{t1}} \\
      \aver{\frac{1}{2}\rho U_nV_{t2}} \\
      0
  \end{array} \right),
     \label{eq:riemann:central-cons}
 \end{equation}
 and
 \begin{equation}
\sum_{m=1}^{5}   \left(\ssvec{\Phi}_{m}W_{m}\right)^{\diamondsuit}\cdot\svec{n}
   = \sum_{m=1}^{5}\ssvec{\Phi}_{m}\aver{W_m}\cdot\svec{n} = \smat{T}^{T}\left(\begin{array}{c}
   0 \\
   \frac{1}{2}\rho U_n\aver{U_n}+C\aver{\mu} \\
   \frac{1}{2}\rho U_n\aver{V_{t1}} \\
   \frac{1}{2}\rho U_n\aver{V_{t2}} \\
   \rho_0 c_0^2\aver{U_n}
   \end{array}\right).
   \label{eq:riemann:central-noncons}
 \end{equation}
 In \eqref{eq:riemann:central-cons} and \eqref{eq:riemann:central-noncons}, $\aver{u}$ 
 represents the average operator,
 \begin{equation}
 \aver{u} = \frac{u_L + u_R}{2}
 \end{equation}
 and $\ssvec{\Phi}_m$ (without average) represents the element local state at the 
 face.

 \subsubsection{Inviscid fluxes: entropy--stable exact Riemann solver (ERS)}
The second option uses the solution of the exact Riemann problem derived in \cite{2017:Bassi}.
The star region solution of the incompressible NSE (i.e. without the Cahn--Hilliard equation) is,
\begin{equation}
\begin{split}
U_n^{\star} &= \frac{P_{L}-P_{R} + \rho_L U_{nL}\lambda_{L}^{+}-\rho_{R}U_{nR}\lambda_{R}^{-}}{\rho_L\lambda_L^{+}-\rho_R \lambda_R^{-}}, ~~ P^{\star} = P_L + \rho_L \lambda^{+}_L\left(U_{nL} - U_n^\star\right),\\
\rho^{\star} &= \left\{\begin{array}{ccc}\rho_{L}^{\star} & \text{ if } & U_n^\star \geqslant 0 \\\rho_{R}^{\star} & \text{ if } & U_n^\star < 0 
\end{array}\right.,~~ \rho_L^{\star} = \frac{\rho_L \lambda_L^{+}}{U_n^{\star}-\lambda_L^-},~~\rho_R^{\star} = \frac{\rho_R \lambda_R^{-}}{U_n^{\star}-\lambda_R^+},~~V_{ti}^{\star}=\left\{\begin{array}{ccc}V_{tiL} & \text{ if } & U_n^\star \geqslant 0 \\V_{tiR} & \text{ if } & U_n^\star < 0 
\end{array}\right. ,
\end{split}
\label{eq:Riemann:star-region-solution}
\end{equation}
with eigenvalues,
\begin{equation}
\lambda_L^{\pm} = \frac{U_{nL} \pm a_L}{2},~~\lambda_R^{\pm} = \frac{U_{nR} \pm a_R}{2},~~ a = \sqrt{U_n^2 + \frac{4\rho_0c_0^2}{\rho}}.
\end{equation}

We evaluate the numerical flux with the star region solution with the normal state $\stvec{Q}_{n}$ in 
\eqref{eq:Riemann:rotational-invariance} for the momentum, and perform 
the average for the Cahn-Hilliard part,
\begin{equation}
  \ssvec{F}_{e}^{\star}\cdot\svec{n} = \smat{T}^{T}\left(\begin{array}{c}
  \aver{C U_n} \\
  \frac{1}{2}\rho^{\star}U_n^{\star,2} + p^{\star} \\
  \frac{1}{2}\rho^{\star}U_n^{\star}V_{t1}^{\star} \\
  \frac{1}{2}\rho^{\star}U_n^{\star}V_{t2}^{\star} \\
  0
  \end{array}\right).
  \label{eq:riemann:ers-cons}
\end{equation}
For non--conservative terms, we choose the diamond fluxes,
\begin{equation}
\begin{split}
  \sum_{m=1}^{5}\left(\ssvec{\Phi}_{m}W_{m}\right)^{\diamondsuit}\cdot\svec{n} &= 
  \ssvec{\Phi}_{1}\aver{W_{1}}\cdot\svec{n} + \sum_{m=2}^{5}\left(\ssvec{\Phi}_{m}^{\star}W^{\star}_m+\ssvec{\Phi}_m W_m - 
  \ssvec{\Phi}_{m}^{\star}W_m\right)\cdot\svec{n}\\
  &=\smat{T}^{T}\left(\begin{array}{c}
  0 \\
  \frac{1}{2}\rho^{\star}U_{n}^{\star,2} + \frac{1}{2}\rho U_{n}^{2} - 
  \frac{1}{2}\rho^{\star}U_n^{\star}U_n + C\aver{\mu} \\
  \frac{1}{2}\rho^{\star}U_{n}^{\star}V_{t1}^{\star} + \frac{1}{2}\rho U_{n}V_{t1} - 
  \frac{1}{2}\rho^{\star}U_n^{\star}V_{t1} \\
    \frac{1}{2}\rho^{\star}U_{n}^{\star}V_{t2}^{\star} + \frac{1}{2}\rho U_{n}V_{t2} - 
  \frac{1}{2}\rho^{\star}U_n^{\star}V_{t2} \\
  \rho_0 c_0^2 U_n^{\star}
  \end{array}\right)=\smat{T}^{T}\stvec{F}_{\Phi W}^{\mathrm{ERS}},
  \end{split}
  \label{eq:riemann:ers-noncons}
\end{equation}
where $\ssvec{\Phi}_{m}^{\star}$ and $W_m^{\star}$ refer to non--conservative 
coefficients and entropy variables evaluated with the star region solution \eqref{eq:Riemann:star-region-solution}.
This choice is justified by the stability analysis.

\subsubsection{Viscous and chemical potential fluxes: Bassi--Rebay 1 (BR1) method}

For the viscous fluxes and the chemical potential we use the Bassi--Rebay 1 (BR1) scheme \cite{BR1}, 
which averages entropy variables and fluxes between the adjacent elements. However, for the Cahn--Hilliard equation, we include the possibility to add interface stabilization. Although it is not a requirement for stability, we have found it enhances accuracy when the flow configuration is under--resolved \cite{2019:Manzanero-CH}. The viscous and Cahn--Hilliard fluxes are,
\begin{equation}
\stvec{W}^{\star} = \aver{\stvec{W}},~~{\ssvec{F}}_{v}^{\star} = 
\aver{{\ssvec{F}}_{v}}+\beta\jump{\mu}\stvec{e}_{1}\svec{n}_{L},~~C^{\star} = \aver{C},~~\svec{G}_{c}^{\star} = 
\aver{\svec{G}_{c}}+\beta\jump{C}\svec{n}_{L},
\label{eq:Riemann:BR1}
\end{equation}
where $\beta$ is the penalty parameter, computed in this work following \cite{2018:Manzanero,2019:Manzanero-CH},
\begin{equation}
\beta=\kappa_{\beta}\frac{N(N+1)}{2}|\mathcal J_{f}|\aver{\mathcal J^{-1}},
\label{eq:dg:penalty-parameter}
\end{equation}
with $\kappa_{\beta}$ a dimensionless free parameter (in this work we only use $\kappa_{\beta}=0$ to disable interface stabilization, and $\kappa_{\beta}=1$ to enable interface stabilization), we define 
the inter--element jumps as $\jump{\bullet}=\bullet_{R}-\bullet_{L}$, and $\svec{n}_{L}$ is the outward normal vector to the left element. In \eqref{eq:dg:penalty-parameter} $N$ is the polynomial order, $|\mathcal J_{f}|$ is the surface Jacobian of the face, and $\aver{\mathcal J^{-1}}$ is the average of the inverse of the Jacobians of left and right elements.

\subsection{Boundary conditions}\label{sec:BoundaryConditions}

The approximation is completed with the addition of boundary conditions. 
Here we show how we impose free-- and no--slip wall boundary conditions. The boundary conditions are imposed through the numerical and diamond fluxes at the physical boundary faces. We consider the inviscid and viscous fluxes individually. 

\subsubsection{Inviscid flux}

The inviscid numerical flux controls the normal velocity $\svec{u}\cdot\svec{n}=0$  
for both free-- and no--slip boundary conditions.  Again, as at interior faces, we 
provide an entropy conserving option with central fluxes, and an entropy stable 
version with the exact Riemann problem solution. Either way, we apply central fluxes 
\eqref{eq:riemann:central-cons}--\eqref{eq:riemann:central-noncons} or the ERS
\eqref{eq:riemann:ers-cons}--\eqref{eq:riemann:ers-noncons} to the interior 
state $\stvec{Q}_{n}^{i}$ and a mirrored ghost state $\stvec{Q}_{n}^{e}$,
\begin{equation}
  \stvec{Q}_{n}^{i} = \left(\begin{array}{c}C \\ \sqrt{\rho}U_{n} \\ \sqrt{\rho}V_{t1} \\ \sqrt{\rho}V_{t2} \\ P\end{array}\right),~~ \stvec{Q}_{n}^{e} = \left(\begin{array}{c}C \\ -\sqrt{\rho}U_{n} \\ \sqrt{\rho}V_{t1} \\ \sqrt{\rho}V_{t2} \\ 
  P\end{array}\right).
  \label{eq:BCs:inviscid-wall-ghost}
\end{equation}
The fluxes differ depending on which numerical flux is used.
\begin{enumerate}
  \item \emph{The entropy conserving approximation with central fluxes.} The 
  numerical \eqref{eq:riemann:central-cons} and diamond \eqref{eq:riemann:central-noncons} 
  central fluxes for \eqref{eq:BCs:inviscid-wall-ghost} 
  are,
\begin{equation}
  \ssvec{F}_{e}^{\star}\cdot\svec{n} =\smat{T}^{T}\left(\begin{array}{c}0 \\
  \frac{1}{2}\rho U_n^2 + P \\ 0 \\ 0 \\ 0 \end{array}\right),~~\sum_{m=1}^{5}\left(\ssvec{\Phi}_{m}W_{m}\right)^{\diamondsuit}\cdot\svec{n}
  =\smat{T}^{T}\left(\begin{array}{c}0 \\
C\mu \\
\frac{1}{2}\rho U_n V_{t1}\\
\frac{1}{2}\rho U_n V_{t2}\\
0
  \end{array}\right).
  \label{eq:BCs:inviscid-central}
\end{equation}

\item \emph{The entropy stable approximation with the exact Riemann 
solver.} The star region values computed from the states \eqref{eq:BCs:inviscid-wall-ghost} are
\begin{equation}
  U_{n}^{\star}=0,~~P^{\star}=P+\rho\lambda^{+}U_{n},~~\lambda^{+} = 
  \frac{U_n+a}{2}.
\end{equation}
With these star values the numerical \eqref{eq:riemann:ers-cons} and diamond \eqref{eq:riemann:ers-noncons} fluxes become
\begin{equation}
  \ssvec{F}_{e}^{\star}\cdot\svec{n} = \smat{T}^{T}\left(\begin{array}{c}0\\
  P+\rho\lambda^{+}U_n \\ 0 \\ 0 \\ 0 
  \end{array}\right),~~\sum_{m=1}^{5}\left(\ssvec{\Phi}_{m}W_{m}\right)^{\diamondsuit}\cdot\svec{n} 
  = \smat{T}^{T}\left(\begin{array}{c}0 \\
  \frac{1}{2}\rho U_n^{2}+C\mu \\
  \frac{1}{2}\rho U_n V_{t1} \\
  \frac{1}{2}\rho U_n V_{t2} \\
 0
\end{array}  \right).
  \label{eq:BCs:inviscid-ers}
\end{equation}
\end{enumerate}

\subsubsection{Viscous and Cahn--Hilliard fluxes}\label{subsubsec:DGSEM:BC-wall-viscous}

For the free--slip wall boundary condition, we impose zero normal stress 
$\tens{S}\cdot\svec{n}=0$ (Neumann), whereas for the no--slip wall boundary condition, we 
impose zero velocity $\svec{u}=0$ (Dirichlet). Furthermore, we 
implement the discrete version of the homogeneous and non--homogeneous Neumann 
boundary conditions \eqref{eq:continuous:CahnHilliardBCs} for the Cahn--Hilliard equation. The concentration, 
concentration gradient, entropy variables, and viscous fluxes at the boundaries 
are:
\begin{enumerate}
  \item \emph{Free--slip wall boundary condition.}
  \begin{equation}
    C^{\star}=C,~~\svec{G}_{c}^{\star}=g'(C),~~\stvec{W}^{\star}=\stvec{W},~~\ssvec{F}^{\star}_{v}\cdot\svec{n} 
    = 0.
    \label{eq:BCs:viscous-free-slip}
  \end{equation}
  
  \item \emph{No--slip wall boundary condition.}
  \begin{equation}
    C^{\star}=C,~~\svec{G}_{c}^{\star}=g'(C),~~
    \stvec{W}^{\star}=\left(\mu,0,0,0,P\right),~~
    \ssvec{F}_{v}^{\star}\cdot\svec{n}=\left(\begin{array}{c} 0 \\ 
    2\eta\mathcal{S} \\
    0
    \end{array}\right)\cdot\svec{n}.
    \label{eq:BCs:viscous-no-slip}      
  \end{equation}
\end{enumerate}

\subsection{Implicit--Explicit (IMEX) time discretization}

The fully--discrete scheme is completed with 
the discretization of the time derivatives $\stvec{Q}_{t}$ in \eqref{eq:dg:semi-discrete-approx}. 
The numerical stiffness induced by the fourth order derivatives present in the Cahn--Hilliard equation prevents us from practically using fully 
explicit time marching. Therefore, we use the 
IMplicit--EXplicit (IMEX) 
Backward Differentiation Formula (BDF) time integrator described in \cite{2018:Dong}. 
We revisit the continuous setting \eqref{eq:governing:iNS/CH-system}
\begin{equation}
  \begin{split}
   & c_{t} + \svec{\nabla}\cdot\left(c\svec{u}\right) = M_0\svec{\nabla}^{2}\left(f_0'(c)- 
    \frac{3}{2}\sigma\varepsilon\svec{\nabla}^{2}c\right),\\
    &\sqrt{\rho}\left(\sqrt{\rho}\svec{u}\right)_{t} + \svec{\nabla}\cdot\left(\frac{1}{2}\rho\svec{u}\svec{u}\right)
    +\frac{1}{2}\rho \svec{u}\cdot\svec{\nabla}\svec{u} + 
    c\svec{\nabla}\mu =-\svec{\nabla}p +  2\eta\tens{S} + \rho\svec{g},\\
    &p_{t} + \rho_0 c_0^{2}\svec{\nabla}\cdot\svec{u} = 0,
  \end{split}
\end{equation}
 to describe the approximation
in time.

The IMEX version of the BDF integrator for a time derivative $dy/dt$ is
\begin{equation}
  \gamma_0 = \left\{\begin{array}{ccc}
  1 & \text{if} & J=1  \\
  \frac{3}{2} & \text{if} & J=2
  \end{array} \right.,~~\hat{y} =\left\{\begin{array}{ccc}
y^{n} & \text{if} & J=1 \\  
2y^{n}-\frac{1}{2}y^{n-1} & \text{if} & J=2 \end{array}\right.
 ,~~y^{n+1,e} = \left\{\begin{array}{ccc} y^{n}  & \text{if} 
  & J=1 \\
  2y^{n}-y^{n-1} & \text{if} & J=2 \end{array}\right.,
\end{equation}
where $y^{n}=y(t_{n})$, $\gamma_0$ and $\hat{y}$ yield a $J$--th order approximation of the time 
derivative,
\begin{equation}
y_{t}^{n+1} \approx \frac{\gamma_0 y^{n+1}-\hat{y}}{\Delta t},
\end{equation}
and $y^{n+1,e}$ is a $J$--th order explicit approximation of $y^{n+1}$. The IMEX 
scheme used in \cite{2018:Dong} was constructed to solve the fourth order derivative implicitly, and the 
remaining terms explicitly,
\begin{equation}
  \begin{split}
    \frac{\gamma_0 c^{n+1}-\hat{c}}{\Delta t} + \svec{\nabla}\cdot\left(c\svec{u}\right)^{n+1,e} &= M_0\svec{\nabla}^{2}\left(f_0'(c^{n+1,e})- 
    \frac{3}{2}\sigma\varepsilon\svec{\nabla}^{2}c^{n+1}\right),\\
    \sqrt{\rho^{n+1,e}}\frac{\gamma_0\sqrt{\rho}\svec{u}^{n+1}-\reallywidehat{\sqrt{\rho}\svec{u}}}{\Delta t} +&\left(\svec{\nabla}\cdot\left(\frac{1}{2}\rho\svec{u}\svec{u}\right)
    +\frac{1}{2}\rho \svec{u}\cdot\svec{\nabla}\svec{u} + 
    c\svec{\nabla}\mu\right)^{n+1,e} \\
 &=\left(-\svec{\nabla}p +  2\eta\tens{S} + \rho\svec{g}\right)^{n+1,e},\\
    \frac{\gamma_0 p^{n+1}-\hat{p}}{\Delta t} + \rho_0 c_0^{2}\svec{\nabla}\cdot\svec{u}^{n+1,e} &= 0.
  \end{split}
\end{equation}
As in the vast majority of BDF implementations, the first step is first order ($J=1$), and the 
following steps can be maintained first or use second order ($J=2$). Note that the Jacobian 
matrix associated with the fourth order derivative is constant in time. Hence, only one 
computation of its solution is required. In this work, we perform a single LU decomposition, and afterwards use Gauss 
substitution to solve the linear system at each time step, but other direct or iterative solvers could be used.

One can use a fully explicit time integrator if the time step restriction allows it (e.g. when the chemical characteristic time is large enough). We have also implemented a low--storage third--order Runge--Kutta method to explicitly solve the system \cite{1980:Williamson}.

\section{Semi--discrete stability analysis}\label{sec:Stability}

We prove the stability of the spatial discretization 
\eqref{eq:dg:semi-discrete-approx}, by showing that it mimics a discrete version of 
the entropy conservation law \eqref{eq:continuous:entropy-balance}. As in the continuous analysis,
we do not study the effect of lower order source terms ($\stvec{S}=0$). The 
entropy analysis reproduces the continuous analysis steps 
discretely. Therefore, we use proper choices for the test functions, which 
contract the system of equations into a single equation, to be 
shaped into an entropy equation. The latter is produced in four steps, as we 
study the time derivative terms, inviscid terms, viscous terms, and finally 
interior and physical boundary terms.

First, we follow \cite{2019:Manzanero-CH} 
and take the time derivative of \eqref{eq:dg:semi-discrete-approx:c},
\begin{equation}
\begin{split}
  \left\langle \mathcal J\svec{G}_{c,t},\svec{\varphi}_{c}\right\rangle_{E,N}&=\int_{\partial e,N}\left(C^{\star}_{t}-C_{t}\right)\svec{\varphi}_{c}\cdot\svec{n}\diff S + \left\langle \svec{\tilde 
\varphi}_{c},\svec{\nabla}_{\xi}C_{t}\right\rangle_{E,N},
\end{split}
\end{equation}
and then replace the test functions. Following \cite{2017:Gassner,2019:Manzanero-iNS}, we insert $\stvecg{\varphi}_{q}=\stvec{W}$, 
and $\ssvecg{\varphi}_{g}=\ssvec{F}_{v}$,  and following \cite{2019:Manzanero-CH}, $\varphi_{\mu}=C_{t}$, and 
$\svec{\varphi}_{c}=\frac{3}{2}\sigma\varepsilon\svec{G}_{c}$,
\begin{subequations}\label{eq:stability:replace-test}
	\begin{align}
\left\langle \mathcal J\smat{M}\stvec{Q}_{t},\stvec{W}\right\rangle_{E,N} &+\int_{\partial e,N}\stvec{W}^{T}\left(\ssvec{F}_{e}^{\star}-\ssvec{F}_{e}+\sum_{m=1}^{5}\left(\ssvecg{\Phi}_{m}W_m\right)^{\diamondsuit}-\ssvecg{\Phi}_{m}W_{m}\right)\cdot\svec{n}\diff S \nonumber\\
&+\left\langle \stvec{W},\svec{\nabla}_{\xi}\cdot\cssvec{F}_{e}\right\rangle_{E,N}+\sum_{m=1}^{5}\left\langle \stvec{W},\cssvecg{\Phi}_{m}\cdot\svec{\nabla}_{\xi}W_{m}\right\rangle_{E,N} \nonumber\\
=&\int_{\partial e,N}\stvec{W}^{T}\ssvec{F}_{v}^{\star}\cdot\svec{n}\diff S- \left\langle\cssvec{F}_{v},\snabla_{\xi}\stvec{W}\right\rangle_{E,N},\label{eq:stability:replace-test:q}\\
\left\langle \mathcal J\ssvec{G},\ssvec{F}_{v}\right\rangle_{E,N}=& \int_{\partial e,N}\left(\stvec{W}^{\star,T}-\stvec{W}^{T}\right)\left(\ssvec{F}_{v}\cdot\svec{n}\right)\diff S + \left\langle \cssvec{F}_{v},\svec{\nabla}_{\xi}\stvec{W}\right\rangle_{E,N},\label{eq:stability:replace-test:g}\\
	\left\langle \mathcal J\mu,C_{t}\right\rangle_{E,N}=&\left\langle \mathcal JF_0',C_{t}\right\rangle_{E,N}-\int_{\partial e,N}\frac{3}{2}\sigma\varepsilon C_{t}\svec{{G}}_{c}^{\star}\cdot\svec{n}\diff S+\left\langle\frac{3}{2}\sigma\varepsilon\svec{\tilde{G}}_{c},\svec{\nabla}_{\xi}C_{t}\right\rangle_{E,N}, \label{eq:stability:replace-test:mu}\\
\frac{3}{2}\sigma\varepsilon\left\langle \mathcal J\svec{G}_{c,t},\svec{G}_{c}\right\rangle_{E,N}=&\int_{\partial e,N}\frac{3}{2}\sigma\varepsilon\left(C^{\star}_{t}-C_{t}\right)\svec{G}_{c}\cdot\svec{n}\diff S + \left\langle \frac{3}{2}\sigma\varepsilon\svec{\tilde 
G}_{c},\svec{\nabla}_{\xi}C_{t}\right\rangle_{E,N}.\label{eq:stability:replace-test:c}
	\end{align}
\end{subequations}
Next, we replace the last term in \eqref{eq:stability:replace-test:q} by the last term 
in \eqref{eq:stability:replace-test:g}, and the last term in \eqref{eq:stability:replace-test:mu} 
by the last term in \eqref{eq:stability:replace-test:c},
\begin{subequations}\label{eq:stability:combine-2}
	\begin{align}
\left\langle \mathcal J\smat{M}\stvec{Q}_{t},\stvec{W}\right\rangle_{E,N} &+\int_{\partial e,N}\stvec{W}^{T}\left(\ssvec{F}_{e}^{\star}-\ssvec{F}_{e}+\sum_{m=1}^{5}\left(\ssvecg{\Phi}_{m}W_m\right)^{\diamondsuit}-\ssvecg{\Phi}_{m}W_{m}\right)\cdot\svec{n}\diff S \nonumber\\
&+\left\langle \stvec{W},\svec{\nabla}_{\xi}\cdot\cssvec{F}_{e}\right\rangle_{E,N}+\sum_{m=1}^{5}\left\langle \stvec{W},\cssvecg{\Phi}_{m}\cdot\svec{\nabla}_{\xi}W_{m}\right\rangle_{E,N} \nonumber\\
=&\int_{\partial e,N}\left(\stvec{W}^{T}\ssvec{F}_{v}^{\star}+\stvec{W}^{\star,T}\ssvec{F}_{v}-\stvec{W}^{T}\ssvec{F}_{v}\right)\cdot\svec{n}\diff S- \left\langle \mathcal J\ssvec{G},\ssvec{F}_{v}\right\rangle_{E,N},\label{eq:stability:combine-2:q}\\
	\left\langle \mathcal J\mu,C_{t}\right\rangle_{E,N}=&\left\langle \mathcal JF_0',C_{t}\right\rangle_{E,N}+\frac{3}{2}\sigma\varepsilon\left\langle \mathcal J\svec{G}_{c,t},\svec{G}_{c}\right\rangle_{E,N} \nonumber \\
	&-\int_{\partial e,N}\frac{3}{2}\sigma\varepsilon \left(C_{t}\svec{{G}}_{c}^{\star}+C_{t}^{\star}\svec{G}_{c} - C_{t}\svec{G}_{c}\right)\cdot\svec{n}\diff S,\label{eq:stability:combine-2:mu+c}
	\end{align}
\end{subequations}
to combine four equations into two. Under the assumption of exactness in time,
we can use the chain rule on \eqref{eq:stability:combine-2:mu+c},
\begin{equation}
\begin{split}
  \left\langle \mathcal JF_0',C_{t}\right\rangle_{E,N}+\frac{3}{2}\sigma\varepsilon\left\langle \mathcal J\svec{G}_{c,t},\svec{G}_{c}\right\rangle_{E,N}  
  =&\left\langle \mathcal J\left(F_{0}+\frac{3}{4}\sigma\varepsilon|\svec{G}_{c}|^{2}\right)_{t},1\right\rangle_{E,N}  
  \\
  =&\left\langle \mathcal J \mathcal 
  F_{t},1\right\rangle_{E,N},
  \end{split}
\end{equation}
to obtain the time derivative of the discrete approximation of the free--energy \eqref{eq:governing:freeenergy-def}.

We complete the analysis by examining each term in \eqref{eq:stability:combine-2} separately to get the 
semi--discrete entropy law. The rest of the section computes the terms as follows: the time derivative coefficients in 
Sec.~\ref{subsec:stability:time-derivative}, inviscid volume integrals in 
Sec.~\ref{subsec:stability:Inviscid-volume}, viscous volume integrals in 
Sec.~\ref{subsec:stability:Viscous-volume}, and lastly boundary integrals in 
Sec.~\ref{subsec:Stability:BoundaryTerms}. The discrete entropy law is 
completed in Sec.~\ref{subsec:stability:Summary}.

\subsection{Time derivative coefficients}\label{subsec:stability:time-derivative}

We first study the discrete inner product containing time derivatives  found in \eqref{eq:stability:combine-2:q}.
We expand the inner product,
\begin{equation}
  \left\langle \mathcal J\smat{M}\stvec{Q}_{t},\stvec{W}\right\rangle_{E,N} =
  \left\langle\mathcal JC_{t},\mu\right\rangle_{E,N} + \left\langle\mathcal 
  J\left(\sqrt{\rho}\svec{U}\right)_{t},\sqrt{\rho}\svec{U}\right\rangle_{E,N} + 
  \frac{1}{\rho_0 c_0^2}\left\langle \mathcal J P_{t},P\right\rangle_{E,N},
  \label{eq:stability:time-derivative-contraction-prep}
\end{equation}
into which we insert the discrete approximation of $\left\langle \mathcal J C_{t},\mu\right\rangle_{E,N}$ 
found in \eqref{eq:stability:combine-2:mu+c}, and apply the chain rule in time for the kinetic and artificial compressibility
terms,
\begin{equation}
\begin{split}
  \left\langle \mathcal J\smat{M}\stvec{Q}_{t},\stvec{W}\right\rangle_{E,N} =&\left\langle \mathcal J\mathcal F_{t},1\right\rangle_{E,N}-\int_{\partial e,N}\frac{3}{2}\sigma\varepsilon \left(C_{t}\svec{{G}}_{c}^{\star}+C_{t}^{\star}\svec{G}_{c} - C_{t}\svec{G}_{c}\right)\cdot\svec{n}\diff 
  S\\
  &+\left\langle\mathcal 
  J\left(\frac{1}{2}\rho V_{tot}^{2}\right)_{t},1\right\rangle_{E,N} + 
\left\langle \mathcal 
 J\left(\frac{P^2}{2\rho_0 c_0^2}\right)_{t},1\right\rangle_{E,N}\\
  =&\left\langle \mathcal J\left(\mathcal F + \mathcal K + \mathcal E_{\mathrm{AC}}\right)_{t},1\right\rangle_{E,N}\\
  &-\int_{\partial e,N}\frac{3}{2}\sigma\varepsilon \left(C_{t}\svec{{G}}_{c}^{\star}+C_{t}^{\star}\svec{G}_{c} - C_{t}\svec{G}_{c}\right)\cdot\svec{n}\diff 
  S \\
    =&\left\langle \mathcal J \mathcal E_{t},1\right\rangle_{E,N}-\int_{\partial e,N}\frac{3}{2}\sigma\varepsilon \left(C_{t}\svec{{G}}_{c}^{\star}+C_{t}^{\star}\svec{G}_{c} - C_{t}\svec{G}_{c}\right)\cdot\svec{n}\diff 
  S.
  \end{split}
  \label{eq:stability:time-derivative-contraction}
\end{equation}
Eq.~\eqref{eq:stability:time-derivative-contraction} is the discrete version of \eqref{eq:continuous:entropy-time-contraction}. As a result, we get the discrete version of \eqref{eq:continuous:entropy-time-contraction} for the entropy,
\begin{equation}
  \mathcal E = \mathcal F + \mathcal K + \mathcal E_{\mathrm AC} = F_0 + \frac{3}{4}\sigma\varepsilon|\svec{G}_{c}|^{2} + 
  \frac{1}{2}\rho V_{tot}^{2} + \frac{P^2}{2\rho_0c_0^2},
\end{equation}
and for the time entropy flux $\svec{F}^{\mathcal E}_{t}$, which appears as the argument of the surface integral in \eqref{eq:stability:time-derivative-contraction}.

\subsection{Inviscid volume terms}\label{subsec:stability:Inviscid-volume}

We now show that the contraction of the inviscid volume terms into a boundary entropy flux holds 
discretely as in the 
continuous analysis \eqref{eq:continuous:inviscid-fluxes-contraction}, i.e., 
\begin{equation}
  \left\langle \stvec{W},\svec{\nabla}_{\xi}\cdot\cssvec{F}_{e}\right\rangle_{E,N}+\sum_{m=1}^{5}\left\langle \stvec{W},\cssvecg{\Phi}_{m}\cdot\svec{\nabla}_{\xi}W_{m}\right\rangle_{E,N}
  = \int_{\partial e,N}\svec{F}^{\mathcal E}\cdot\svec{n}\diff S.
  \label{eq:stability:volume-contraction}
\end{equation}
Two properties were used in the continuous analysis: the 
integration by parts (which holds discretely by the discrete Gauss law 
\eqref{eq:dg:discreteGaussLaw}), and the relationship between inviscid fluxes 
and non--conservative terms \eqref{eq:continuous:entropy-condition}, which
also holds for contravariant components,
\begin{equation}
\stvec{e}_{m}^{T}\cssvec{F}_{e}= \stvec{e}_{m}^{T}\left(\mathcal 
  M^{T}\ssvec{F}_{e}\right)=  \tens{M}^{T}\left(\stvec{e}_{m}^{T}\ssvec{F}_{e}\right)  = 
  \tens{M}^{T}\left(\stvec{W}^{T}\ssvecg{\Phi}_{m}\right)=\stvec{W}^{T}\left(\mathcal 
  M^{T}\ssvecg{\Phi}_{m}\right)=\stvec{W}^{T}\cssvecg{\Phi}_{m}.
  \label{eq:stability:volume-relation}
\end{equation}
In \eqref{eq:stability:volume-relation}, space and state multiplications commute 
since the matrix \eqref{eq:metrics:metrics-matrix} is built from identity matrices state--wise. 

Eq.~\eqref{eq:stability:volume-relation} 
together with the discrete Gauss law \eqref{eq:dg:discreteGaussLaw} is sufficient to 
prove \eqref{eq:stability:volume-contraction}. We use the discrete Gauss law for 
the first term in \eqref{eq:stability:volume-contraction}, and explicitly write  the 
scalar product of the resulting integral approximation as the sum of the components,
\begin{equation}
\begin{split}
  \left\langle \stvec{W},\svec{\nabla}_{\xi}\cdot\cssvec{F}_{e}\right\rangle_{E,N} 
  =& \int_{\partial e,N}\stvec{W}^{T}\ssvec{F}_{e}\cdot\svec{n}\diff S - 
  \left\langle\cssvec{F}_{e},\svec{\nabla}_{\xi}\stvec{W}\right\rangle_{E,N} \\
  =& \int_{\partial e,N}\stvec{W}^{T}\ssvec{F}_{e}\cdot\svec{n}\diff S - 
  \sum_{m=1}^{5}\left\langle\stvec{e}_m^{T}\cssvec{F}_{e},\svec{\nabla}_{\xi}W_{m}\right\rangle_{E,N}\\
  =& \int_{\partial e,N}\svec{F}^{\mathcal E}\cdot\svec{n}\diff S - 
  \sum_{m=1}^{5}\left\langle\stvec{W}^{T}\cssvec{\Phi}_{m},\svec{\nabla}_{\xi}W_{m}\right\rangle_{E,N}.
  \end{split}
  \label{eq:Stability:volumeProof}
\end{equation}
In the last identity, we used \eqref{eq:continuous:entropy-flux} and 
\eqref{eq:stability:volume-relation}. Thus, when we move the non--conservative terms to the left hand side of \eqref{eq:Stability:volumeProof}, we prove 
\eqref{eq:stability:volume-contraction}.

\subsection{Viscous volume terms}\label{subsec:stability:Viscous-volume}

Analogously to \eqref{eq:continuous:viscous-volume-terms}, viscous volume terms are dissipative if %\eqref{eq:stability:viscous-volume}
\begin{equation}
\left\langle \mathcal J\ssvec{G},\ssvec{F}_{v}\left(\ssvec{G}\right)\right\rangle_{E,N} 
\geqslant 0.
\label{eq:stability:viscous-volume}
\end{equation}
Thus, we closely follow the continuous steps,
\begin{equation}
\begin{split}
  \left\langle \mathcal J\ssvec{G},\ssvec{F}_{v}\right\rangle_{E,N} &= 
  M_0\left\langle \mathcal J\svec{G}_{\mu},\svec{G}_{\mu}\right\rangle_{E,N} + 
  \left\langle \mathcal J\svec{G}_{\svec{u}},2\eta\mathcal 
  S\right\rangle_{E,N}\\
  &=M_0\left\langle \mathcal J\svec{G}_{\mu},\svec{G}_{\mu}\right\rangle_{E,N}  
  + \left\langle \mathcal J\left(\frac{\svec{G}_{\svec{u}}+\svec{G}_{\svec{u}}^{T}}{2}\right),2\eta\mathcal 
  S\right\rangle_{E,N}\\
&=\left\langle \mathcal J\left(M_0|\svec{G}_{\mu}|^2+2\eta\mathcal 
  S:\mathcal S\right),1\right\rangle_{E,N} \geqslant 0,
  \end{split}
  \label{eq:stability:viscous-volume-expression}
\end{equation}
where $\mathcal S = \tens{S}\left(\svec{G}_{\svec{u}}\right)$ is the strain tensor 
\eqref{eq:governing:strain-tensor}
evaluated with the discrete gradients $\svec{G}_{\svec{u}}$. As a result, we 
find that discrete viscous volume terms are dissipative.

\subsection{Boundary terms}\label{subsec:Stability:BoundaryTerms}

Boundary term stability only makes sense when considering all elements 
in the domain. Thus, we update \eqref{eq:stability:combine-2} with 
\eqref{eq:stability:time-derivative-contraction}, 
\eqref{eq:stability:volume-contraction}, and 
\eqref{eq:stability:viscous-volume-expression}, and sum over all elements, creating
\begin{equation}
  \frac{\diff \bar{\mathcal E}}{\diff t} + \mathrm{IBT} + \mathrm{PBT} = -\sum_{e}\left\langle \mathcal J\left(M_0|\svec{G}_{\mu}|^2+2\eta\mathcal 
  S:\mathcal S\right),1\right\rangle_{E,N},
  \label{eq:stability:entropy-all-elements}
\end{equation}
where $\bar{\mathcal E}$ is the total entropy,
\begin{equation}
  \bar{\mathcal E} = \sum_{e}\left\langle \mathcal J\mathcal 
  E,1\right\rangle_{E,N},
\end{equation}
IBT is the contribution from interior faces,
\begin{equation}
\begin{split}
  \mathrm{IBT} =& \mathrm{IBT}_{e} + \mathrm{IBT}_{v} + \mathrm{IBT}_{ch} \\
  =& -\sum_{\interiorfaces}\int_{f,N}\left(\jump{\stvec{W}^{T}}\ssvec{F}_{e}^{\star} +\sum_{m=1}^{5}\jump{\stvec{W}^{T}\left(\ssvec{\Phi}_{m}W_{m}\right)^{\diamondsuit}}-\jump{\stvec{W}^{T}\ssvec{\Phi}_{m}W_{m}}\right)\cdot\svec{n}_{L}\mathrm{d}S\\
  & +\sum_{\interiorfaces}\int_{f,N}\left(\jump{\stvec{W}^{T}}\ssvec{F}_{v}^{\star} + \jump{\stvec{W}^{T}\ssvec{F}_{v}}- 
  \stvec{W}^{\star,T}\jump{\svec{F}_{v}}\right)\cdot\svec{n}_{L}\mathrm{d}S\\
  &+\frac{3}{2}\sigma\varepsilon\sum_{\interiorfaces}\int_{f,N}\left(\jump{C_{t}}\svec{G}_{c}^{\star}+\jump{\svec{G}_{c}}C_{t}^{\star}-\jump{C_{t}\svec{G}_{c}}\right)\cdot\svec{n}_{L}\mathrm{d}S,\\
  \end{split}
  \label{eq:stability:IBT}
\end{equation}
and PBT is the contribution from physical boundary faces,
\begin{equation}
  \begin{split}
    \mathrm{PBT} =
    & \mathrm{PBT}_{e} + \mathrm{PBT}_{v}+ \mathrm{PBT}_{ch}\\
    =& \sum_{\boundaryfaces}\int_{f,N}\left(\stvec{W}^{T}\ssvec{F}_{e}^{\star} 
    +\sum_{m=1}^{5}\stvec{W}^{T}\left(\left(\ssvec{\Phi}_{m}W_{m}\right)^{\diamondsuit}-\ssvec{\Phi}_{m}W_{m}\right)\right)\cdot\svec{n}\mathrm{d}S\\
    &-\sum_{\boundaryfaces}\int_{f,N}\left(\stvec{W}^{T}\left(\ssvec{F}_{v}^{\star}-\ssvec{F}_{v}\right)+\stvec{W}^{\star,T}\ssvec{F}_{v}\right)\cdot\svec{n}\mathrm{d}S\\
    &-\frac{3}{2}\sigma\varepsilon\sum_{\boundaryfaces}\int_{f,N}\left(C_t\svec{G}_{c}^\star+C^\star_t \svec{G}_{c}-C_t \svec{G}_{c}\right)\cdot\svec{n}\diff 
    S.
  \end{split}
  \label{eq:stability:PBT}
\end{equation}
In the inviscid fluxes, the entropy flux definition $\svec{F}^{\mathcal E}=\stvec{W}^{T}\ssvec{F}_{e}$ was taken into account to reduce the number of terms. 
Moreover, in the interior boundary terms, we used the 
left element as the reference (hence the $\svec{n}_{L}$ vector).

We now analyze the contributions of the interior and boundary terms in 
\eqref{eq:stability:entropy-all-elements} separately.

\subsubsection{Inviscid interior boundary terms: entropy conserving scheme with central fluxes}\label{subsubsec:stability:boundary:interior-inviscid-central}

We address the stability of the inviscid interior boundary terms when using the central 
fluxes for the conservative numerical \eqref{eq:riemann:central-cons} and non--conservative 
 diamond fluxes \eqref{eq:riemann:central-noncons}. Inserting the central fluxes into the inviscid $\mathrm{IBT}_{e}$ \eqref{eq:stability:IBT} terms,
\begin{equation}
\begin{split}
\mathrm{IBT}_{e}=&  -\sum_{\interiorfaces}\int_{f,N}\left(\jump{\stvec{W}^{T}}\ssvec{F}_{e}^{\star} +\sum_{m=1}^{5}\jump{\stvec{W}^{T}\left(\ssvec{\Phi}_{m}W_{m}\right)^{\diamondsuit}}-\jump{\stvec{W}^{T}\ssvec{\Phi}_{m}W_{m}}\right)\cdot\svec{n}_{L}\mathrm{d}S\\
  =&  -\sum_{\interiorfaces}\int_{f,N}\left(\jump{\stvec{W}^{T}}\aver{\ssvec{F}_{e}} +\sum_{m=1}^{5}\jump{\stvec{W}^{T}\ssvec{\Phi}_{m}}\aver{W_{m}}-\jump{\stvec{W}^{T}\ssvec{\Phi}_{m}W_{m}}\right)\cdot\svec{n}_{L}\mathrm{d}S\\
    =&  -\sum_{\interiorfaces}\int_{f,N}\sum_{m=1}^{5}\left(\jump{{W}_{m}}\aver{\stvec{e}^{T}_{m}\ssvec{F}_{e}} -\jump{W_{m}}\aver{\stvec{W}^{T}\ssvec{\Phi}_{m}}\right)\cdot\svec{n}_{L}\mathrm{d}S=0.\\
  \end{split}
  \label{eq:stability_IBTe}
\end{equation}
Thus we find that inviscid interior boundaries exactly satisfy $\mathrm{IBT}_{e}=0$. In the second equality, we explicitly wrote the first product as the sum 
of the product of its components, and in the last equality, we used the algebraic relationship satisfied by jump and average operators,
\begin{equation}
  \jump{uv} = \aver{u}\jump{v} + \jump{u}\aver{v}
\end{equation}
for the second and third products. Additionally,  both terms in the last sum of \eqref{eq:stability_IBTe} are identical because of the condition \eqref{eq:continuous:entropy-condition}. 
We conclude that inviscid interior boundary contribution to the entropy 
equation is identically zero, i.e. the interface approximation is entropy conserving.

\subsubsection{Inviscid interior boundary terms: entropy stable scheme with the exact Riemann solver}\label{subsubsec:stability:boundary:interior-inviscid-ers}

The exact Riemann solver was proven to be stable in \cite{2019:Manzanero-iNS} 
for the incompressible Navier--Stokes equations. Now we extend the proof to 
the iNS/CH system, which is done by a  consistent choice of the diamond fluxes.
The proof rests on the exact Riemann solver 
star solution satisfying
\begin{equation}
\begin{split}
  \Delta_{e}^{\mathrm{ERS}} =& \frac{1}{2}\rho^{\star}\jump{V_{tot}^2}U_n^{\star} - \rho^{\star}U_n^{\star}\left(U_n^{\star}\jump{U_n}+V_{t1}^{\star}\jump{V_{t1}}+V_{t2}^{\star}\jump{V_{t2}}\right)\\
  &- P^{\star}\jump{U_n}- \jump{P}U_n^{\star} +\jump{P 
  U_n}\geqslant 0.
  \end{split}
  \label{eq:stability:ERS-iNS}
\end{equation}
Eq.~\eqref{eq:stability:ERS-iNS} implies that the ERS is dissipative for the incompressible NSE. We use the same 
star region solution for the inviscid numerical fluxes \eqref{eq:riemann:ers-cons} 
and the choice for the diamond fluxes given in 
\eqref{eq:riemann:ers-noncons}. The result is that the argument of the inviscid interior boundary 
terms is
\begin{equation}
\begin{split}
  \Delta_{\mathrm{e}}^{\mathrm{IBT}} =&\left(-\jump{\stvec{W}^{T}}\ssvec{F}_{e}^{\star}-\sum_{m=1}^{5}\jump{\stvec{W}^{T}\left(\ssvec{\Phi}_{m}W_{m}\right)^{\diamondsuit}}
 +\jump{\stvec{W}^{T}\ssvec{\Phi}_{m}W_{m}}\right)\cdot\svec{n}_{L} \\
 =&-\jump{\stvec{W}^{T}}\smat{T}^{T}\stvec{F}_{e}^{\star}-\jump{\stvec{W}^{T}\smat{T}^{T}\stvec{F}_{\Phi W}^{\mathrm{ERS}}}
 +\jump{\stvec{W}^{T}\smat{T}^{T}\stvec{F}_{\Phi W}} \\
  =&-\jump{\stvec{W}_{n}^{T}}\stvec{F}_{e}^{\star}-\jump{\stvec{W}_{n}^{T}\left(\stvec{F}_{\Phi W}^{\mathrm{ERS}}-\stvec{F}_{\Phi 
  W}\right)}.
 \end{split}
 \label{eq:stability:boundary:interior-inviscid-ibt}
\end{equation}
We write the first term in the last line of \eqref{eq:stability:boundary:interior-inviscid-ibt} as
\begin{equation}
  -\jump{\stvec{W}_{n}^{T}}\stvec{F}_{e}^{\star} = -\jump{\mu}\aver{CU_n}
  -\frac{1}{2}\rho^{\star}U_n^{\star}\left(U_n^{\star}\jump{U_n}+V_{t1}^{\star}\jump{V_{t1}}+V_{t2}^{\star}\jump{V_{t2}}\right)-P^{\star}\jump{U_n},
\end{equation}
and the second term as
\begin{equation}
\begin{split}
  -\jump{\stvec{W}_{n}^{T}\left(\stvec{F}_{\Phi W}^{\mathrm{ERS}}-\stvec{F}_{\Phi 
  W}\right)} =& 
  -\frac{1}{2}\rho^{\star}U_n^{\star,2}\jump{U_n}+\frac{1}{2}\rho^{\star}U_n^{\star}\jump{U_n^2}-\aver{\mu}\jump{CU_n}+\jump{CU_n\mu}\\
  &-\frac{1}{2}\rho^{\star}U_n^{\star}V_{t1}^{\star}\jump{V_{t1}}+\frac{1}{2}\rho^{\star}U_n^{\star}\jump{V_{t1}^2}\\
  &-\frac{1}{2}\rho^{\star}U_n^{\star}V_{t2}^{\star}\jump{V_{t2}}+\frac{1}{2}\rho^{\star}U_n^{\star}\jump{V_{t2}^2}\\
  &-U_n^{\star}\jump{P} + \jump{PU_n} \\
  =&\frac{1}{2}\rho^{\star}U_n^{\star}\jump{V_{tot}^2} +\jump{\mu}\aver{CU_n}\\
  &-\frac{1}{2}\rho^{\star}U_n^{\star}\left(U_n^{\star}\jump{U_n}+V_{t1}^{\star}\jump{V_{t1}}+V_{t2}^{\star}\jump{V_{t2}}\right)\\
  &-U_n^{\star}\jump{P} + \jump{PU_n},
  \end{split}
\end{equation}
which when inserted into \eqref{eq:stability:boundary:interior-inviscid-ibt} yields,
\begin{equation}
\begin{split}
    \Delta_{\mathrm{e}}^{\mathrm{IBT}} =& 
    \frac{1}{2}\rho^{\star}U_n^{\star}\jump{V_{tot}}^{2} - 
    \rho^{\star}U_n^{\star}\left(U_n^{\star}\jump{U_n}+V_{t1}^{\star}\jump{V_{t1}}+V_{t2}^{\star}\jump{V_{t2}}\right)
    \\
    &-P^{\star}\jump{U_n}-U_n^{\star}\jump{P}+\jump{PU_n} = \Delta_{e}^{\mathrm{ERS}}\geqslant 0.
    \end{split}
\end{equation}
The non--negativity of $\Delta_{\mathrm{e}}^{\mathrm{IBT}}$ confirms that the incompressible NSE exact Riemann solver \eqref{eq:riemann:ers-cons}, 
with the associated diamond flux choice \eqref{eq:riemann:ers-noncons} introduced in this work produces 
an entropy stable interface approximation for the iNS/CH system.

\subsubsection{Viscous and Cahn--Hilliard interior terms: Bassi--Rebay 1 method}\label{subsubsec:stability:boundary:interior-viscous-ch}

The use of the standard BR1 method identically cancels the boundary integrals (see \cite{2017:Gassner} for the compressible NSE,
\cite{2019:Manzanero-iNS} for the incompressible NSE, \cite{2018:Bohm} for the MHD equations, and \cite{2019:Manzanero-CH} for the Cahn--Hilliard 
equation). Then, we prove that the extra interface stabilization terms are dissipative. In this work, we obtain the same integral approximation as in \cite{2019:Manzanero-iNS} 
for the incompressible NSE, and as in \cite{2019:Manzanero-CH} for the Cahn--Hilliard part,
\begin{equation}
\begin{split}
  \mathrm{IBT}_{v} =& \sum_{\interiorfaces}\int_{f,N}\left(\jump{\stvec{W}^{T}}\ssvec{F}_{v}^{\star} - 
  \jump{\stvec{W}^{T}\ssvec{F}_{v}}+
  \stvec{W}^{\star,T}\jump{\ssvec{F}_{v}}\right)\cdot\svec{n}_{L}\mathrm{d}S\\
  =& \sum_{\interiorfaces}\int_{f,N}\left(\jump{\stvec{W}^{T}}\aver{\ssvec{F}_{v}}+\beta\jump{\mu}^2\svec{n}_{L}- \jump{\stvec{W}^{T}\ssvec{F}_{v}}+ 
  \aver{\stvec{W}}^{T}\jump{\ssvec{F}_{v}}\right)\cdot\svec{n}_{L}\mathrm{d}S\\
  =& \sum_{\interiorfaces}\int_{f,N}\beta\jump{\mu}^2\mathrm{d}S\geqslant 0,
  \end{split}
\end{equation}
whereas for the Cahn--Hilliard terms,
\begin{equation}
\begin{split}
  \mathrm{IBT}_{ch}=& 
  \frac{3}{2}\sigma\varepsilon\sum_{\interiorfaces}\int_{f,N}\left(\jump{C_{t}}\svec{G}_{c}^{\star}+\jump{\svec{G}_{c}}C_{t}^{\star}-\jump{C_{t}\svec{G}_{c}}\right)\cdot\svec{n}_{L}\mathrm{d}S\\
=&  \frac{3}{2}\sigma\varepsilon\sum_{\interiorfaces}\int_{f,N}\left(\jump{C_{t}}\aver{\svec{G}_{c}}+\jump{\svec{G}_{c}}\aver{C_{t}}+\beta\jump{C}\jump{C_{t}}\svec{n}_{L}-\jump{C_{t}\svec{G}_{c}}\right)\cdot\svec{n}_{L}\mathrm{d}S\\  
=&\frac{3}{4}\sigma\varepsilon\frac{\diff}{\diff t}\sum_{\interiorfaces}\int_{f,N}\beta\jump{C}^2\mathrm{d}S,\\
\end{split}
\label{eq:stability:CH-IBT}
\end{equation}
which says that viscous and Cahn--Hilliard interior boundary contributions do not contribute to the entropy 
balance if $\beta=0$, since the integrals are identically zero. This means that interface stabilization is not a requirement for stability. When $\beta>0$, the integrals in $\mathrm{IBT}_{v}$ are positive (i.e. dissipative), and the integral in $\mathrm{IBT}_{ch}$, which is not positive per se because of the time derivative, will be added to the entropy, since the integral argument is positive.

\subsubsection{Interior boundary terms: Summary}

The contribution to the interior boundary terms from viscous and Cahn--Hilliard terms are 
identically zero when using the BR1 method without interface stabilization ($\beta=0$), and dissipative otherwise ($\beta>0$). As for the inviscid terms, interior 
boundary terms are conservative with central fluxes ($\mathrm{IBT}_{e}=0$), and dissipative with the 
exact Riemann solver ($\mathrm{IBT}_{e}\geqslant 0$).

\subsubsection{Physical boundary terms: free-- and no--slip walls}\label{subsubsec:stability:physical-boundary-wall}

We now address the stability of physical wall boundary terms 
\eqref{eq:stability:PBT}. For inviscid fluxes $\mathrm{PBT}_{e}$, we use 
the entropy conserving approach with central fluxes 
\eqref{eq:BCs:inviscid-central},
\begin{equation}
\begin{split}
  \mathrm{PBT}_{e} &=  \sum_{\boundaryfaces}\int_{f,N}\left(\stvec{W}^{T}\ssvec{F}_{e}^{\star}
    +\sum_{m=1}^{5}\stvec{W}^{T}\left(\left(\ssvec{\Phi}_{m}W_{m}\right)^{\diamondsuit}-\ssvec{\Phi}_{m}W_{m}\right)\right)\cdot\svec{n}\mathrm{d}S\\
    &=\sum_{\boundaryfaces}\int_{f,N}\left(U_n\left(\frac{1}{2}\rho U_n^2 + P\right)+U_n\left(\mu C-\frac{1}{2}\rho U_n -\mu C\right)-PU_n\right)\diff 
    S=0,
  \end{split}
\end{equation}
or the entropy stable counterpart with the exact Riemann solver 
\eqref{eq:BCs:inviscid-ers},
\begin{equation}
\begin{split}
  \mathrm{PBT}_{e} &=  \sum_{\boundaryfaces}\int_{f,N}\left(\stvec{W}^{T}\ssvec{F}_{e}^{\star}
    +\sum_{m=1}^{5}\stvec{W}^{T}\left(\left(\ssvec{\Phi}_{m}W_{m}\right)^{\diamondsuit}-\ssvec{\Phi}_{m}W_{m}\right)\right)\cdot\svec{n}\mathrm{d}S\\
    &=\sum_{\boundaryfaces}\int_{f,N}\left(U_n\left(\frac{1}{2}\rho \lambda^{+} U_n + P\right)-PU_n\right)\diff 
    S=\sum_{\boundaryfaces}\int_{f,N}\left(\frac{1}{2}\rho \lambda^{+} U_n^2\right)\diff 
    S\geqslant 0.
  \end{split}
\end{equation}
Either way, the inviscid physical boundary contribution is stable. Furthermore, for the
viscous terms, the free--slip wall boundary condition \eqref{eq:BCs:viscous-free-slip},
\begin{equation}
  \begin{split}
    \mathrm{PBT}_{v}=&-\sum_{\boundaryfaces}\int_{f,N}\left(\stvec{W}^{T}\left(\ssvec{F}_{v}^{\star}-\ssvec{F}_{v}\right)+\stvec{W}^{\star,T}\ssvec{F}_{v}\right)\cdot\svec{n}\mathrm{d}S\\    
    &-\sum_{\boundaryfaces}\int_{f,N}\left(\stvec{W}^{T}\left({0}-\ssvec{F}_{v}\right)+\stvec{W}^{T}\ssvec{F}_{v}\right)\cdot\svec{n}\mathrm{d}S=0,    
  \end{split}
\end{equation}
and for the no--slip wall boundary condition \eqref{eq:BCs:viscous-no-slip},
\begin{equation}
  \begin{split}
    \mathrm{PBT}_{v}=&-\sum_{\boundaryfaces}\int_{f,N}\left(\stvec{W}^{T}\left(\ssvec{F}_{v}^{\star}-\ssvec{F}_{v}\right)+\stvec{W}^{\star,T}\ssvec{F}_{v}\right)\cdot\svec{n}\mathrm{d}S\\    
    &-\sum_{\boundaryfaces}\int_{f,N}\left(-\mu\svec{G}_{\mu}+\mu\svec{G}_{\mu}+\svec{U}\mathcal S-\svec{U}\mathcal 
    S\right)\cdot\svec{n}\mathrm{d}S=0,
  \end{split}
\end{equation}
are dissipative. Lastly, the Cahn--Hilliard physical boundary 
condition for both free-- and no--slip walls 
\eqref{eq:BCs:viscous-free-slip}--\eqref{eq:BCs:viscous-no-slip} gives,
\begin{equation}
\begin{split}
\mathrm{PBT}_{ch}=&-\frac{3}{2}\sigma\varepsilon\sum_{\boundaryfaces}\int_{f,N}\left(C_t\svec{G}_{c}^\star+C^\star_t \svec{G}_{c}-C_t \svec{G}_{c}\right)\cdot\svec{n}\diff 
    S \\
    =&-\sum_{\boundaryfaces}\int_{f,N}\left(-C_t f_{w}'(C)+\frac{3}{2}\sigma\varepsilon\left(C_t \svec{G}_{c}-C_t \svec{G}_{c}\right)\right)\cdot\svec{n}\diff 
    S \\
    =&\frac{\diff}{\diff t}\sum_{\boundaryfaces}\int_{f,N}f_{w}(C)\diff 
    S,
\end{split}
\end{equation}
which represents the discrete version of the surface free--energy time 
derivative \eqref{eq:continuous:entropy-balance-w-bcs}.

\subsection{Semi--discrete stability: Summary}\label{subsec:stability:Summary}

We constructed a DG approximation of the iNS/CH system \eqref{eq:dg:semi-discrete-approx} whose analysis mimics the continuous stability analysis of Sec. 
\ref{subsec:Governing:iNS/CHEntropy} semi--discretely. 
As a result, we have shown that the discrete entropy satisfies a balance law,
\begin{equation}
    \frac{\diff}{\diff t}\left(\bar{\mathcal E}^{\beta}+\sum_{\boundaryfaces}\int_{f,N} f_{w}(C)\diff 
    S\right) \leqslant 
    -\sum_{e}\left\langle \mathcal J\left(M_0|\svec{G}_{\mu}|^2+2\eta\mathcal 
  S:\mathcal S\right),1\right\rangle_{E,N},
  \label{eq:stability:entropy-eq-bc-final}
\end{equation}
for both free-- and no--slip wall boundary conditions. The entropy has been extended to consider the interface stabilization in the Cahn--Hilliard equation,
\begin{equation}
\bar{\mathcal E}^{\beta} = \sum_{e}\left\langle \mathcal J\mathcal E,1\right\rangle_{E,N}+\frac{3}{4}\sigma\varepsilon\sum_{\interiorfaces}\int_{f,N}\beta\jump{C}^2\mathrm{d}S\geqslant 0.
\end{equation}
Note that if $\beta=0$ we recover $\bar{\mathcal E}^{\beta=0}=\bar{\mathcal E}=\sum_{e}\left\langle\mathcal J\mathcal E,1\right\rangle_{E,N}$, and that a continuous solution (i.e. $\jump{C}=0$) also recovers the original entropy. Otherwise, the extended entropy $\bar{\mathcal E}^{\beta}$ is also always positive, such that its definition remains valid and the scheme is entropy--stable. 
The justification to include this term in the entropy, is that the free--energy measures the fluid interfaces using the solution gradients $\frac{3}{4}\sigma\varepsilon|\svec{G}_{c}|^2$. However, if the discrete solution is allowed to be discontinuous at the inter--elements interface, it enables the creation of a fluid interface, which is not reflected in the original free--energy. The addition of this penalization, allows to account for inter--element discontinuities as part of the fluid interface, so that their associated energy is also weighted in the entropy. In practice, we have found that if the flow is under--resolved and $\beta=0$, the solution minimizes the entropy allowing sharp discontinuities at the inter--elements interface. Since that effect has not been taken into account in the entropy with $\beta=0$, the flow wrongly enables the discontinuity mechanism to minimize the entropy. This adverse effect (from the point of view of accuracy) is mitigated when $\beta> 0$ \cite{2016:Gassner,2019:Manzanero-CH}. Since now the inter--elements discontinuities belong to the entropy, their size can be kept bounded and controlled. In any case, the use of interface stabilization does not compromise the stability of the scheme.

Finally, the form of \eqref{eq:stability:entropy-eq-bc-final} depends on the choice of the inviscid numerical and diamond fluxes, and the Cahn--Hilliard interface stabilization: an equality for central fluxes and $\beta=0$, 
and an inequality for the exact Riemann solver and $\beta>0$. For either choice, the discrete entropy plus the surface free energy
remains bounded by the initial condition.

\section{Numerical experiments}\label{sec:NumericalExperiments}

We now explore the capabilities of the new DGSEM with numerical 
simulations. First, we validate the implementation with a convergence study on 
a manufactured solution. Second, we explore the robustness of the scheme by integrating in time from a 
random initial condition, where we compare an entropy--stable scheme with a Gauss point 
formulation that is not provably stable. Third, we solve a static bubble as a 
benchmark for steady--state accuracy, and a transient rising bubble. 
Finally, we present solutions for a three--dimensional pipe flow in the annular regime.

\subsection{Convergence study}\label{subsec:num:conv}

We first assess the spatial and temporal accuracy of the implementation of the approximation. We test with the same two--dimensional manufactured solution and configuration used in 
other Navier--Stokes/Cahn--Hilliard works, specifically \cite{2018:Dong},
\begin{equation}
\begin{split}
c_m(x,y;t) &= \frac{1}{2}\left(1 + \cos\left(\pi x\right)\cos\left(\pi y\right)\sin\left(t\right)\right), \\
u_m(x,y;t) &=2\sin\left(\pi x\right) \cos\left(\pi z\right)\sin\left(t\right), \\
v_m(x,y;t) &=-2\cos\left(\pi x\right)\sin\left(\pi y\right) \sin\left(t\right), \\
p_m(x,y;t) &=2\sin\left(\pi x\right)\sin\left(\pi z\right)\cos\left(t\right),
\end{split}
\label{eq:num:conv:mansol}
\end{equation}
on the domain $(x,y)\in[-1,1]^2$. The final time is $t_F=0.1$, and all physical parameters are presented in Table \ref{tab:num:conv:param}. 
{\begin{table}[h]
		\centering
		\caption{List of the parameter values used with the manufactured solution \eqref{eq:num:conv:mansol}}
		\label{tab:num:conv:param}
%		\resizebox{\textwidth}{!}{%  
			\begin{tabular}{llllllll}
				\hline
				$\rho_1$& $\rho_2$ ($\text{kg}/\text{m}^3$)   & $\eta_1$ & $\eta_2$ (Pa$\cdot$s) & $\varepsilon$ (m)& $t_{CH}$ (s)& $c_0^2$ (m/s$^2$)$^2$& $\sigma$ (N/m) \\ \hline
				1.0 & 2.0 & 1.0E-3 & 1.0E-3 & $1/\sqrt{2}$ & 1.0E3 & 1.0E3 & 6.236E-3 \\
				\hline
		\end{tabular}%}
	\end{table}}
%
%For convenience, the source term is written completely in Appendix \ref{app:conv-source}. 

We measure L$^2$ 
errors using the discrete inner product,
\begin{equation}
  \mathrm{error} = \Vert \phi - \phi_m \Vert_{\mathcal J,N} = \sqrt{\sum_{e}\left\langle\mathcal J\left(\Phi - \Phi_m\right),\left(\Phi-\Phi_m\right)\right\rangle_{E,N}}.
  \label{eq:num:L2}
\end{equation}

Regarding the configuration of the scheme, we use the exact Riemann solver, and we vary the polynomial order $N$, mesh size
(we use an equally spaced Cartesian mesh $M^2$), and time--step 
size $\Delta t$ of the BDF2 scheme ($J=2$) to evaluate the convergence rates as the resolution is increased either by increasing the polynomial order or decreasing the mesh size.
%We present the polynomial order convergence study in 
%Sec.~\ref{subsubsec:Numerical:Conv:Polynomial}, and a mesh size convergence 
%study in Sec.~\ref{subsubsec:Numerical:Conv:Mesh}
%
%
%
%
%	
%\subsubsection{Polynomial order convergence}\label{subsubsec:Numerical:Conv:Polynomial}

\begin{figure}[h]
  \centering
  \subfigure[$\Delta t=10^{-4}$]{\includegraphics[width=0.49\textwidth]{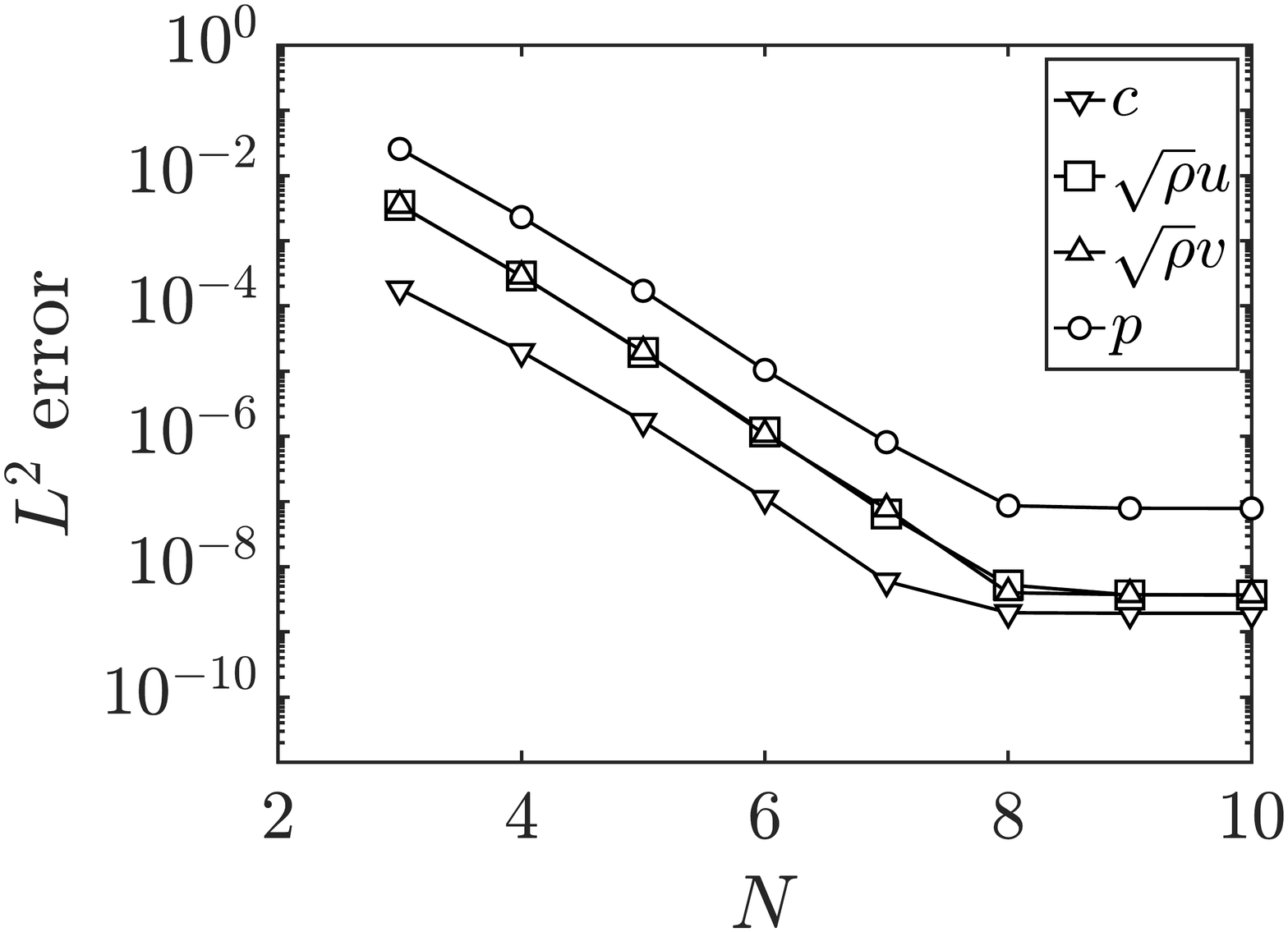}\label{fig:num:conv:pConv-4}}
  \subfigure[$\Delta t=10^{-5}$]{\includegraphics[width=0.49\textwidth]{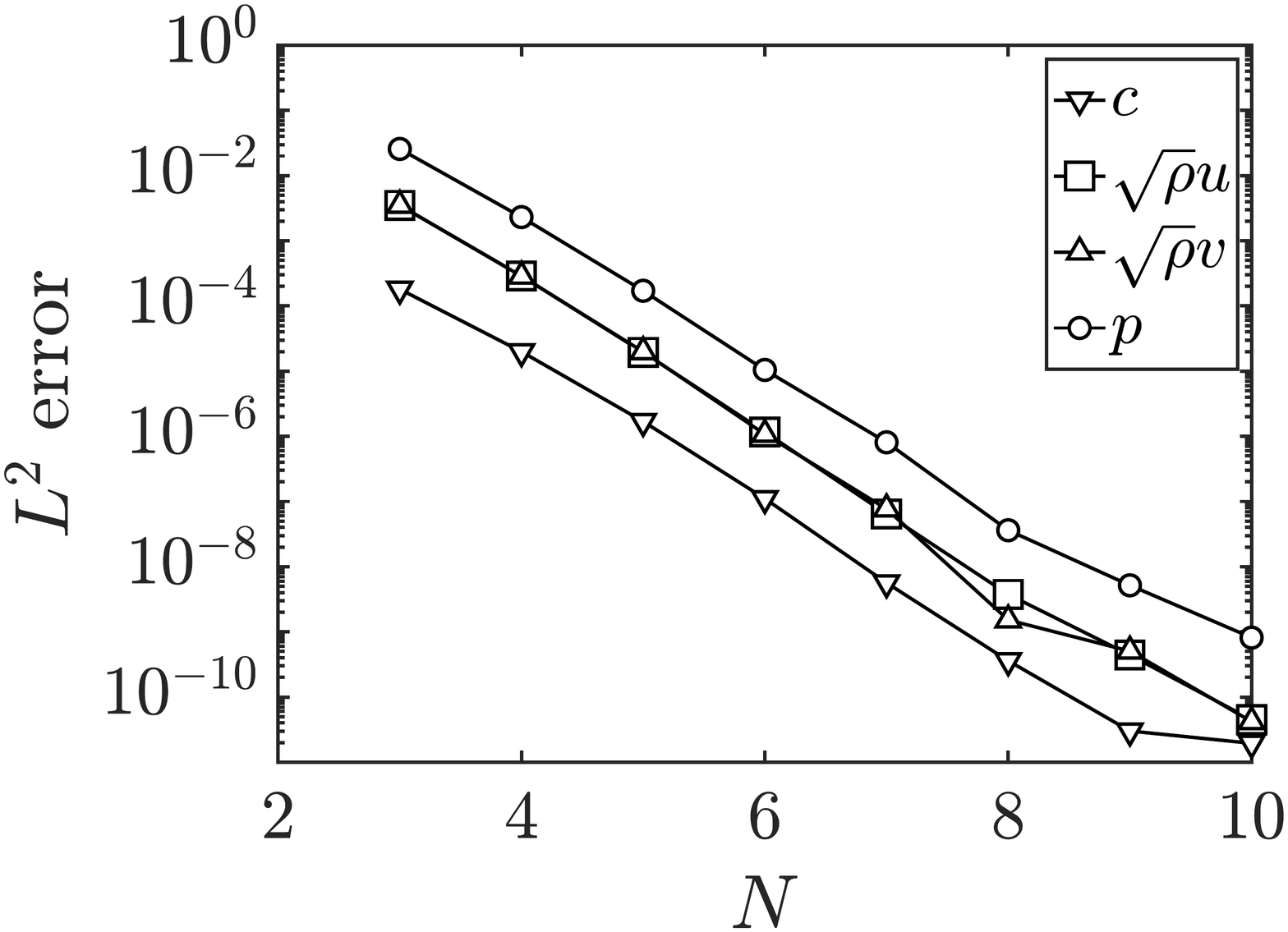}\label{fig:num:conv:pConv-5}}  
  \caption{Exponential convergence of the manufactured solution \eqref{eq:num:conv:mansol} as the polynomial order increases. 
  We represent the L$^2$ errors in concentration, $x$-- and $y$--momentum, and pressure. 
  %DAK what's the ()?
  The polynomial order ranges from 3 to 10, and we integrate in time until $t_{F}=0.1$ with two time step sizes: $\Delta t=10^{-4}$ and $10^{-5}$. All physical parameters are given in Table \ref{tab:num:conv:param}}
  \label{fig:num:conv:pConv}
\end{figure}

We first use a Cartesian $4^2$ element mesh and vary the polynomial order from $N=3$ to 
10. The L$^2$ errors \eqref{eq:num:L2} are presented in Fig.~\ref{fig:num:conv:pConv}, using $\Delta t=10^{-4}$ 
in Fig.~\ref{fig:num:conv:pConv-4} and $\Delta t=10^{-5}$ in 
Fig.~\ref{fig:num:conv:pConv-5}. 

Figs. \ref{fig:num:conv:pConv-4} and \ref{fig:num:conv:pConv-5} show typical and expected error behavior: an 
approximately linear semi-log convergence for lower polynomial orders implying exponential convergence of the error (i.e. under--resolved solution in 
space), and error stagnation for higher orders (under--resolved solution in 
time). The stagnation error threshold can be controlled with $\Delta t$, as shown 
in Fig.~\ref{fig:num:conv:pConv-5}. 
	
%\subsubsection{Mesh refinement convergence}\label{subsubsec:Numerical:Conv:Mesh}

For the next test configuration, we use polynomial orders $N=2,3,4$ and $5$, and meshes with $4^{2},6^{2},8^{2},12^{2}$ 
and $16^2$ elements for an $h-$type convergence study. The time--step is $\Delta t=5\cdot 10^{-5}$. We solve the 
manufactured solution problem \eqref{eq:num:conv:mansol}, and summarize the 
errors in Table~\ref{tab:num:conv:h-conv}. An estimate of the scheme's order 
of convergence is provided there. We note that for the momentum components, the convergence orders are 
as expected, as they are close to $N+1$. For the pressure, the order of convergence is 
somewhat smaller (closer to $N$), which was also experienced in 
\cite{2017:Bassi,2019:Manzanero-iNS}. The 
convergence is also not optimal for the concentration; for even polynomials, the convergence can 
surpass the theoretical $N+1$ value, while for odd polynomials, the convergence rate
is roughly $N$ in most of the cases. This effect was also noticed in 
\cite{2019:Manzanero-CH}, whose nature remains unclear.
{\begin{table}[h]
\centering
  \caption{Manufactured solution \eqref{eq:num:conv:mansol} convergence analysis: we use $4^3$, $8^3$, and $16^3$  meshes, and $N=2,3,4$ and 5.
The final time is $t_{F}=0.1$, and we use the IMEX BDF2 scheme
  with $\Delta t=5\cdot 10^{-5}$}
  \label{tab:num:conv:h-conv}
%\resizebox{\textwidth}{!}{%  
\begin{tabular}{llllllllll}
\hline
      & Mesh   & $c$ error & order & $\sqrt{\rho} u$ error & order & $\sqrt{\rho} v$ error & order & $p$ error & order \\ \hline
N=2 & $4^2$ & 1.35E-03 & -- & 3.39E-02 & -- & 3.39E-02 & -- & 2.22E-01 & -- \\ 
 & $6^2$ & 8.27E-04 & 1.21 & 1.24E-02 & 2.47 & 1.24E-02 & 2.47 & 9.02E-02 & 2.22 \\ 
 & $8^2$ & 4.40E-04 & 2.19 & 5.96E-03 & 2.56 & 5.95E-03 & 2.56 & 4.68E-02 & 2.28 \\ 
 & $12^2$ & 1.30E-04 & 3.01 & 2.04E-03 & 2.65 & 2.03E-03 & 2.65 & 1.78E-02 & 2.39 \\ 
 & $16^2$ & 5.35E-05 & 3.08 & 9.23E-04 & 2.75 & 9.23E-04 & 2.75 & 8.67E-03 & 2.50 \\ 
N=3 & $4^2$ & 1.81E-04 & -- & 3.49E-03 & -- & 3.50E-03 & -- & 2.57E-02 & -- \\ 
 & $6^2$ & 5.71E-05 & 2.85 & 8.28E-04 & 3.55 & 8.29E-04 & 3.55 & 6.81E-03 & 3.27 \\ 
 & $8^2$ & 2.76E-05 & 2.53 & 2.90E-04 & 3.65 & 2.91E-04 & 3.65 & 2.57E-03 & 3.39 \\ 
 & $12^2$ & 9.53E-06 & 2.62 & 6.34E-05 & 3.75 & 6.37E-05 & 3.74 & 6.17E-04 & 3.52 \\ 
 & $16^2$ & 4.36E-06 & 2.72 & 2.10E-05 & 3.83 & 2.12E-05 & 3.83 & 2.17E-04 & 3.64 \\ 
N=4 & $4^2$ & 1.99E-05 & -- & 2.89E-04 & -- & 2.87E-04 & -- & 2.31E-03 & -- \\ 
 & $6^2$ & 5.27E-06 & 3.28 & 4.40E-05 & 4.64 & 4.40E-05 & 4.62 & 3.94E-04 & 4.36 \\ 
 & $8^2$ & 7.50E-07 & 6.78 & 1.12E-05 & 4.75 & 1.12E-05 & 4.75 & 1.09E-04 & 4.47 \\ 
 & $12^2$ & 8.46E-08 & 5.38 & 1.55E-06 & 4.89 & 1.55E-06 & 4.89 & 1.70E-05 & 4.59 \\ 
 & $16^2$ & 1.94E-08 & 5.12 & 3.61E-07 & 5.05 & 3.61E-07 & 5.05 & 4.42E-06 & 4.68 \\ 
N=5 & $4^2$ & 1.69E-06 & -- & 1.94E-05 & -- & 1.98E-05 & -- & 1.71E-04 & -- \\ 
 & $6^2$ & 1.91E-07 & 5.38 & 1.96E-06 & 5.65 & 1.97E-06 & 5.70 & 1.87E-05 & 5.46 \\ 
 & $8^2$ & 4.67E-08 & 4.89 & 3.76E-07 & 5.74 & 3.77E-07 & 5.74 & 3.78E-06 & 5.57 \\ 
 & $12^2$ & 6.43E-09 & 4.89 & 3.64E-08 & 5.76 & 3.65E-08 & 5.76 & 3.79E-07 & 5.67 \\ 
 & $16^2$ & 1.63E-09 & 4.77 & 7.08E-09 & 5.69 & 7.11E-09 & 5.69 & 7.47E-08 & 5.64 \\ 
\hline
\end{tabular}%}
\end{table}
}

%In conclusion, the scheme and its implementation show the expected convergence behavior.

\subsection{Assessing robustness: random initial condition}\label{subsec:num:Random}

We also assess the robustness of the scheme by introducing a random initial condition. Because the scheme developed here is entropy stable even for severely under--resolved solutions, the numerical implementation should not crash, as long as the time--step is small enough to stay within the time integrator stability region and the positivity condition is satisfied.

 As a point of reference, and only in this section, we also consider the DG scheme using Gauss points. The standard DG method with Gauss points does not satisfy the SAT property \cite{2009:Kopriva} and the proofs developed for semi--discrete stability do not carry over. As a result, the standard scheme is not provably stable. To make the comparison even more stark,  we compare robustness of the Gauss scheme with upwind flux dissipation to the new Gauss--Lobatto method with only central fluxes.

We solve the problem on the domain $\Omega=[0,1]^3$, and on a three--dimensional $4^3$ element Cartesian mesh. We initialize the solution using uniform random numbers in $[0,1]$ for concentration, and in $[-1,1]$ for velocities and pressure. We use a high density ratio ($\rho_1/\rho_2=1000$), and Reynolds number $Re =10^6$ so that physical viscosity will have only a weak stabilizing effect. The rest of the parameters are given in Table~\ref{tab:num:random:param}.
{\begin{table}[h]
		\centering
		\caption{Physical parameters for the random initial condition test problem}
		\label{tab:num:random:param}
		%\resizebox{\textwidth}{!}{%  
		\begin{tabular}{lllllllll}
			\hline
			$\rho_1$& $\rho_2$ ($\text{kg}/\text{m}^3$)   & $\eta_1$ & $\eta_2$ (Pa$\cdot$s) & $\varepsilon$ (m)& $t_{CH}$ (s)& $c_0^2$ (m/s$^2$)$^2$& $\sigma$ (N/m) \\ \hline
			1000.0 & 1.0 & $10^{-3}$ &$ 10^{-4} $& $0.75$ & 10.0 & 1.0E2 & 1.0 \\
			\hline
		\end{tabular}%}
\end{table}}

With Gauss--Lobatto points, the discrete entropy balance \eqref{eq:stability:entropy-eq-bc-final} shows that the free--energy time derivative $\mathcal E_{t}$ should be always negative as a result of viscous and chemical potential dissipation, for either central fluxes or the ERS. Moreover, the entropy remainder, defined as
\begin{equation}
\mathcal R\left(\mathcal E\right)=\frac{\diff}{\diff t}\left(\sum_{e}\left\langle \mathcal J \mathcal E,1\right\rangle_{E,N}+\sum_{\boundaryfaces}\int_{f,N} f_{w}(C)\diff 
S\right)  +\left\langle \mathcal J\left(M_0|\svec{G}_{\mu}|^2+2\eta\mathcal 
S:\mathcal S\right),1\right\rangle_{E,N}\leqslant 0,
\label{eq:num:random:entropy-remainder}
\end{equation}
should be zero at each time step for central fluxes, and negative for the ERS. None of the statements regarding $\mathcal E_{t}$ or $\mathcal R\left(\mathcal E\right)$ being negative are guaranteed to hold for Gauss points by theory.

\begin{figure}[h]
	\centering
	\subfigure[$N=2$]{\includegraphics[width=0.32\textwidth]{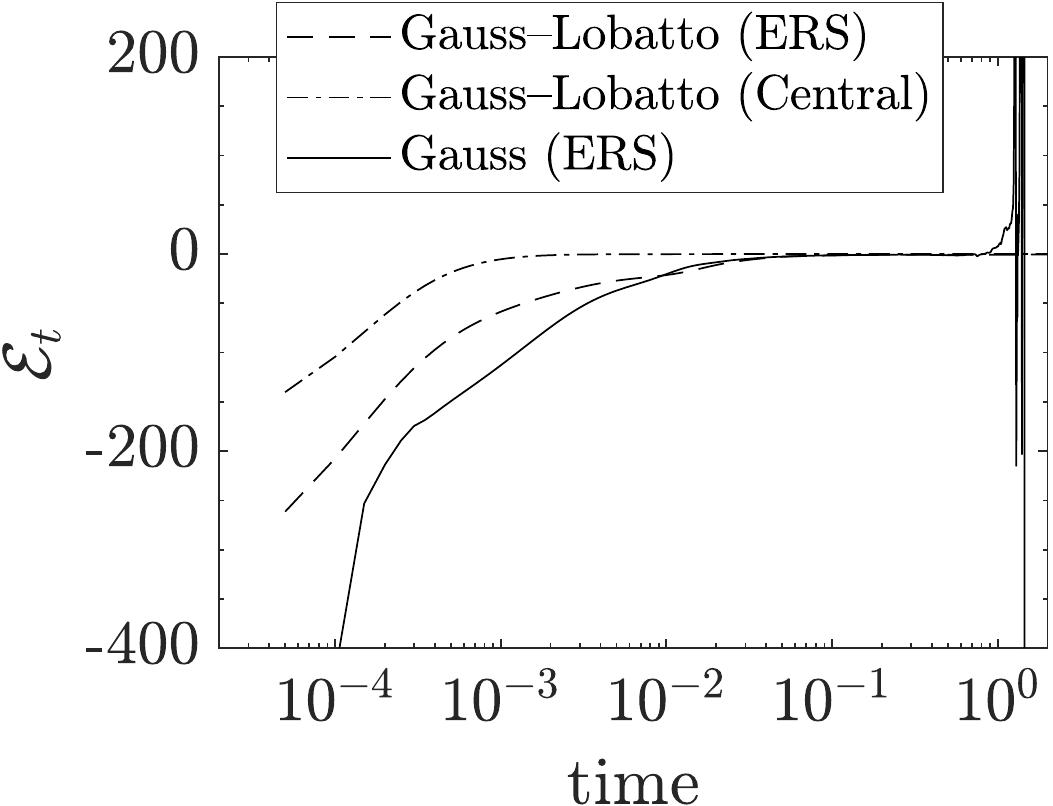}}
	\subfigure[$N=3$]{\includegraphics[width=0.32\textwidth]{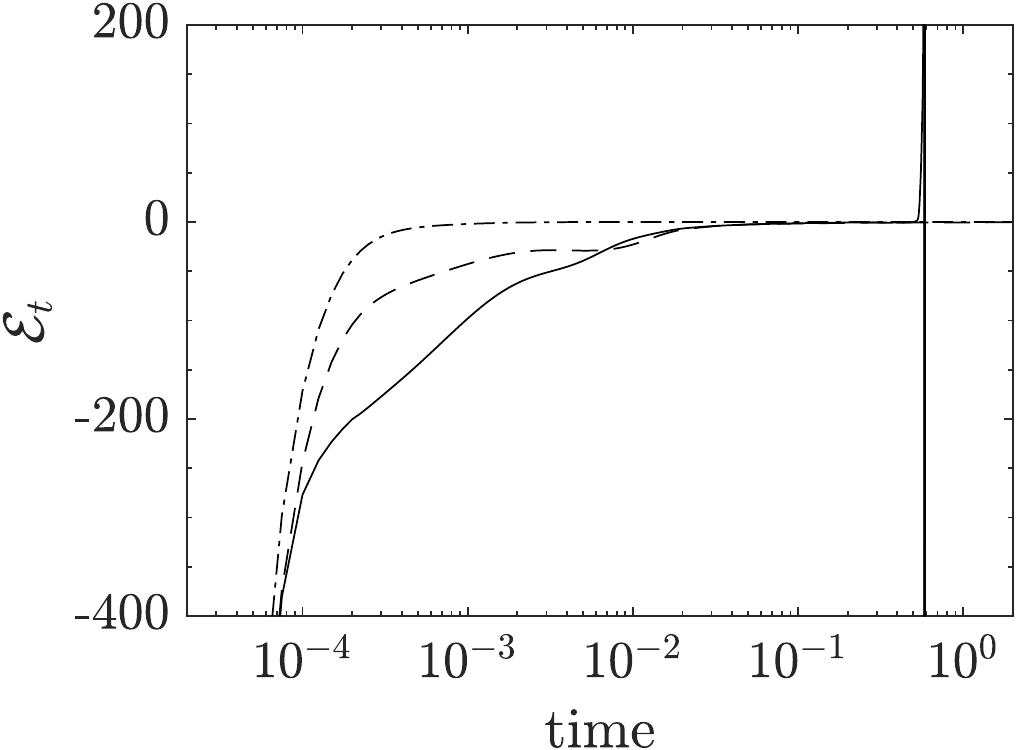}}
	\subfigure[$N=4$]{\includegraphics[width=0.32\textwidth]{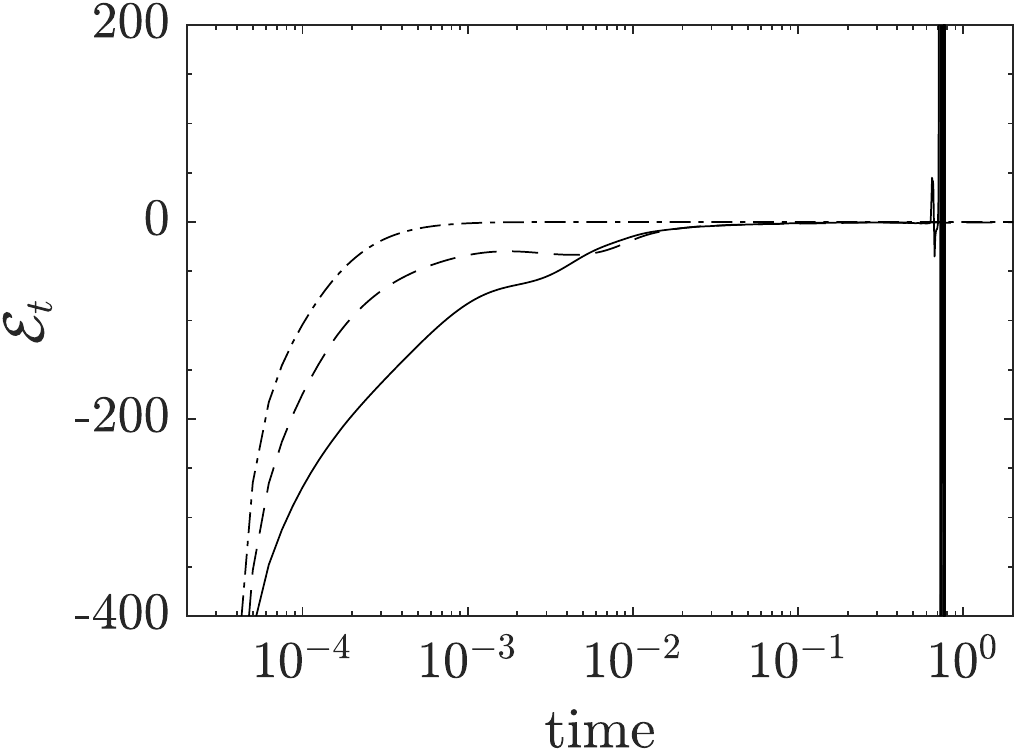}}	
	\caption{Random chosen initial condition: evolution of the entropy time derivative $\mathcal E_{t}$ for three polynomial orders $N=2,3,4$. We consider three schemes: Gauss--Lobatto points with the ERS and central fluxes, and Gauss points with the ERS. For all the polynomial orders, using Gauss--Lobatto points is an entropy preserving scheme. The dissipation is a result of physical viscosity and chemical potential for central fluxes, to which the dissipation by the exact Riemann solver can added. On the contrary, the Gauss approximation entropy time derivative is not always negative, which numerically diverges as a result}
	\label{fig:num:random:Entropy}	
\end{figure}

\begin{figure}[h]
	\centering
	\subfigure[$N=2$]{\includegraphics[width=0.32\textwidth]{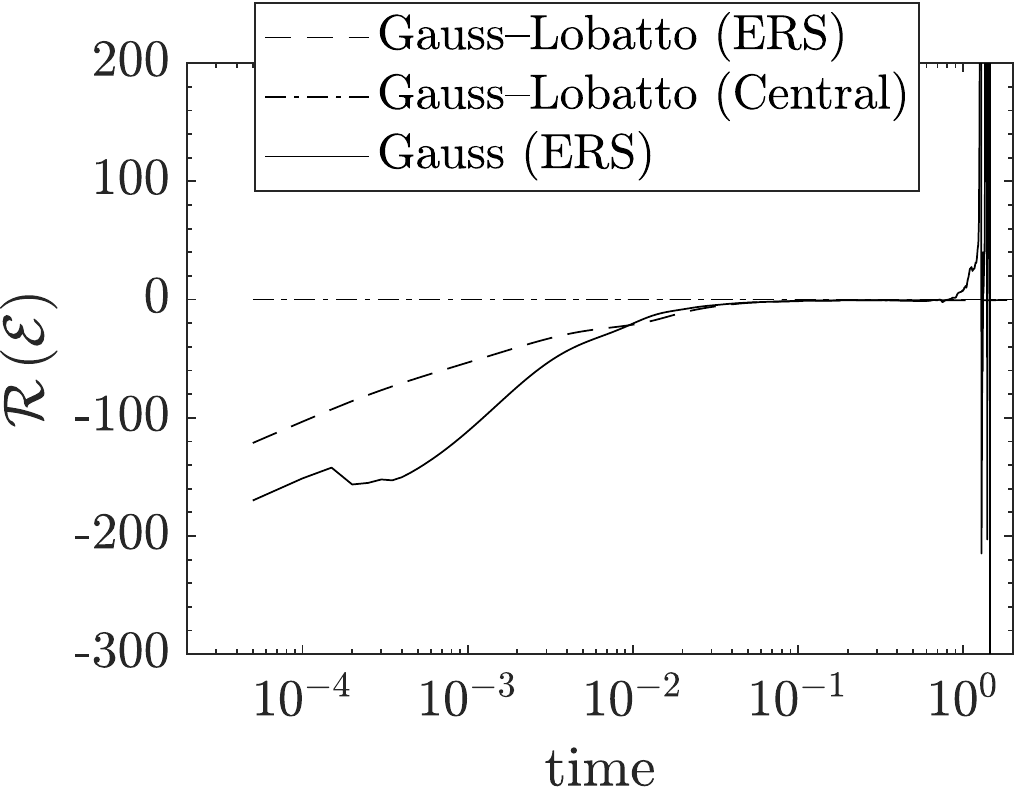}}
	\subfigure[$N=3$]{\includegraphics[width=0.32\textwidth]{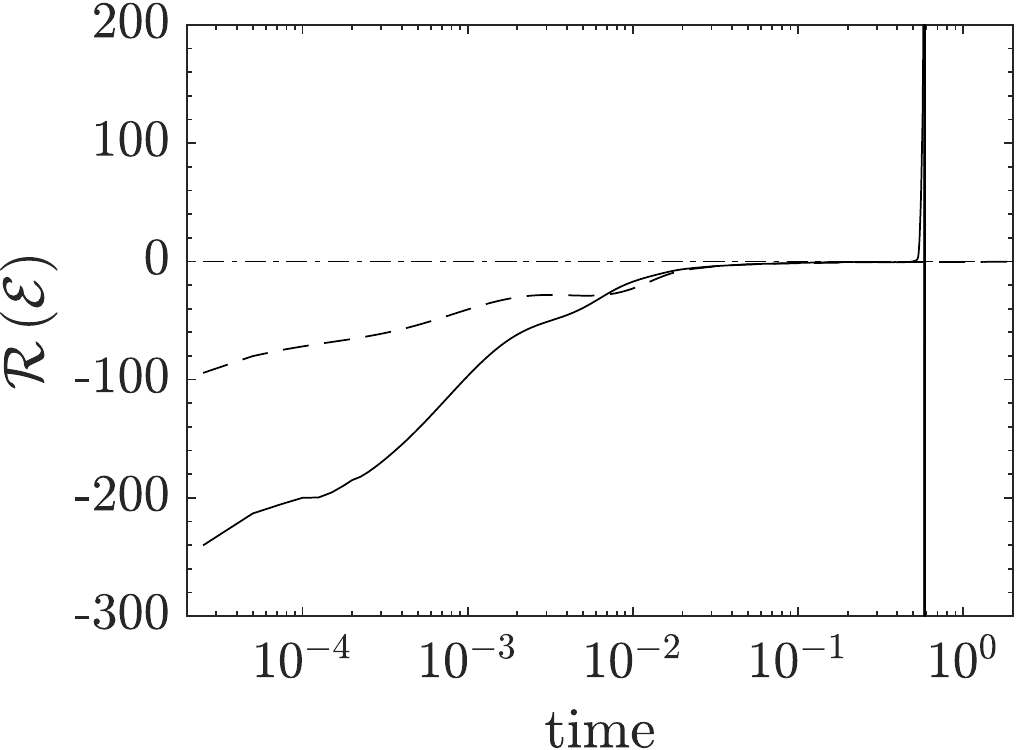}}
	\subfigure[$N=4$]{\includegraphics[width=0.32\textwidth]{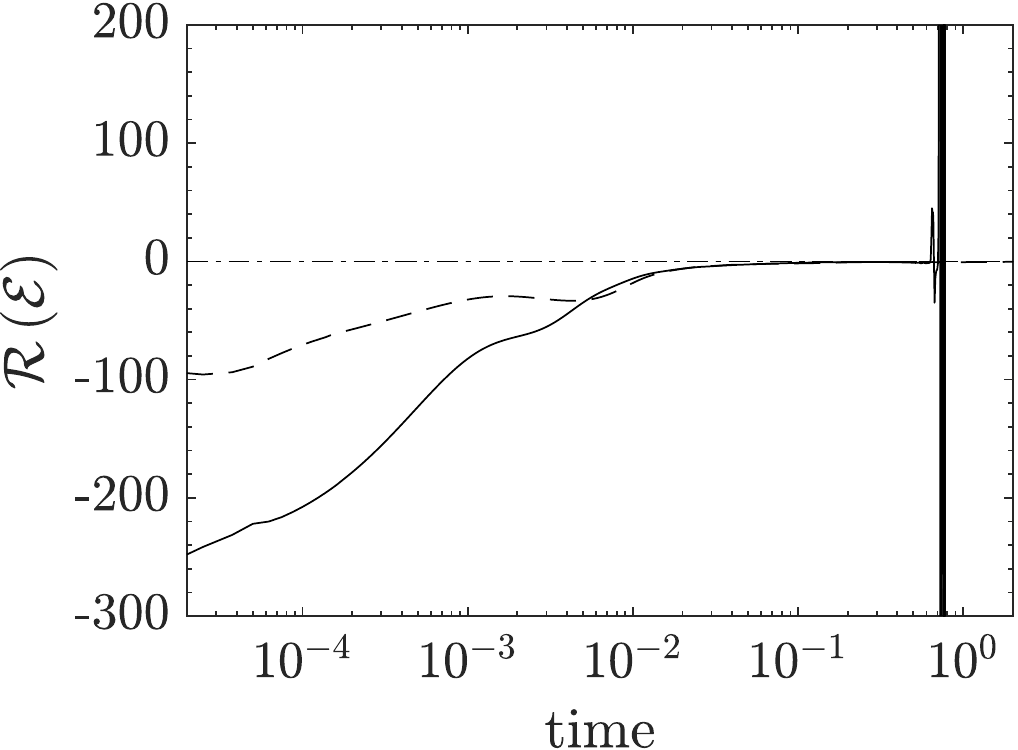}}	
	\caption{Random initial condition: evolution of the entropy remainder \eqref{eq:num:random:entropy-remainder} for three polynomial orders $N=2,3,4$. We consider three schemes: Gauss--Lobatto points with the ERS and central fluxes, and Gauss points with the ERS. The remainder studies the dissipation introduced by the scheme without considering physical viscosity and chemical potential. Thus, the Gauss--Lobatto scheme with the ERS is dissipative (always negative), the Gauss--Lobatto scheme with central fluxes is conserving (always zero), and the scheme with Gauss points is unstable and crashes}
	\label{fig:num:random:Remainder}
\end{figure}

The entropy time derivative evolution is presented in Fig.~\ref{fig:num:random:Entropy}, and the entropy remainder \eqref{eq:num:random:entropy-remainder} in Fig.~\ref{fig:num:random:Remainder}.
%DAK don't forget markers
Both are presented for three polynomial orders ($N=2,3,4$), and for the three configurations described: Gauss--Lobatto points with ERS,  Gauss--Lobatto points with central fluxes, and Gauss points with ERS. We find that the solutions follow the same pattern for the different polynomial orders. Gauss--Lobatto points are always stable, and the ERS enhances the stability (i.e. the entropy time derivative is always smaller or equal when compared to central fluxes). This is confirmed by looking at the entropy remainder in Fig.~\ref{fig:num:random:Remainder}, which is machine precision for central fluxes, and negative for the ERS. Whereas for Gauss points, the solution is divergent (i.e. crashes) when $t\geq 1.53$. We see that for this case, Gauss points are more dissipative than Gauss--Lobatto points for the first steps, but then the opposite occurs when $t> 10^{-2}$.

\begin{figure}[h]
\centering
\includegraphics[width=0.5\textwidth]{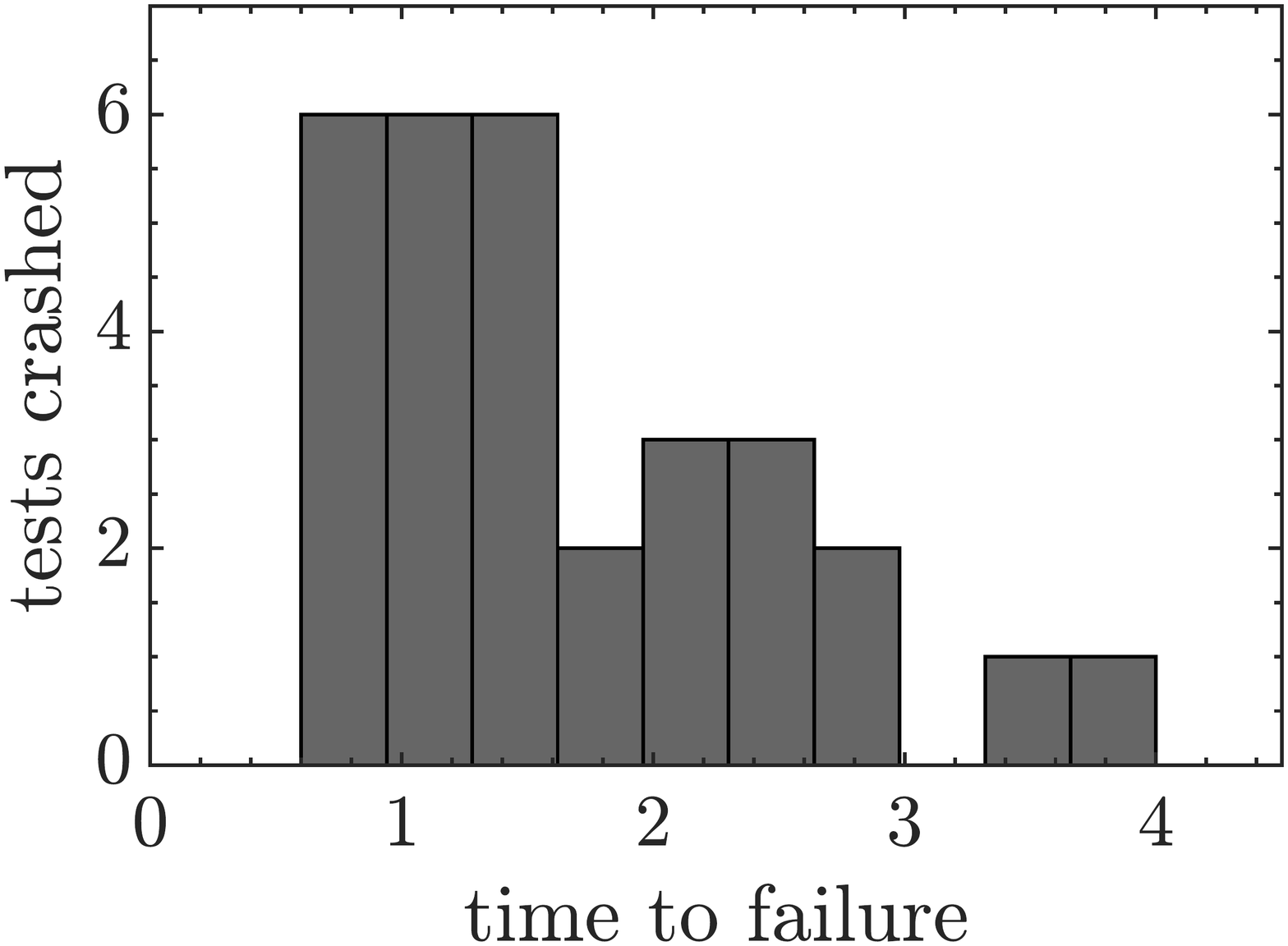}
\caption{Random initial condition: we run 100 simulations modifying the random numbers sequence, of which 30 crashed with Gauss points. In this figure, the physical time in which the simulation diverged is represented. Although the other test passed the final time $t_{F}=10.0$, their entropy time derivative was found to be positive in several time--steps. None of the 100 simulations crashed with Gauss--Lobatto}
\label{fig:num:random:Crash}
\end{figure}

With a random initial condition, the Gauss approximation is prone to diverge, but might not. We ran another $100$ random initial conditions with different seeds until a final time $t_{F}=10.0$, of which 30/100 crashed. For the 70 simulations that did not crash, most of them experienced growth of their entropy at some points in time. We represent in Fig.~\ref{fig:num:random:Crash} the number of simulations that crashed before the physical time on the horizontal axis.
Most often crashes occur early, i.e. times $t\approx 1.0$, after which the simulation is less likely to crash. None of the simulations crashed if they stayed stable for more than 4.0 seconds. In other words, if it can be computed beyond the initial stages, it is more likely to continue. 

We repeated the experiment for the entropy stable Gauss--Lobatto variant, to confirm that none of them crashed.
 There may be advantages to use Gauss points, because of higher accuracy per degree of freedom
  \cite{2011:Gassner,2016:Manzanero}, and indeed the scheme presented herein allows us to use Gauss points  in practice. However, as it is not entropy stable, it might crash under certain conditions, perhaps severely under--resolved simulations like high Reynolds turbulent flows, etc.

\subsection{Static bubble}\label{subsec:num:StaticBubble}

In this test problem we solve a steady two--dimensional bubble and validate the pressure jump that results from surface tension. In the domain $\Omega=[0,1]\times[0,1]$, the initial condition for the concentration is
\begin{equation}
c(x,y;0) = \frac{1-\tanh\left(-\frac{r-0.25}{\zeta}\right)}{2},~~ 
r=\sqrt{(x-0.5)^2+(y-0.5)^2},~~\zeta = 2.5\cdot 10^{-3},
\end{equation}
which approximates a circle with radius $R_0=0.25$. Since we project the initial condition into our polynomial ansatz, the radius obtained, $R_{e}$, subtly differs to $R_0$. The analytical pressure jump  between outside and inside of the bubble $\Delta p_{a}$ is given by the Poisson law,
\begin{equation}
\Delta p_{a} = p_{i} - p_{e} = p(0.5,0.5) - p(1.0,0.5) = \frac{\sigma}{R_{e}}.
\label{eq:num:SB:poisson}
\end{equation}
Note that we used the radius $R_{e}$, which is estimated from the final solution contour $c(x,y)=0.5$, to compute the analytical solution. The rest of the parameters, including the surface tension $\sigma=1$ N/m, are given in Table~\ref{tab:num:SB:param}.
{\begin{table}[h]
		\centering
		\caption{Parameters of the Static bubble test problem}
		\label{tab:num:SB:param}
		%\resizebox{\textwidth}{!}{%  
		\begin{tabular}{llllllll}
			\hline
			Grid & $\rho_{1,2}$ ($\text{kg}/\text{m}^3$)   & $\eta_{1,2}$ (Pa$\cdot$s) & $\varepsilon$ (m)& $t_{CH}$ (s)& $c_0^2$ (m/s$^2$)$^2$& $\sigma$ (N/m)& $\Delta t$ (s) \\ \hline
			$16^2$ &	1.0 & 1.0 & $0.16$ & 7.0 & 1.0E3 &  1.0 & $10^{-5}$\\
%			$23^2$ &	1.0 & 1.0 & $0.113$ & 7.0 & 1.0E3 &  1.0 & $7.07\cdot 10^{-6}$\\			
			$32^2$ &	1.0 & 1.0 & $0.08$ & 7.0 & 1.0E3 & 1.0 & $5\cdot 10^{-6}$\\
						$64^2$ &	1.0 & 1.0 & $0.05$ & 7.0 & 1.0E3 & 1.0 & $2.5\cdot 10^{-6}$\\			
			\hline
		\end{tabular}%}
\end{table}}
We apply periodic boundary conditions at the four physical boundaries, and use the first order IMEX scheme ($J=1$) until the residuals are kept lower than $10^{-7}$.

The results are summarized in Table~\ref{tab:num:SB:results}, which are computed with polynomial order $N=2$ and three meshes $16^2$, $32^2$, and $64^2$. For all meshes, we compute the solution radius $R_{e}$, the analytical pressure jump $\Delta p_{a}$, which are compared to the numerical solution interior pressure, $p_{i}$, exterior pressure $p_{e}$, and pressure jump $\Delta p$. Here we refer to the static pressure, and not the auxiliary pressure \eqref{eq:governing:pressure}. The results show that the pressure jump converges to that given by the Poisson law \eqref{eq:num:SB:poisson} as we refine the grid. 

In the last column of Table~\ref{tab:num:SB:results} we represent the norm of the velocity, which is of the size of the residuals. 
This implies that no parasitic currents are produced using this formulation, which are common in alternative 
Navier--Stokes/Cahn--Hilliard formulations \cite{2009:Lee}. Usually, the size of the parasitic currents is $10^{-3}$ for this problem, and since they do not disappear as the grid is refined, they can dominate in low speed simulations.

{\begin{table}[h]
		\centering
		\caption{Results of the Static bubble test problem}
		\label{tab:num:SB:results}
		%\resizebox{\textwidth}{!}{%  
		\begin{tabular}{llllllll}
			\hline
			Grid & $R_{e}$ (m)   & $\Delta p_{a}$ & $p_{i}$ & $p_{e}$ & $\Delta p$ (Pa)& $|\Delta p-\Delta p_{a}|/\Delta p_{a}$& $\Vert \svec{v}\Vert $\\ \hline
			$16^2$ &	0.220  & 4.546 & 0.033  
			& -4.571 & 4.604 &  $1.27\cdot10^{-2}$ & $4.0\cdot 10^{-10}$\\
%			$23^2$ & 0.232 & 4.306 & 0.04 & -4.179 & 4.219 & $2.02\cdot 10^{-2}$ & $2.5\cdot 10^{-12}$ \\
     		$32^2$ &	0.240 & 4.175 & 0.029 & -4.120 &  4.150 & $6.09\cdot 10^{-3}$ & $2.4\cdot 10^{-8}$\\
			$64^2$ & 0.244 & 4.100 & 0.173 & -4.094 & 4.111 & $2.75\cdot 10^{-3}$ & $3.7\cdot 10^{-9}$ \\			
			\hline
		\end{tabular}%}
\end{table}}

\subsection{Rising bubble}\label{subsec:num:RisingBubble}

The rising bubble problem has been widely used in the multiphase 
flow community to assess the space--time accuracy and robustness of the methods \cite{2009:Hysing,2007:Ding,2017:Hosseini}. The 
problem follows the trajectory of a bubble submerged in a heavier 
fluid as it rises. We consider the domain $\left(x,y\right)\in[0,1]\times[0,2]$, and a 
circular bubble with center in $(0.5,0.5)$ and diameter $0.5$. The bubble is approximated with the initial condition,
\begin{equation}
c(x,y;0) = 1 - \frac{1}{2} \left(\tanh\left(-\frac{2\left(r-0.25\right)}{\varepsilon}\right) + 
1\right),~~r = \sqrt{(x-0.5)^2+(y-0.5)^2},
\end{equation}
and the rest of the variables are initialized to zero. We consider a Cartesian mesh with element size $h=2^{-4}$, and a polynomial order $N=10$. The physical parameters 
are taken from \cite{2009:Hysing} and are summarized in Table~\ref{tab:num:RB:param}.
{\begin{table}[h]
		\centering
		\caption{Physical parameters of the Rising bubble test problem}
		\label{tab:num:RB:param}
		%\resizebox{\textwidth}{!}{%  
		\begin{tabular}{llllllllll}
			\hline
			Test & $\rho_1$& $\rho_2$ ($\text{kg}/\text{m}^3$)   & $\eta_1$ & $\eta_2$ (Pa$\cdot$s) & $\varepsilon$ (m)& $t_{CH}$ (s)& $c_0^2$ (m/s$^2$)$^2$& $\sigma$ (N/m)& $g$ (m/s$^2$) \\ \hline
			1 &	1000.0 & 100.0 & 10.0 & 1.0 & $0.03$ & 1.0E3 & 1.0E3 & 24.5 & 0.98\\
			2 &	1000.0 & 1.0 & 10.0 & 0.1 & $0.04$ & 1.0E4 & 1.0E3 & 1.96 & 0.98\\
			\hline
		\end{tabular}%}
\end{table}}
%
%{\begin{table}[h]
%		\centering
%		\caption{Physical parameters of the Rising bubble test problem}
%		\label{tab:num:RB:param}
%		%\resizebox{\textwidth}{!}{%  
%		\begin{tabular}{lllllllll}
%			\hline
%			 $\rho_1$& $\rho_2$ ($\text{kg}/\text{m}^3$)   & $\eta_1$ & $\eta_2$ (Pa$\cdot$s) & $\varepsilon$ (m)& $t_{CH}$ (s)& $c_0^2$ (m/s$^2$)$^2$& $\sigma$ (N/m)& $g$ (m/s$^2$) \\ \hline
%				1000.0 & 100.0 & 10.0 & 1.0 & $0.03$ & 1.0E3 & 1.0E3 & 24.5 & 0.98\\
%			%			2 &	1000.0 & 1.0 & 10.0 & 0.1 & $0.04$ & 80.0 & 1.0E3 & 1.96 & 0.98\\
%			\hline
%		\end{tabular}%}
%\end{table}}
%
%
We consider 
two tests: the first with moderate density ratio $\rho_1/\rho_2=10$, the second 
with large density ratio $\rho_1/\rho_2=1000$. 
The boundary conditions are 
free--slip walls in $x=0$ and $x=1$, and no--slip walls in $y=0$ and $y=2$. For the Test 1, the chemical characteristic time is not big enough to use the explicit Runge--Kutta, and we use the second order IMEX BDF method, with time--step $\Delta t=7.2\cdot 10^{-6}$. For the Test 2, we use the explicit third order Runge--Kutta, with $\Delta t=7.2\cdot 10^{-7}$. Although we use fixed time--stepping, in both cases the CFL number ($\mathrm{CFL}=16(N+1)\left(|u|+a\right)\Delta t$) is maintained close to 0.2.

Although the model presented in this paper is diffuse interface, we compare it  to a sharp interface model. Moreover, in this test case we compare the artificial compressibility method to a pressure--correction method that enforces the incompressible constraint in a transient simulation. Thus, we can assess the validity of the model and its implementation.

\begin{figure}[h]
	\centering
	\subfigure[Final shape of the bubble. We represent the interface with three contours: $c=0.1$, $c=0.5$, and $c=1.0$. The dots represent the sharp interface reference \cite{2009:Hysing}]{\includegraphics[clip,width=0.45\textwidth,trim=0cm -4cm 0cm 0cm]{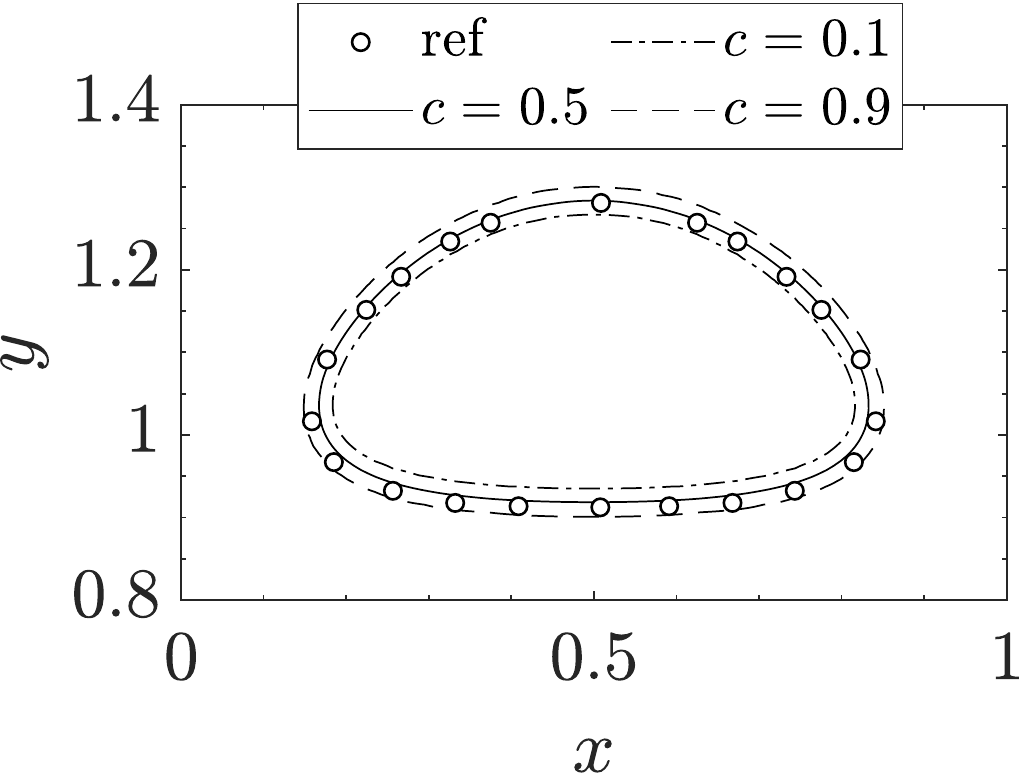}\label{fig:num:RB:FinalBubble}}\hfill
	\subfigure[Concentration contour at the final time $t=3$. The velocity field is also represented, showing the currents that raise the bubble]{\includegraphics[clip,width=0.45\textwidth,trim=20cm 0cm 20cm 0cm]{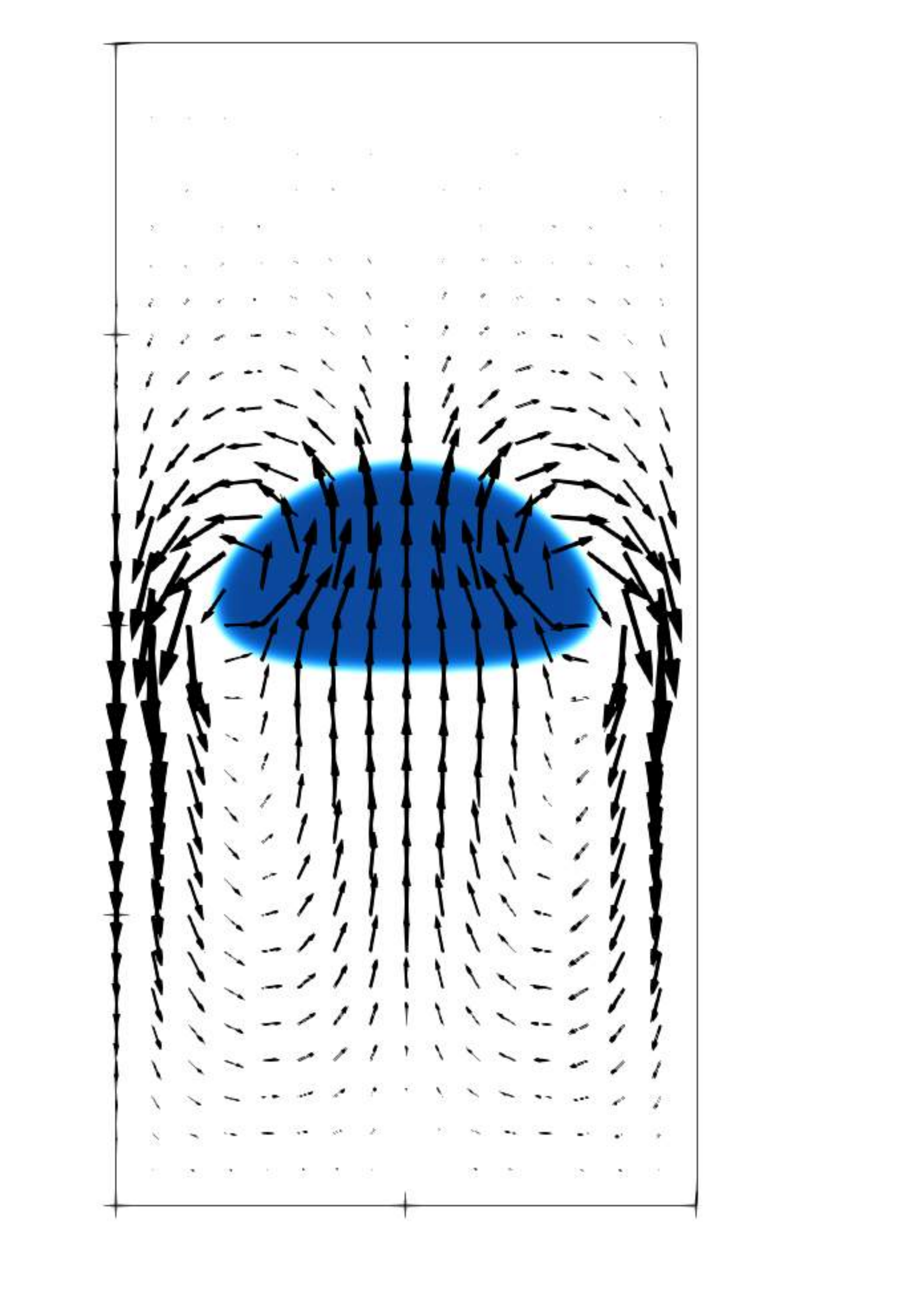}\label{fig:num:RB:FinalContour}}
	\caption{Rising bubble Test 1: position of the bubble and flow configuration at the final time $t=3.0$}
	\label{fig:num:RB:FinalShape}
\end{figure}

\begin{figure}[h]
	\centering
	\subfigure[Center of gravity]{  \includegraphics[width = 
		0.4\textwidth]{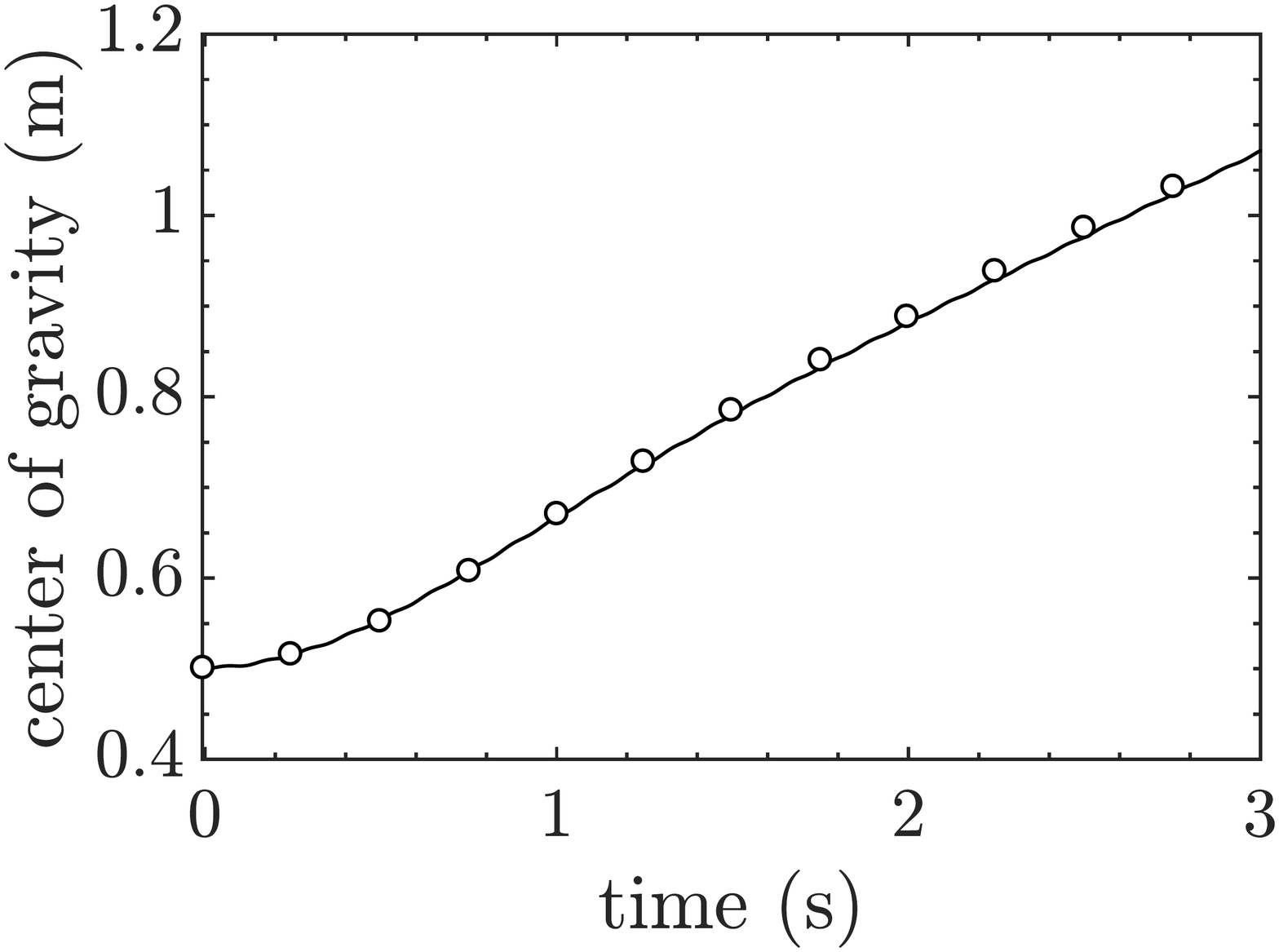}\label{fig:num:RB:CoG}}
	\subfigure[Rise velocity]{  \includegraphics[width = 
		0.4\textwidth]{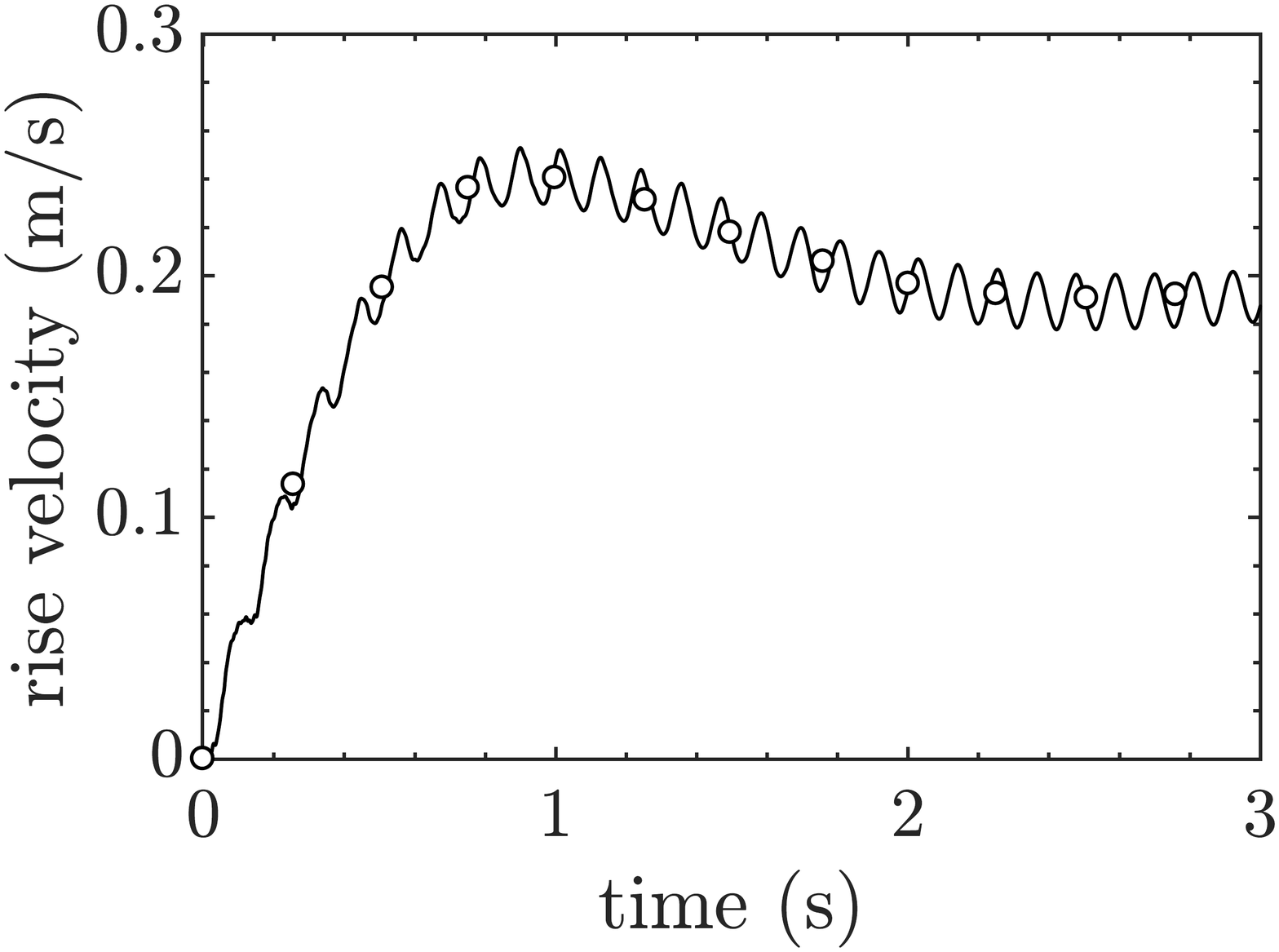}\label{fig:num:RB:RiseVel}}
	\caption{Rising bubble Test 1: evolution of the bubble center of gravity and rise velocity}
	\label{fig:num:RB:CoGandRiseVel}
\end{figure}

The shape of the bubble at the final time $t=3.0$ is represented in Fig.~\ref{fig:num:RB:FinalShape}. First, in Fig.~\ref{fig:num:RB:FinalBubble}, we represent the final shape of the interface by drawing three contour lines: $c=0.1$ in dash--dot, $c=0.5$ as a solid line, and $c=0.9$ as a dashed line. The sharp interface reference solution is represented with black dots. We find a good agreement in both the position and the shape of the bubble, as the $c=0.5$ contour follows the shape of the bubble, and the dots are always found inside the $c=0.1$ and $c=0.9$ contours. Next, in Fig.~\ref{fig:num:RB:FinalContour} we represent the concentration contour, with the velocity vectors on top. The solution obtained is also in agreement with other diffuse interface works \cite{2017:Hosseini}.
In Fig.~\ref{fig:num:RB:CoGandRiseVel} we represent the evolution of the center of gravity (Fig.~\ref{fig:num:RB:CoG}) and the bubble rise velocity (Fig.~\ref{fig:num:RB:RiseVel}),
\begin{equation}
X_{c} = \frac{1}{A}\sum_{e}\left\langle \mathcal J (1-C),X\right\rangle_{E,N},~~V_{c}=\frac{1}{A}\sum_{e}\left\langle \mathcal J (1-C),U\right\rangle_{E,N},~~A=\sum_{e}\left\langle \mathcal J (1-C),1\right\rangle_{E,N}.
\end{equation}
Both center of gravity position and rise velocity show good agreement with the sharp interface reference. 
As for the rise velocity in Fig.~\ref{fig:num:RB:RiseVel}, the oscillations are a result of the artificial compressibility pressure waves. One can reduce the amplitude of these oscillations by increasing the artificial sound speed $c_0$, but we have only found small differences in the evolution when increasing this parameter above the value used here.

\begin{figure}[h]
	\centering
	\subfigure[$t=0$ s]{\includegraphics[clip,height=5cm,trim=8cm 0cm 6cm 0cm]{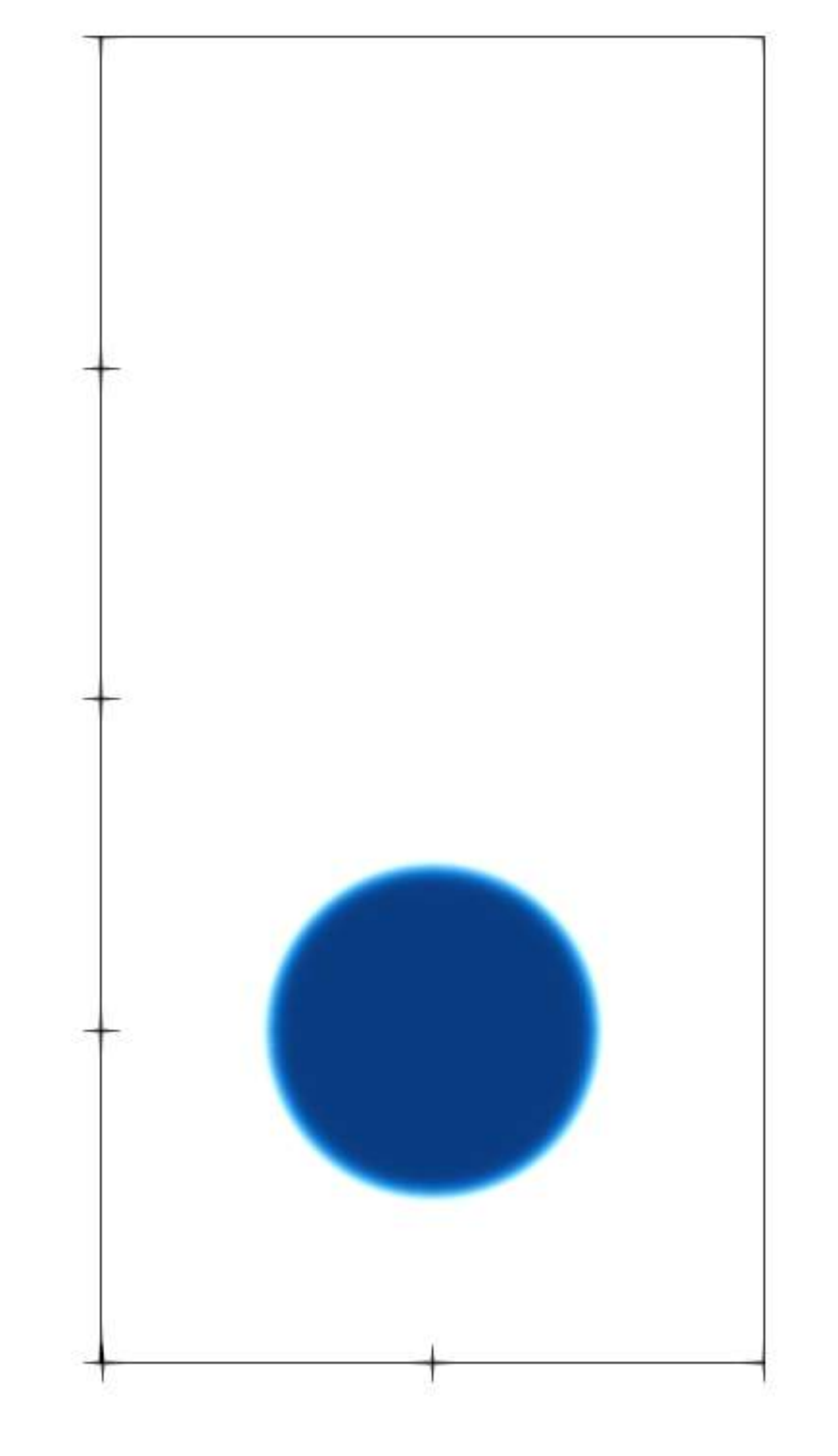}}
	\subfigure[$t=0.5$ s]{\includegraphics[clip,height=5cm,trim=7.5cm 0cm 6cm 0cm]{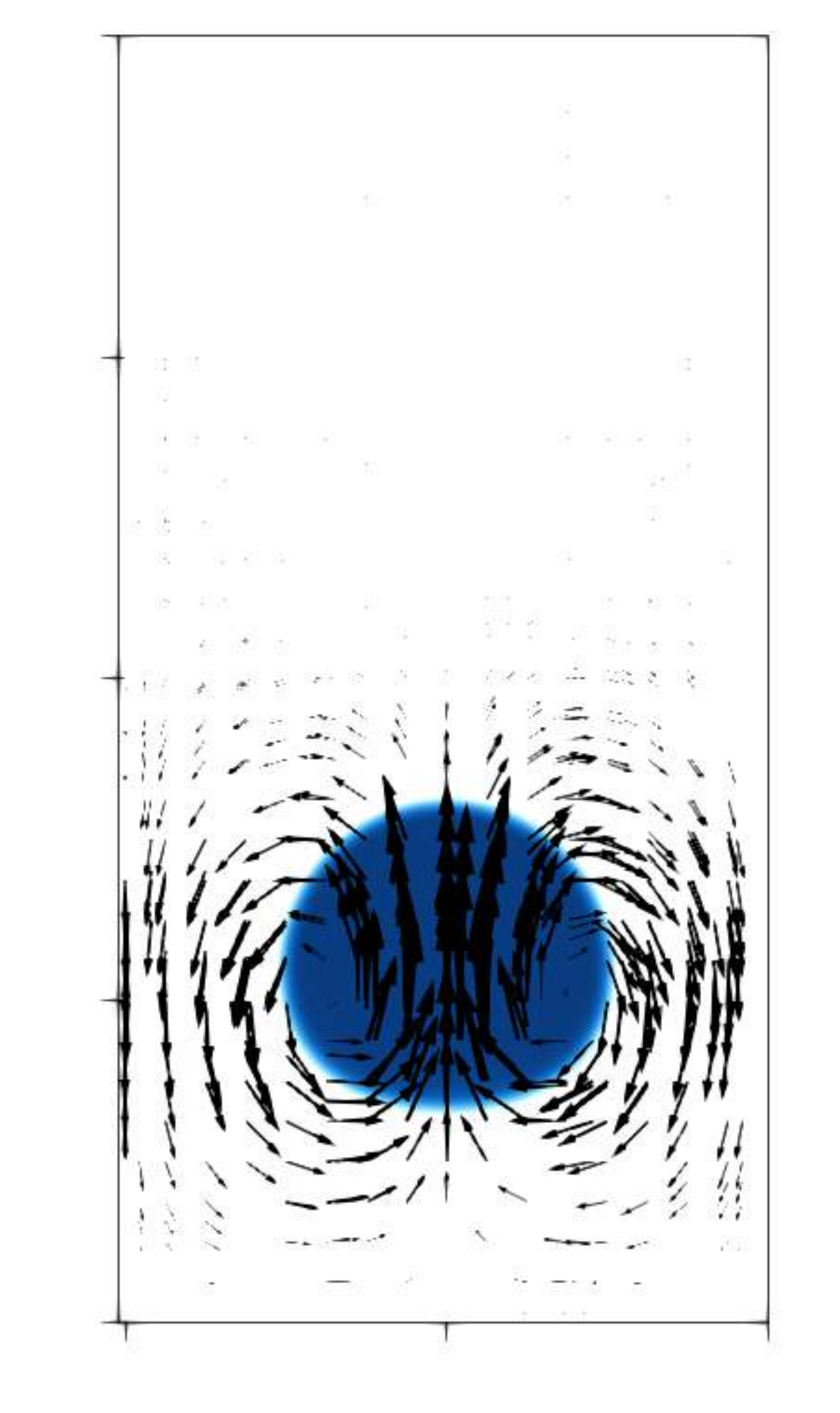}}	
	\subfigure[$t=1.0$ s]{\includegraphics[clip,height=5cm,trim=7.5cm 0cm 6cm 0cm]{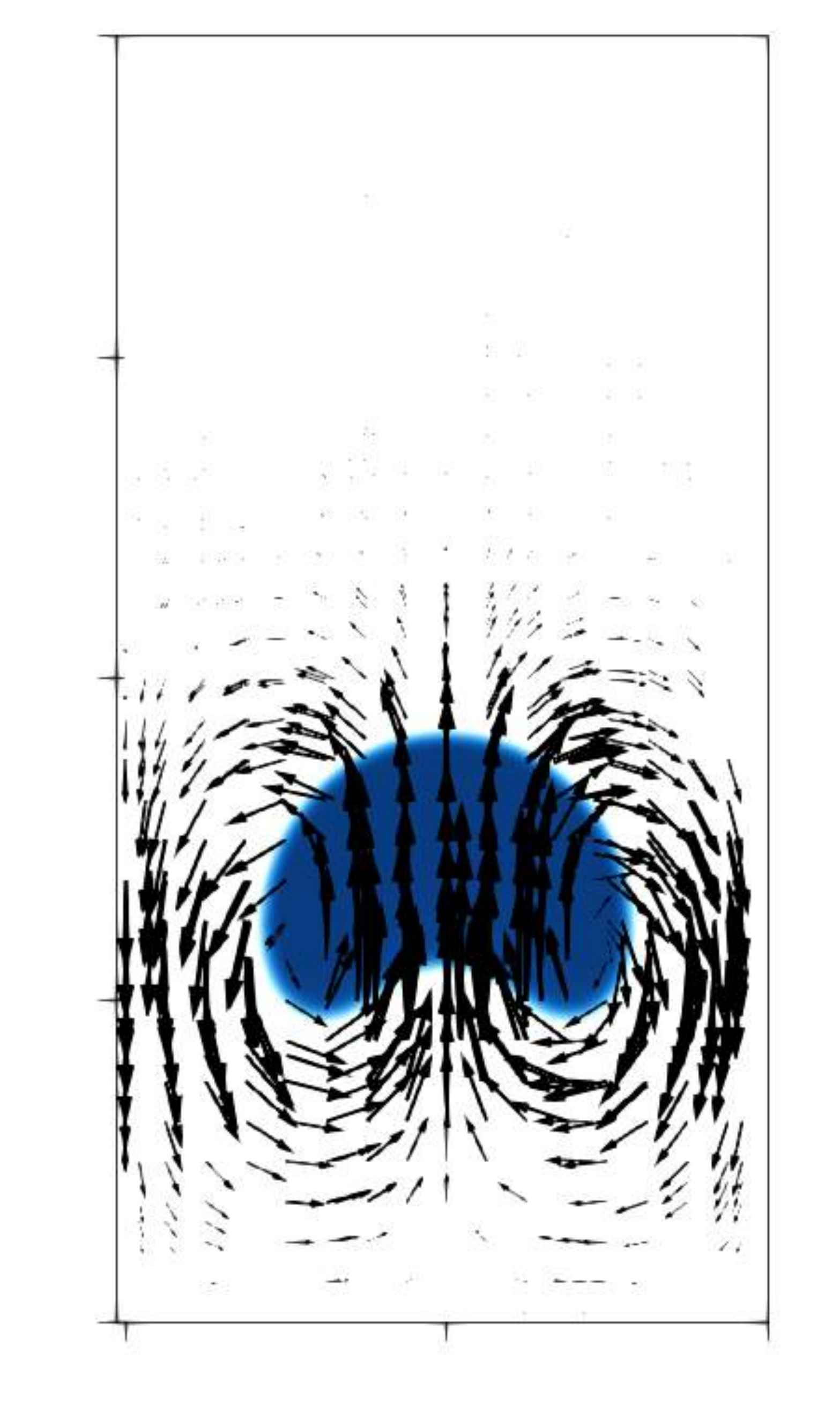}}
	\subfigure[$t=1.5$ s]{\includegraphics[clip,height=5cm,trim=7.5cm 0cm 6cm 0cm]{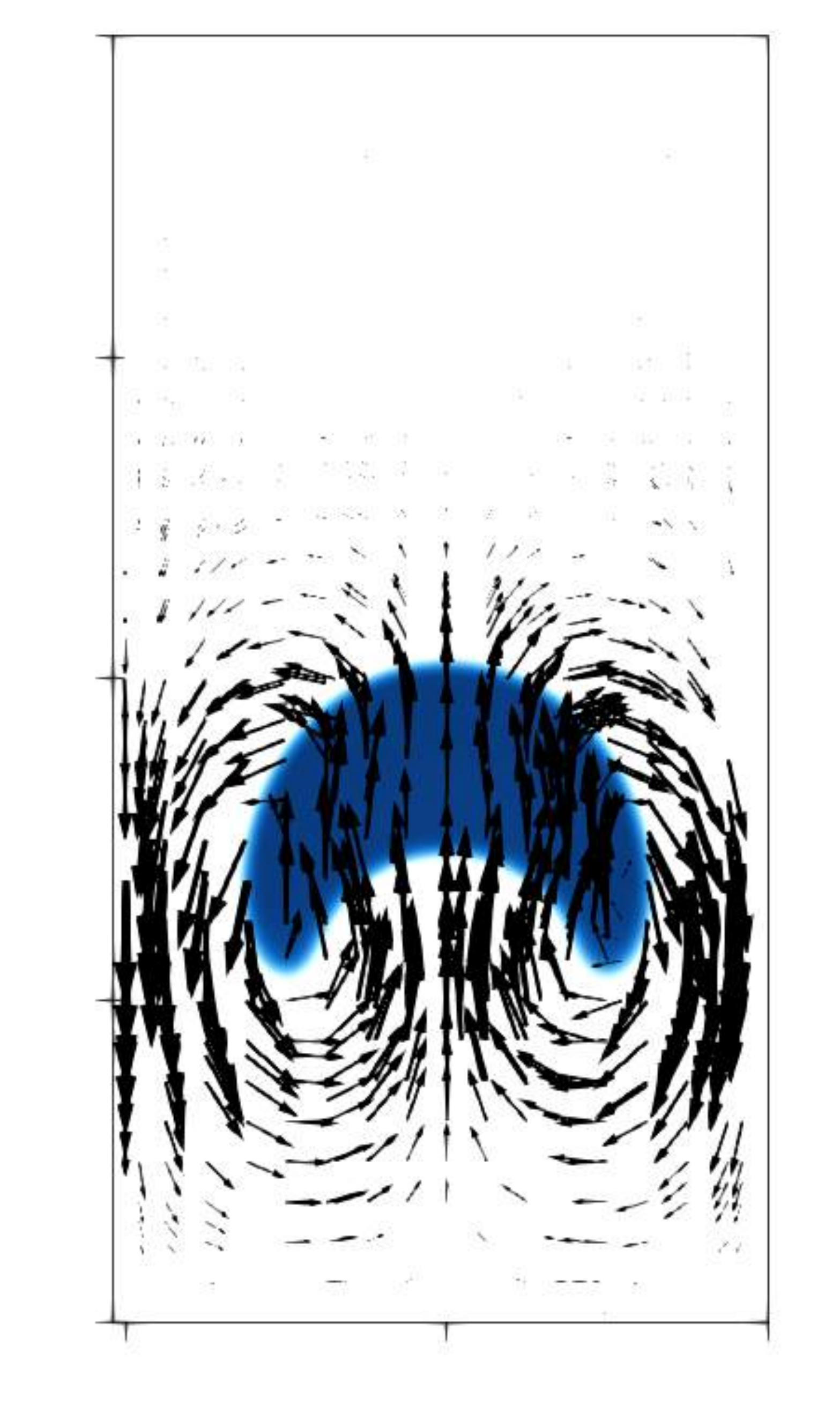}}
	\subfigure[$t=2.0$ s]{\includegraphics[clip,height=5cm,trim=7.5cm 0cm 6cm 0cm]{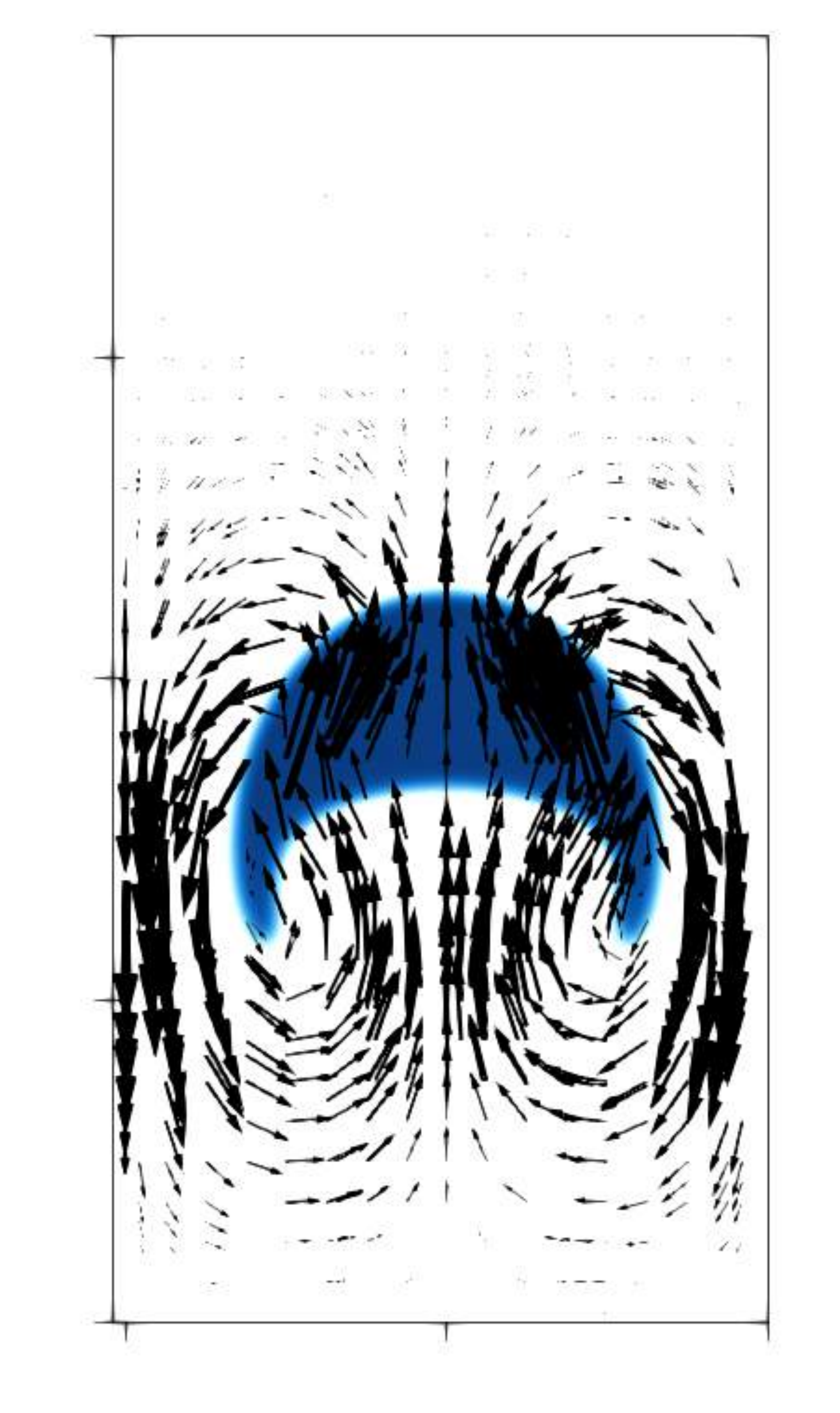}}
	\subfigure[$t=2.5$ s]{\includegraphics[clip,height=5cm,trim=7.5cm 0cm 6cm 0cm]{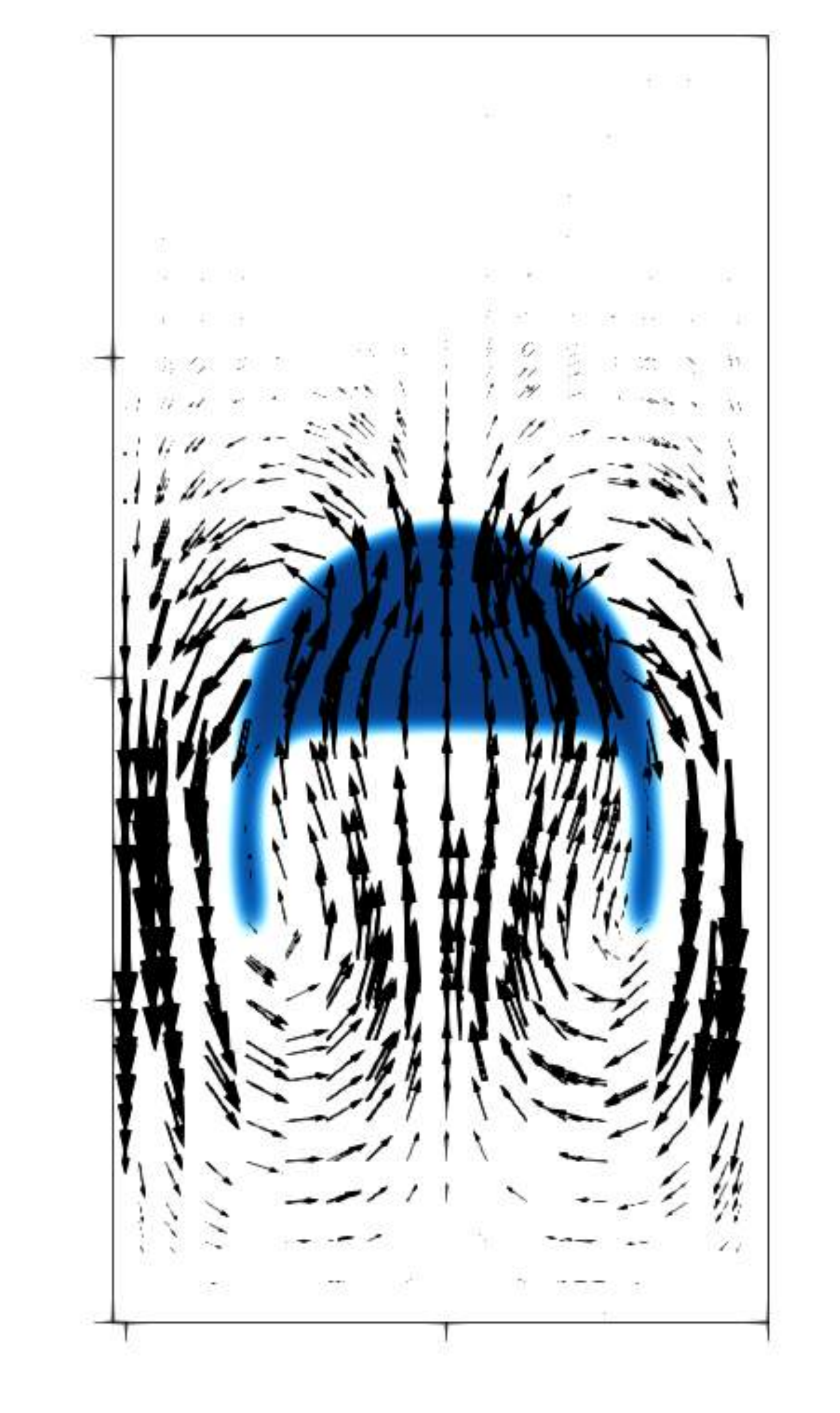}}
	\subfigure[$t=3.0$ s]{\includegraphics[clip,height=5cm,trim=7.5cm 0cm 6cm 0cm]{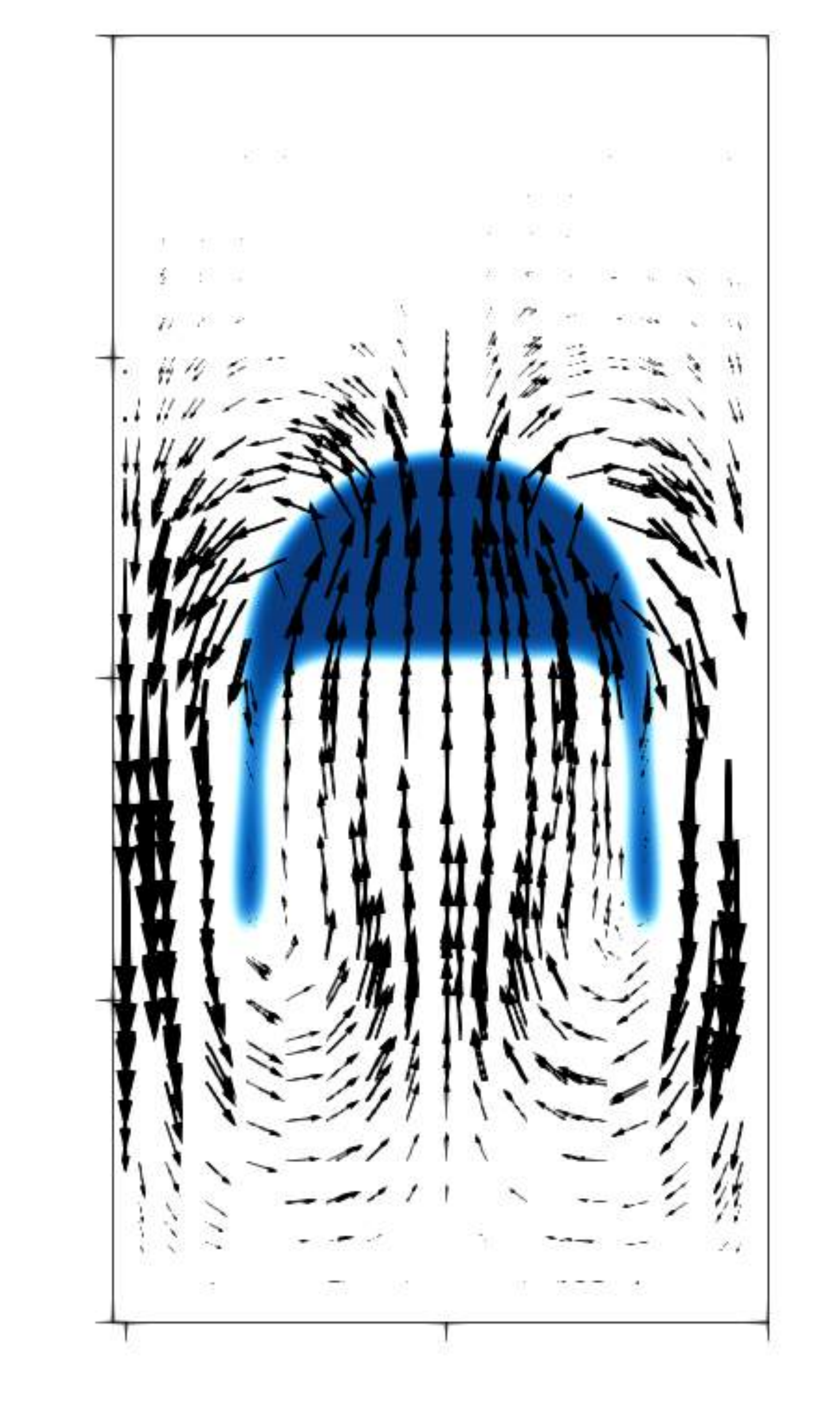}}		
	\caption{Rising bubble Test 2: different snapshots showing the position of the bubble and the flow configuration}			
	\label{fig:num:RB2:Evolution}
\end{figure}

\begin{figure}[h]
	\centering
	\subfigure[Center of gravity]{  \includegraphics[width = 
		0.4\textwidth]{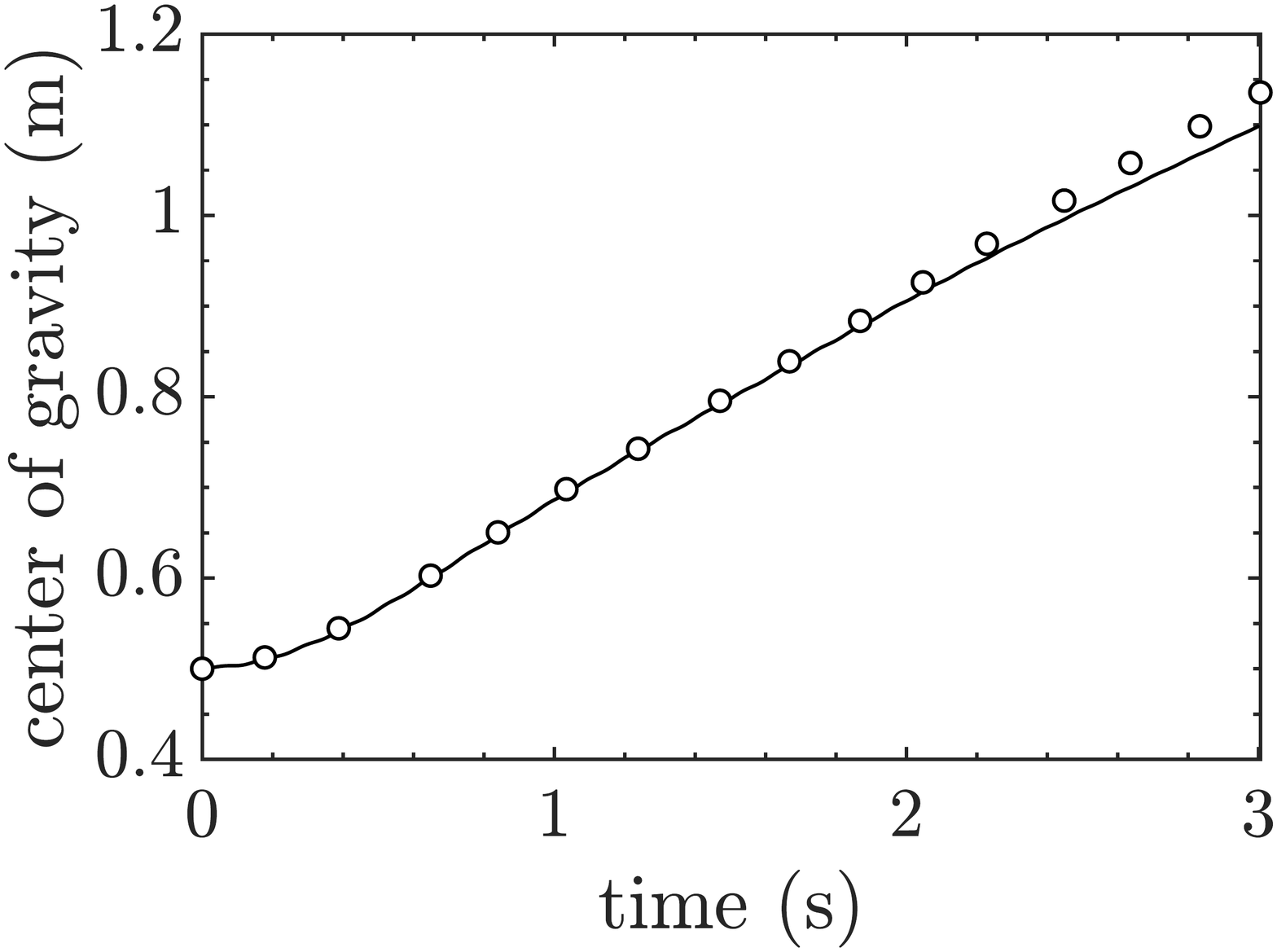}\label{fig:num:RB2:CoG}}
	\subfigure[Rise velocity]{  \includegraphics[width = 
		0.4\textwidth]{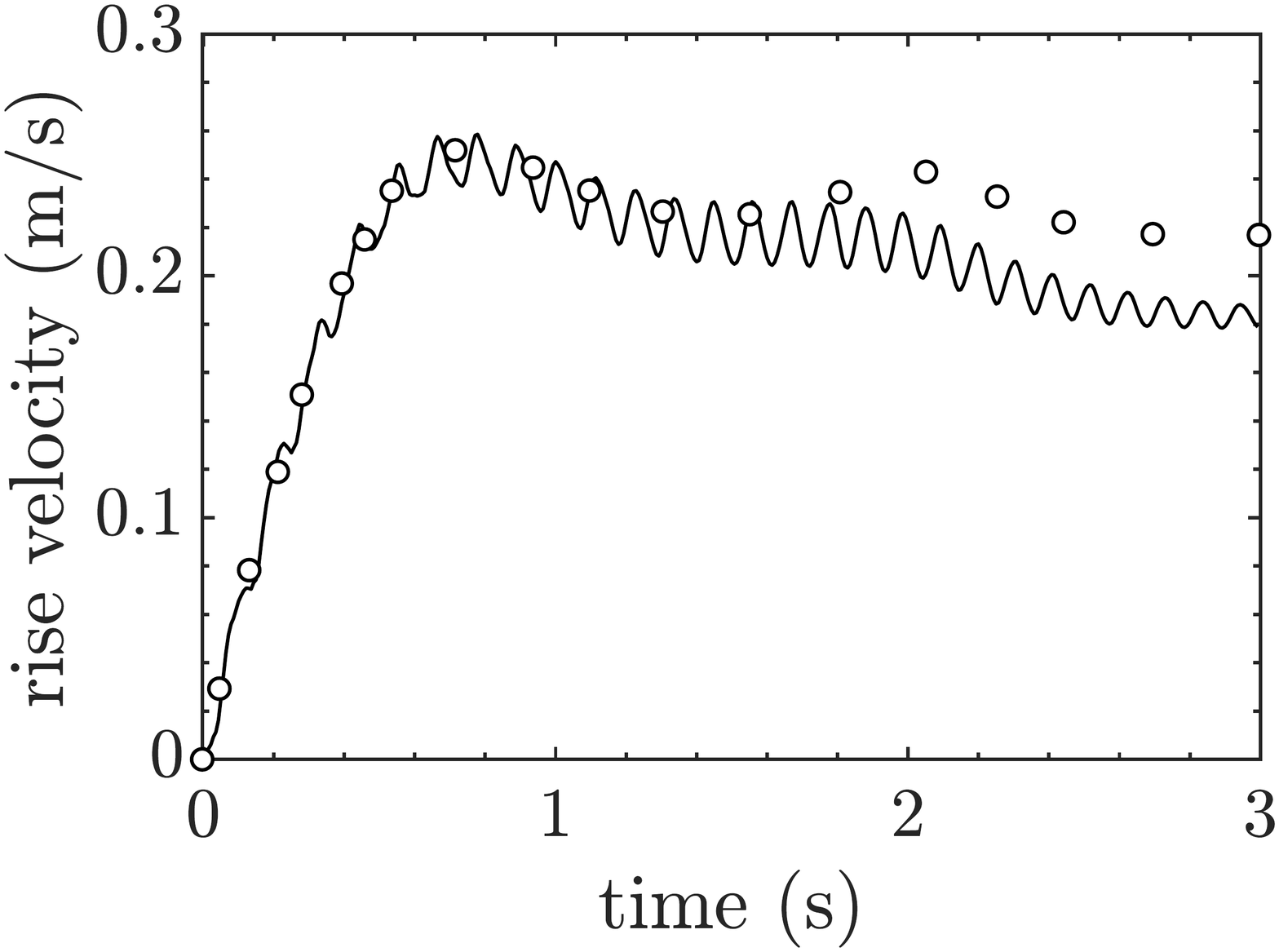}\label{fig:num:RB2:RiseVel}}
	\caption{Rising bubble Test 2: evolution of the bubble center of gravity and rise velocity}
	\label{fig:num:RB2:CoGandRiseVel}
\end{figure}

Next, we solve the more challenging rising bubble Test 2 (see Table~\ref{tab:num:RB:param}) with a higher density ($\rho_1/\rho_2=1000$) and viscosity ($\mu_1/\mu_2=100$) ratios. The evolution of the bubble and the flow configuration are represented in Fig.~\ref{fig:num:RB2:Evolution} at each $0.5$ second.
Both the shape and position of the bubble are in agreement with the sharp interface method \cite{2009:Hysing}, and other diffuse interface Cahn--Hilliard solvers \cite{2017:Hosseini}. Contrary to the rising bubble Test 1, the bubble now leaves behind an elongated skirt, which influences the velocity field. For completeness, we represent the center of gravity position and rise velocity as a function of time in Fig.~\ref{fig:num:RB2:CoGandRiseVel}, compared to the sharp interface solution provided in \cite{2009:Hysing}.
Compared to the sharp interface reference, there are visible differences. However, it has been noted in \cite{2009:Hysing,2017:Hosseini} that the solution with different sets of equations might differ for this more challenging test case. Nonetheless, the evolution of the center of gravity and the rise velocity presented here agree with other diffuse interface solvers \cite{2017:Hosseini}.

\subsection{Three--dimensional annular flow simulation}\label{subsec:num:Annular}

The last numerical experiment we present is a three--dimensional flow in a straight pipe ($L=10$) with circular section ($D=1$). Depending on the flow rates of fluid 1 and 2, the flow regime can be of different nature: stratified, slug, dispersed bubble, and annular. The last is the one we compute here, where one of the fluids behaves as a coating around the pipe surface. More details of the flow can be found in \cite{1976:Taitel,2019:Manzanero-xPipe}.

\begin{figure}[h!]
	\centering
	\subfigure[Top view]{\includegraphics[clip,width=0.9\textwidth,trim=0cm 0cm 0cm 0cm]{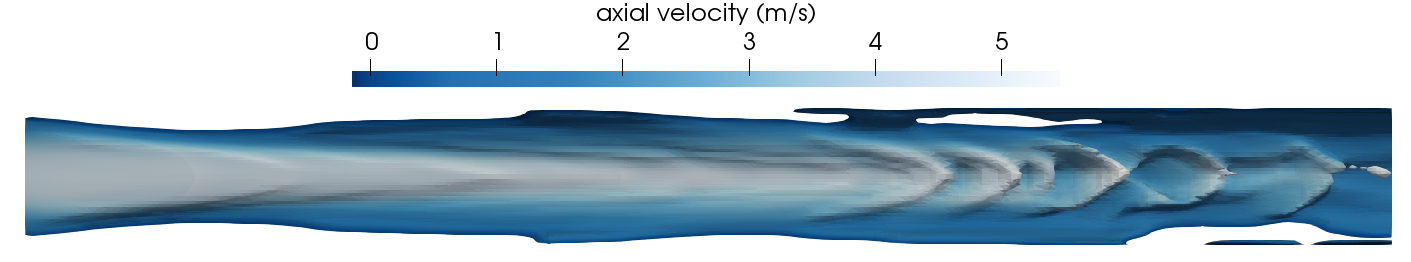}}
	\subfigure[Side view]{\includegraphics[clip,width=0.9\textwidth,trim=0cm 11cm 0cm 10cm]{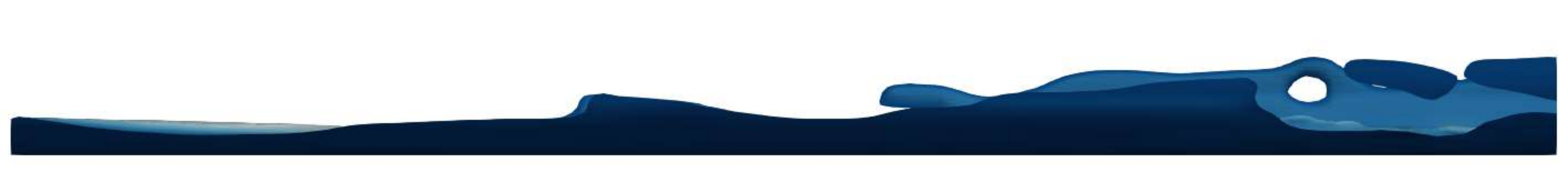}}
	\subfigure[Front view]{\includegraphics[clip,width=0.45\textwidth,trim=5cm 0cm 5cm 0cm]{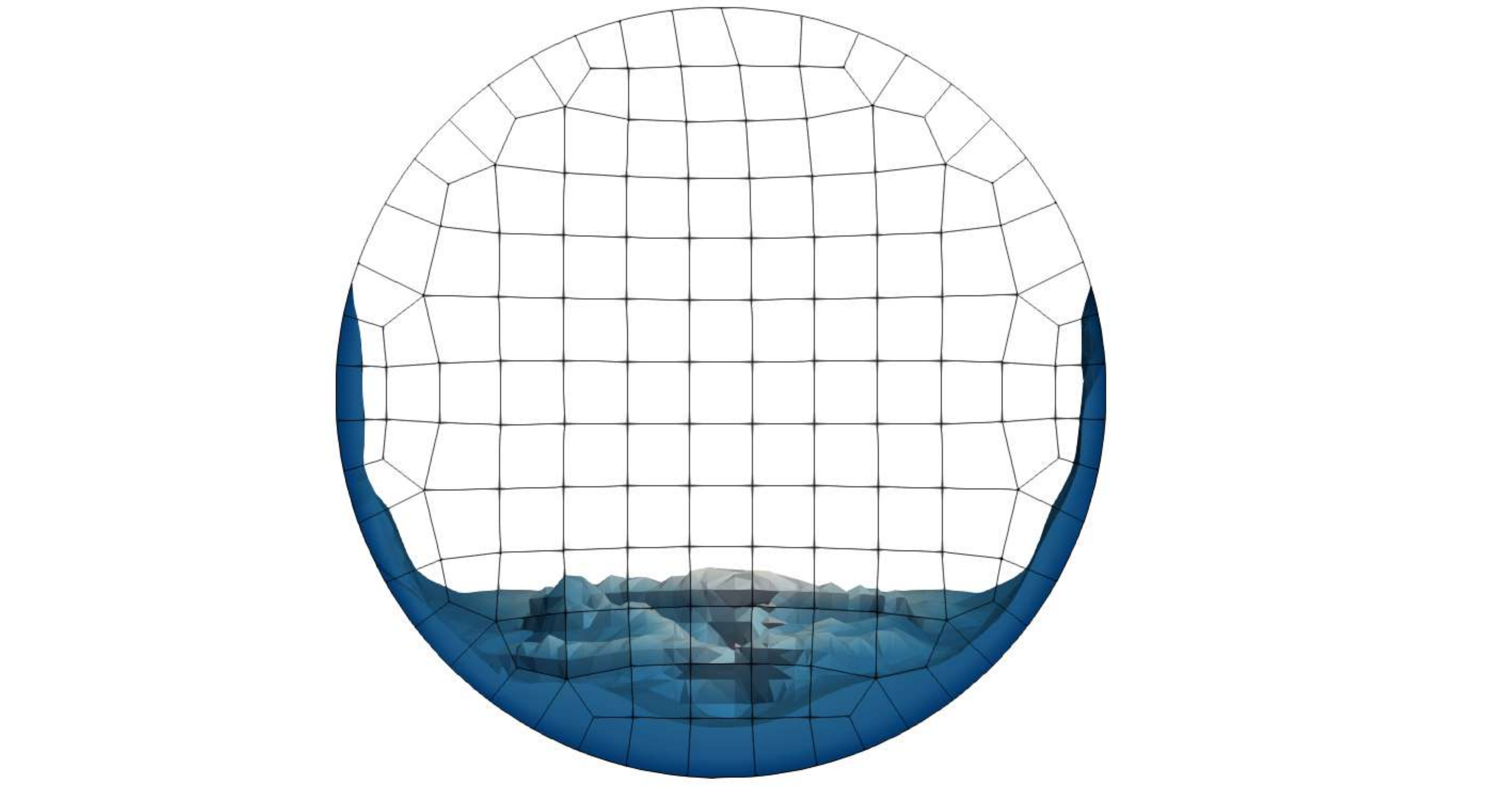}\label{fig:num:annular:front-view}}
	\subfigure[Three dimensional view]{\includegraphics[clip,width=0.45\textwidth,trim=4cm 2cm 4cm 0cm]{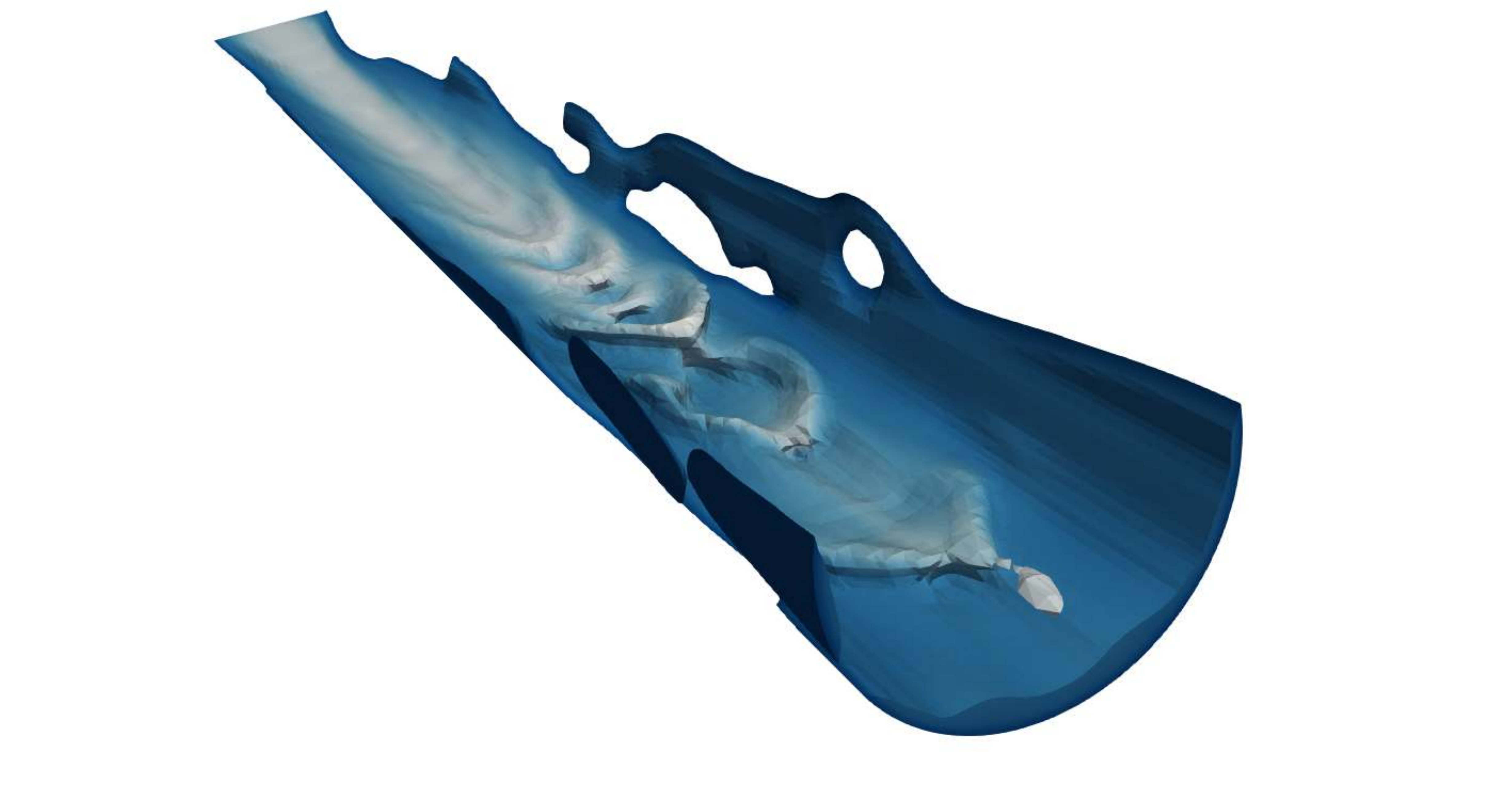}}
	\caption{Annular flow simulation. We present four view angles to the simulation at $t=4.0$ s. In the front view, the mesh cross section is also represented}
	\label{fig:num:annular:results}
\end{figure}

We construct a $z$--oriented mesh of the pipe with 8200 elements, and use a polynomial order $N=3$ (for both solutions and physical boundary representation). Note that this is the first time in this work that we use an unstructured curvilinear mesh. The cross section of the mesh is represented in Fig.~\ref{fig:num:annular:front-view}. The physical parameters are taken from the literature \cite{1976:Taitel} and given in Table~\ref{tab:num:annular:params}.
{\begin{table}[h]
		\centering
		\caption{Physical parameters of the annular flow}
		\label{tab:num:annular:params}
		%\resizebox{\textwidth}{!}{%  
		\begin{tabular}{lllllllllll}
			\hline
			$\rho_1$& $\rho_2$ ($\text{kg}/\text{m}^3$)   & $\eta_1$ & $\eta_2$ (Pa$\cdot$s) & $\varepsilon$ (m)& $t_{CH}$ (s)& $c_0^2$ (m/s$^2$)$^2$& $\sigma$ (N/m)& $g$ (m/s$^2$) & $\theta_{w}$ (deg)\\ \hline
			1.0 & 5.0 & $5\cdot 10^{-3}$ &$ 10^{-2} $& $0.0424$ & 900.0 & 1.0E3 & $2.5\cdot 10^{-4}$ & 1.0 & 45$^{o}$\\
			\hline
		\end{tabular}%}
\end{table}}

For the inflow and outflow boundary conditions, we construct auxiliary ghost states,
\begin{equation}
\stvec{Q}^{in} = \left(\begin{array}{c}
c_{in} \\
\sqrt{\rho\left(c_{in}\right)}\svec{u}_{in} \\
p
\end{array}\right),~~\stvec{Q}^{out} = \left(\begin{array}{c}
c \\
\sqrt{\rho}\svec{u} \\
0
\end{array}\right),
\end{equation}
with
\begin{equation}
\begin{split}
c_{in} &= \frac{1}{2}+\frac{1}{2}\tanh\left(\frac{x-x_0}{\varepsilon} + \frac{1}{10}\sin\left(10 z\right) + \frac{1}{10}\sin\left(20 y\right)\right), \\
u_{in} &= 0, \\
v_{in} &= 0, \\
w_{in} &= \left(V_{1}c_{in} + V_{2}\left(1-c_{in}\right)\right)\left(1-4\frac{x^2+y^2}{D^2}\right),\\
p_{in} &= 0.
\end{split}
\end{equation}
Furthermore, we use $\stvec{Q}^{in}$ as the initial condition. The velocities $V_{1}$ and $V_{2}$, and interface position $x_0$ are set so that the superficial velocities, defined as,
\begin{equation}
V_{s1} = \frac{1}{\pi R^2}\int_{\mathrm{inlet}} V_{1}c_{in}\left(1-4\frac{x^2+y^2}{D^2}\right)\diff S,~~ V_{s2} =  \frac{1}{\pi R^2}\int_{\mathrm{inlet}} V_{2}(1-c_{in})\left(1-4\frac{x^2+y^2}{D^2}\right)\diff S,
\end{equation}
are those found in \cite{1976:Taitel}. 
Precisely, we chose the velocities $V_{s1}=4.9$ m/s and $V_{s2}=0.06$ m/s, 
which produce the annular flow regime. If we place the interface between both fluids at 
$x_0=-0.262$ m (which yields a slip velocity $V_{1}-V_{2}=10$ m/s), the inlet velocities are 
$V_{1}=11.06$ m/s and $V_{2}=1.06$ m/s. 
At the walls, we also impose a contact angle $\theta_{w}=45^{o}$.
With this configuration, we can use the third order Runge--Kutta scheme as time integrator with $\Delta t=10^{-5}$. In this simulation, we take into account the effect of gravity, which points in the negative $x$--direction. 

The annular flow configuration is reached through an unstable mode, which makes the transient problem under--resolved until this mode is non--linearly damped and the final annular flow configuration is found. For that reason, we use a non--zero penalty parameter, $\beta$, on the Cahn--Hilliard equation. Precisely, we use the definition \eqref{eq:dg:penalty-parameter} with $\kappa_{\beta}=1$. Otherwise, we have found the scheme to be stable, but the solution obtained is mostly noise (i.e. it is not accurate). The introduction of the penalty parameter enhances the accuracy when the flow is under--resolved \cite{2017:Ferrer}.

\begin{figure}[h!]
	\centering
	\subfigure[$t=0.1$ s]{\includegraphics[clip,width=0.45\textwidth,trim=5cm 2cm 7cm 0cm]{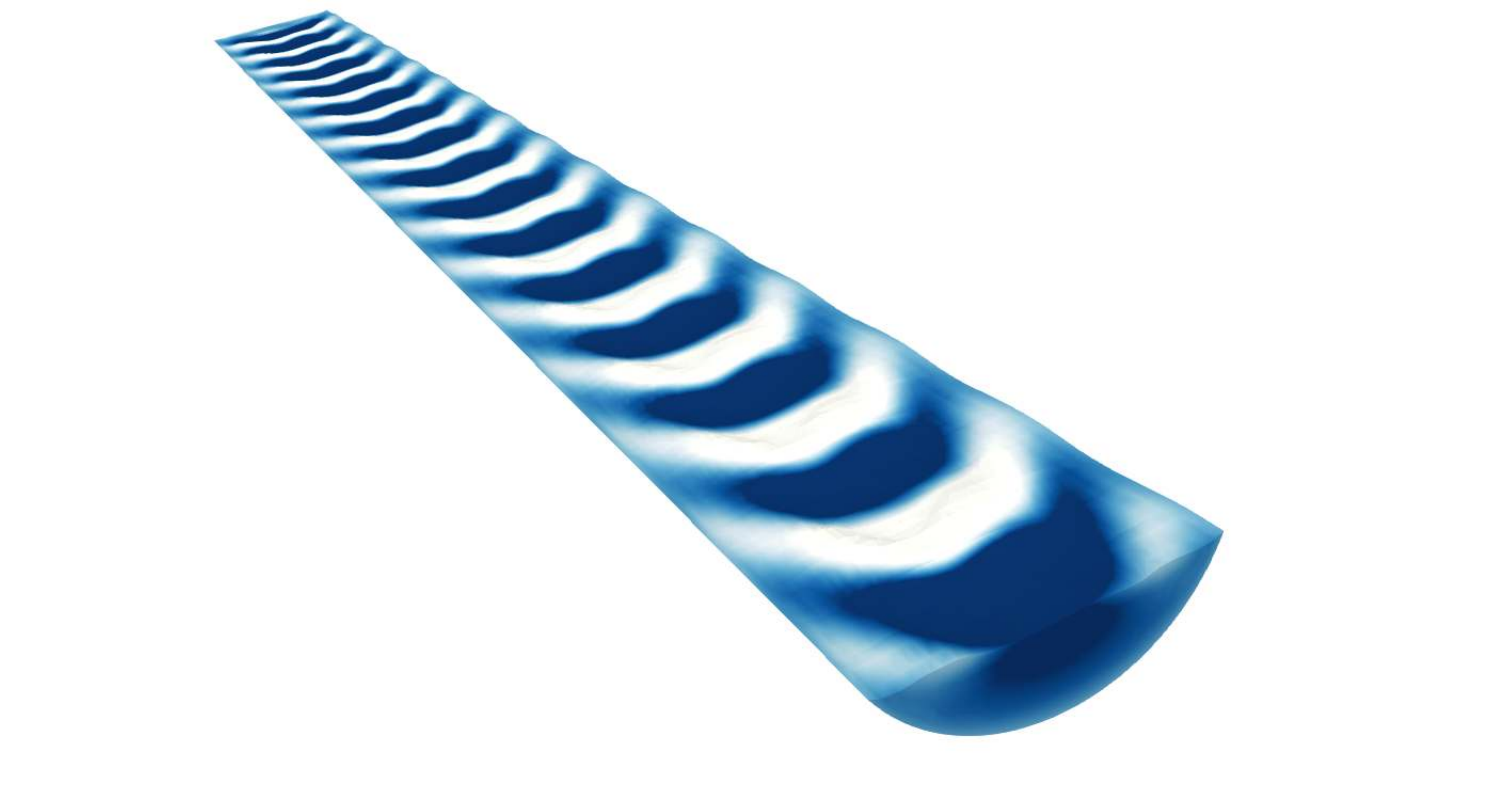}}
	\subfigure[$t=0.3$ s]{\includegraphics[clip,width=0.45\textwidth,trim=5cm 2cm 7cm 0cm]{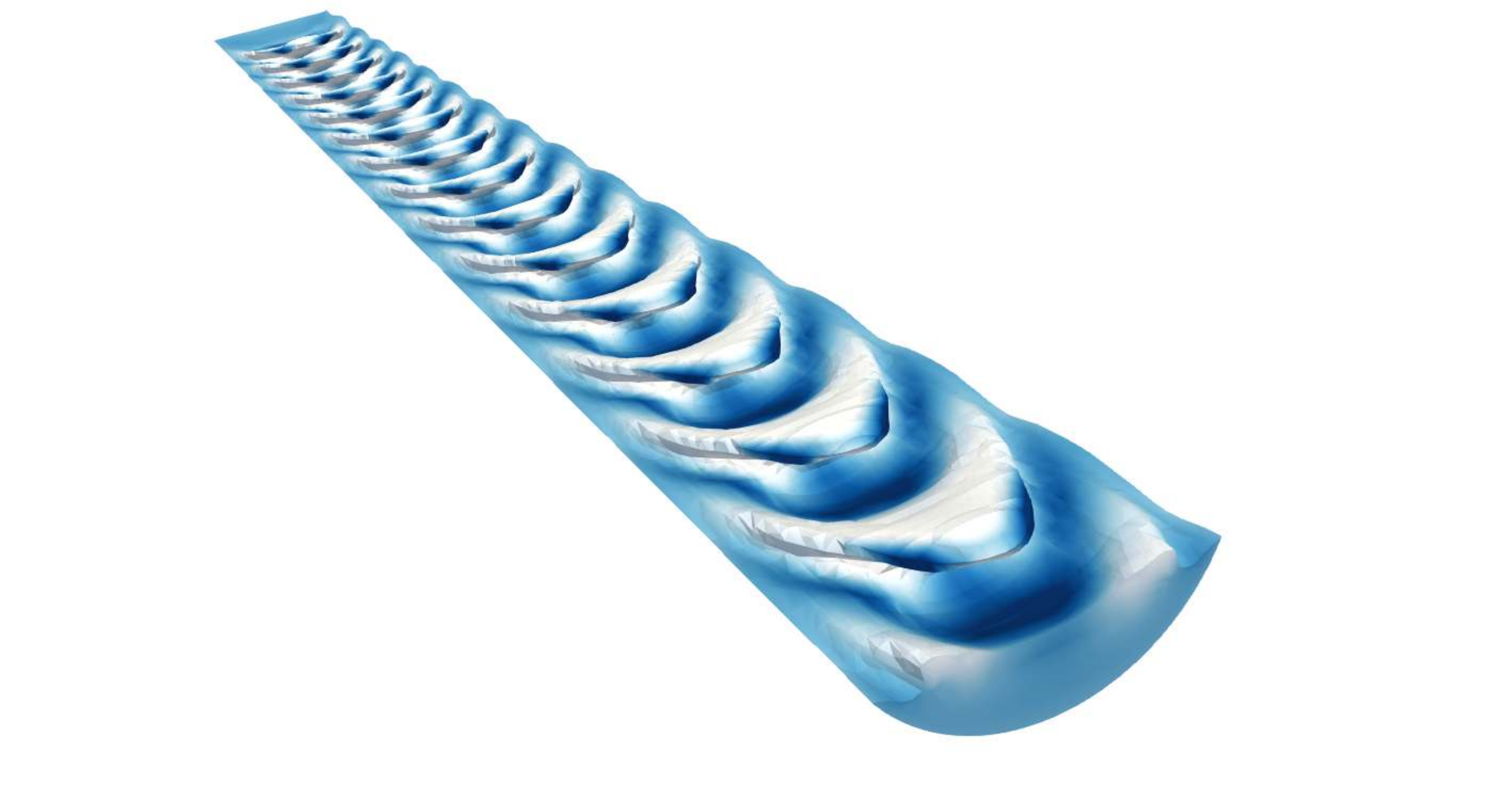}}
	\subfigure[$t=0.5$ s]{\includegraphics[clip,width=0.45\textwidth,trim=5cm 2cm 7cm 0cm]{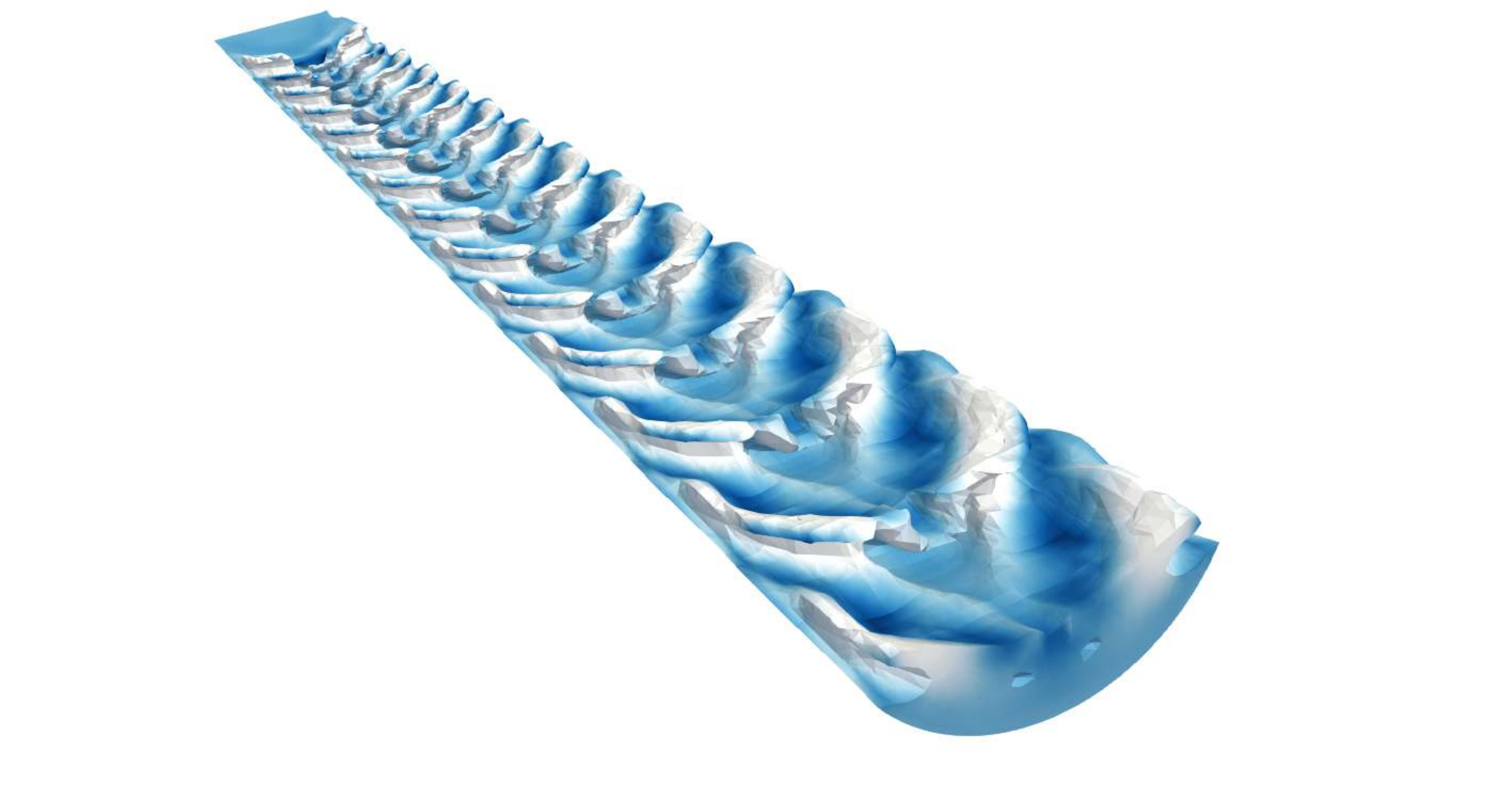}}
	\subfigure[$t=0.7$ s]{\includegraphics[clip,width=0.45\textwidth,trim=5cm 2cm 7cm 0cm]{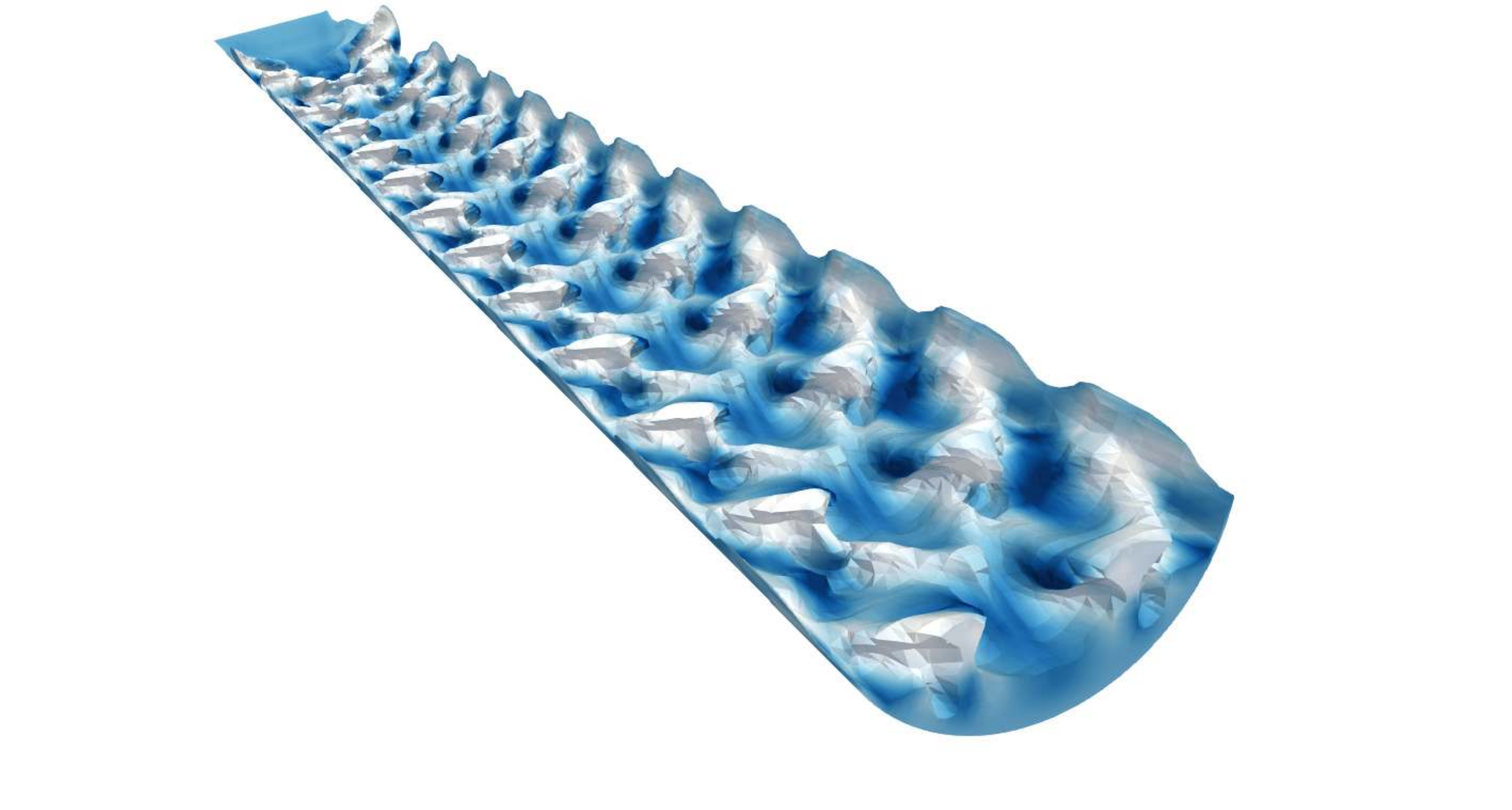}}
	\subfigure[$t=0.9$ s]{\includegraphics[clip,width=0.45\textwidth,trim=5cm 2cm 7cm 0cm]{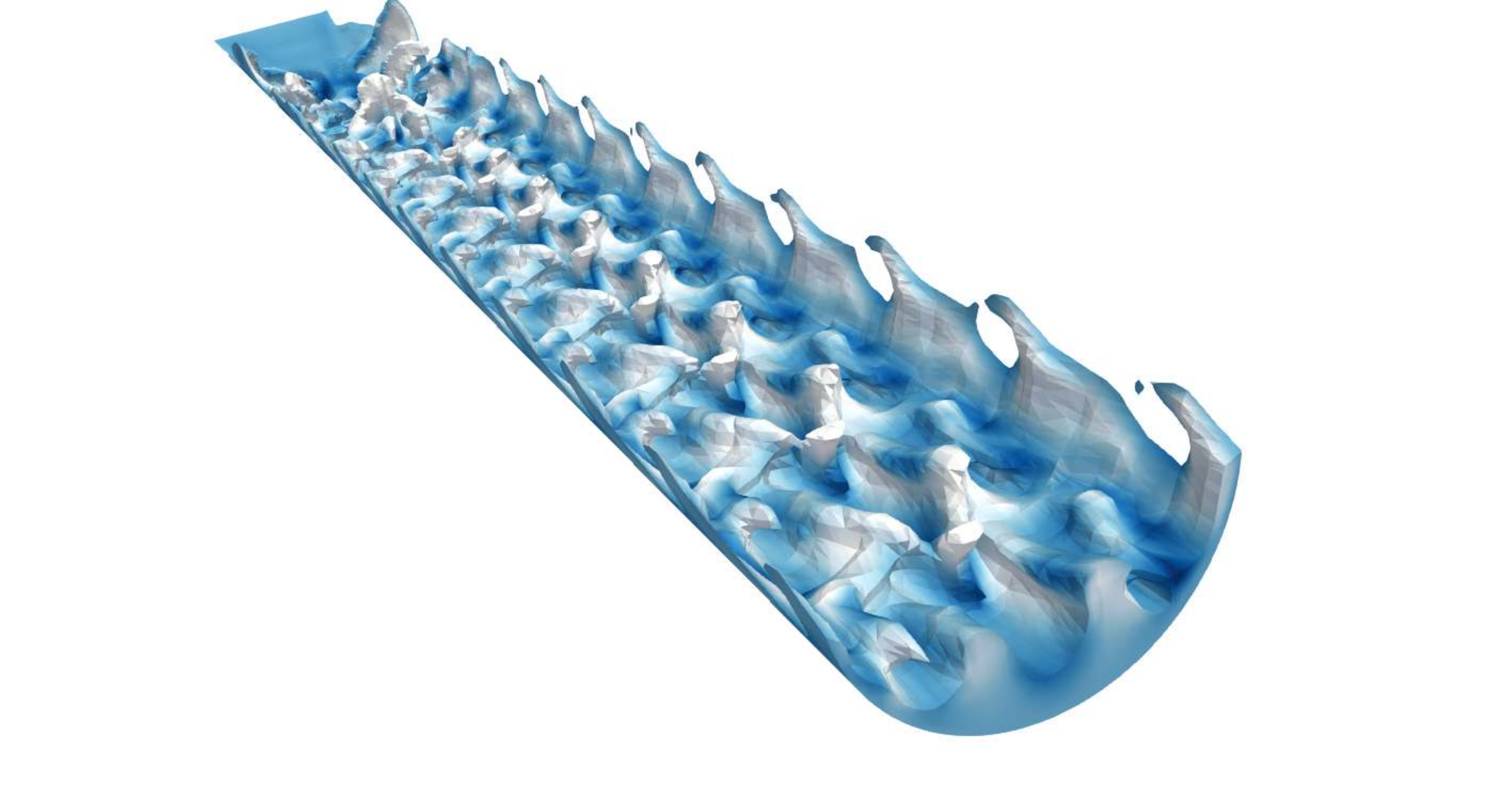}}
	\subfigure[$t=1.1$ s]{\includegraphics[clip,width=0.45\textwidth,trim=5cm 2cm 7cm 0cm]{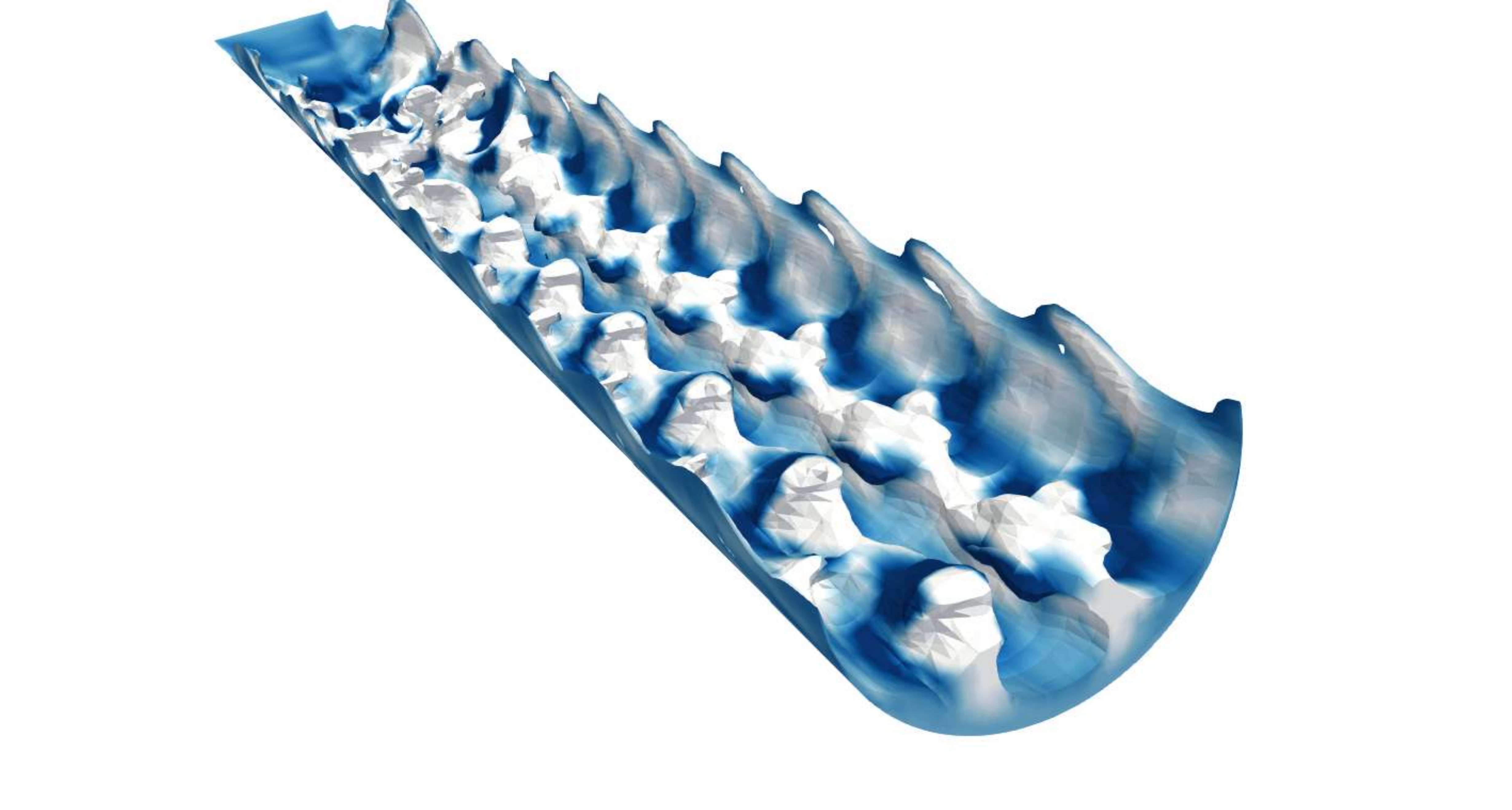}}
	\caption{Evolution of the physical mode that configures the annular flow. We have represented fluid 2 (i.e. $c\leqslant 0.5$), colored with vertical velocity, $u$. We see that during early times, there is an unstable physical mode that triggers the annular flow. When fluid 2 rises above the interface, it is dragged and accelerated by the quicker fluid 1. When it is confined below the interface, it is slowed down by the wall. This evolves non--linearly to the final state}
	\label{fig:num:annular:instability}		
\end{figure}

We represent snapshots of the initial stages of the flow in Fig.~\ref{fig:num:annular:instability}. The fluid 2 is represented ($c\leqslant 0.5$), and colored by the vertical velocity (along the $x$--axis, $u$). We see that in the initial steps there is an unstable mode that grows until it changes the flow configuration. The vertical velocity contour shows that when fluid 2 rises above the interface it is dragged and accelerated by the quicker fluid 1, and the contrary occurs when it is confined below the interface because of the wall. This produces a wavy pattern at $t=0.3$s, which eventually wraps the fluid to the wall and produces the annular flow. In the early stages, the flow is highly under--resolved.

We represent a final snapshot for $t=4.0$ s in Fig.~\ref{fig:num:annular:results}, showing the contour $c\leqslant 0.5$ (i.e. the region occupied by fluid 2), colored by the longitudinal velocity (i.e. along the $z$--axis, $w$).
Although the interface is flat at the inlet, the flow is unstable and breaks down into the 
annular flow configuration. The heavy fluid (fluid 2) is confined to the wall by the light fluid,
which is introduced 82 times faster at the inlet. In Fig.~\ref{fig:num:annular:front-view} 
we represent the front view of the pipe, with the mesh detailed on top of it. We see that the 
flow gets smoother after the numerical and physical dissipation take over the mode breakdown after the under--resolved stages in 
Fig.~\ref{fig:num:annular:results}.

\section{Summary and conclusions}\label{sec:Conclusions}

We have derived a two--phase flow model that combines the Cahn--Hilliard 
equation \eqref{eq:governing:cahn--hilliard}, a skew--symmetric version of the momentum equation
\eqref{eq:governing:momentum-skewsymmetric-sqrtRho}, and an artificial 
compressibility method \eqref{eq:governing:ACM} to get the pressure. 
Among the many available options, the versions were chosen to satisfy an entropy 
inequality,  \eqref{eq:continuous:entropy-balance-w-bcs}.  The entropy 
inequality bounds the mathematical entropy inside a finite domain $\Omega$ limited with wall 
boundary conditions by the entropy of the initial condition. As time marches, the mathematical
entropy decreases as a result of the physical and the chemical potential dissipation.

We then constructed a DG approximation that satisfies the SBP--SAT property, 
which allowed us to mimic the continuous entropy analysis discretely. As is usual in a DG approximation, we had choices for the fluxes at inter--element and physical boundaries. We 
studied two options:
\begin{enumerate}
  \item Entropy conserving, using central fluxes for the advective fluxes and the BR1 scheme for the diffusion. This choice gives a bound that
  discretely mimics the continuous entropy bound 
  \eqref{eq:continuous:entropy-balance-w-bcs}. This scheme is not of 
  practical use, since some amount of numerical dissipation is 
  required to provide accurate solutions when non--linear terms are present 
  \cite{2016:Gassner,2016:Manzanero,2017:Flad,2018:Manzanero-role}. However, it 
  serves as a baseline model to verify the stability proofs and to obtain a 
  dissipation--free scheme.
  
  \item Entropy stable, with an exact Riemann solver for the advective fluxes and the BR1 scheme for diffusive,  with additional interface dissipation for the Cahn--Hilliard equation. This 
  scheme transforms the entropy balance to an entropy inequality in the discrete version 
  of  \eqref{eq:continuous:entropy-balance-w-bcs}, as a result of the numerical 
  dissipation at the inter--element faces. This scheme uses the exact Riemann 
  solver presented in \cite{2017:Bassi} for the incompressible Navier--Stokes 
  equations, with an appropriate choice of the diamond 
  fluxes that arise from the non--conservative terms, and a modification of the discrete entropy to account for solution 
  discontinuities in the concentration as interfacial energy in the numerical solution.
\end{enumerate}
Both the scheme and the stability proofs hold for three--dimensional unstructured meshes with
hexahedral curvilinear elements. 

We selected two choices to march the scheme in time: an explicit third--order 
Runge--Kutta method, and an implicit--explicit BDF with first or second order or 
accuracy. The former is used when the mobility is small enough to not severely restrict the time--step 
size, and the latter otherwise.

We test the scheme, addressing its accuracy with a manufactured solution 
convergence analysis, and its robustness by initializing the flow with random 
initial conditions. We showed that the scheme converges spectrally as 
expected, and that the scheme is robust in the sense that none of one hundred 
simulations crashed from random initial conditions, a high density ratio ($\rho_2/\rho_1 = 
1000$), and a high Reynolds number ($\Rey = 10^{6}$). We compared the
entropy--stable scheme with the more accurate Gauss counterpart (which is not provably entropy
stable) to find that the latter crashes in $30\%$ of the simulations. We also solved commonly used static and 
rising bubble test problems to assess the steady state and transient accuracy of 
the solver. Finally we challenged the method by solving a three--dimensional 
pipe flow in the annular regime.

\acknowledgement{The authors would like to thank Dr. Gustaaf Jacobs of the San Diego State University for his hospitality. This work was supported by a grant from the Simons Foundation ($\#426393$, David Kopriva). This work has been partially supported by Ministerio de Econom\'ia y Competitividad under the research grant (EUIN2017-88294, Gonzalo Rubio). This project has received funding from the European Union's Horizon 2020 research and innovation programme under grant agreement No 785549 (FireExtintion: H2020-CS2-CFP06-2017-01).
The authors acknowledge the computer resources and technical assistance provided by the Centro de Supercomputaci\'on y Visualizaci\'on de Madrid (CeSViMa).
}

\appendix

\section{Stability analysis of the alternative artificial compressibility model}\label{app:second-ACM}

We address the stability of the second artificial compressibility model 
\eqref{eq:governing:ACM-2}. Since changes with respect to the original model \eqref{eq:governing:ACM} 
only affect to time--derivative terms, 
the analysis performed for spatial terms hold in this model, and we need only include 
the new time--derivative terms.

\subsection{Continuous entropy analysis}

The entropy analysis performed in Sec.~\ref{subsec:Governing:iNS/CHEntropy} can be 
 extended to the second artificial compressibility method 
\eqref{eq:governing:ACM-2}. If we maintain the same entropy variables, all  of the 
steps performed for space operators hold, and only the temporal terms contraction needs to be recomputed. Doing so,
\begin{equation}
\begin{split}
& \mu c_t + 
  \sqrt{\rho}\svec{u}\cdot\left(\sqrt{\rho}\svec{u}\right)_t - \frac{p\svec{\nabla}^{2}p_t}{\rho_0 c_0^2} 
  = \mathcal F_{t} + \mathcal K_{t} + \frac{\svec{\nabla}p\cdot\svec{\nabla}p_{t}}{\rho_0 c_0^2}+\svec{\nabla}\cdot\left(\svec{f}_{t}^{\mathcal E} 
  - \frac{p\svec{\nabla}p_{t}}{\rho_0 c_0^2}\right)\\
  &=\mathcal F_{t} + \mathcal K_{t} + \frac{|\svec{\nabla}p|_{t}^{2}}{2\rho_0 c_0^2}+\svec{\nabla}\cdot\left(\svec{f}_{t}^{\mathcal E} -
   \frac{p\svec{\nabla}p_{t}}{\rho_0 c_0^2}\right)=\mathcal E_{2,t} + \svec{\nabla}\cdot\svec{f}_{t}^{\mathcal E_{2}},
  \end{split}
  \label{eq:continuous:second-model:time-derivative-contraction}
\end{equation}
with
\begin{equation}
  \mathcal E_{2} = \mathcal F_{t} + \mathcal K_{t} + 
  \frac{|\svec{\nabla}p|^2_{t}}{2\rho_0c_0^2},~~\svec{f}_{t}^{\mathcal 
  E_{2}}=\frac{3}{2}\sigma\varepsilon c_{t}\svec{\nabla}c - \frac{p\svec{\nabla}p_{t}}{\rho_0 
  c_0^2}.
  \label{eq:continuous:second-model:entropy-def}  
\end{equation}
The second artificial compressibility model is complemented with homogeneous 
Neumann boundary conditions for the pressure time derivative \cite{1996:Shen},
\begin{equation}
  \svec{\nabla}p_{t}\cdot\svec{n}\biggr|_{\partial \Omega} = 0,
\end{equation}
which makes the pressure term in $\svec{f}_{t}^{\mathcal E_{2}}$ vanish at 
the boundaries. Hence, the entropy equation with physical boundary conditions for the second artificial 
compressibility method is identical to the original model \eqref{eq:continuous:entropy-balance-w-bcs}
\begin{equation}
\frac{\diff}{\diff t}\left(\bar{\mathcal E}_{2}+\int_{\partial\Omega}f_{w}(c)\diff S\right) = -\int_{\Omega}\left(M_0|\svec{\nabla}\mu|^2 + 2\eta\tens{S}:\tens{S}\right)\diff\svec{x} \leqslant 
0,~~\bar{\mathcal E}_{2} = \int_{\Omega}\mathcal E_{2}\diff\svec{x}.
\label{eq:continuous:second-model:entropy-balance-w-bcs}
\end{equation}

\subsection{Discretization}

The discretization of the original artificial compressibility method \eqref{eq:dg:semi-discrete-approx}
translates to the second method, where only the pressure time derivative terms 
need to be updated. We write the first order term using 
${g}_{p}=\svec{\nabla}p$, transform the operators to local coordinates, 
integrate by parts, and replace the interface fluxes by a numerical flux and  integrals by  quadratures to get
\begin{equation}
\begin{split}
  \left\langle J\svec{\nabla}\cdot\svec{
  G}_{p,t},\varphi_{q}\right\rangle_{E}=  \left\langle \svec{\nabla}_{\xi}\cdot\svec{\tilde 
  G}_{p,t},\varphi_{q}\right\rangle_{E} 
  \approx  \int_{\partial 
  e,N}\varphi_{q}\svec{G}_{p,t}^{\star}\cdot\svec{n}\diff S - \left\langle\svec{\tilde 
  G}_{p,t}, \svec{\nabla}_{\xi}\varphi_{q}\right\rangle_{E,N}.
\end{split}
\end{equation}
For the interface fluxes, we use the BR1 method in interior faces, 
\begin{equation}
P^{\star}=\aver{P},~~\svec{G}_{p,t}^{\star}=\aver{\svec{G}_{p,t}},
\label{eq:ACM-2:BR1}
\end{equation}
and a homogeneous Neumann boundary condition in physical boundaries,
\begin{equation}
P^{\star}=P,~~\svec{G}_{p,t}^{\star}\cdot\svec{n}=0.
\label{eq:ACM-2:Neumann}
\end{equation}

\subsection{Semi--discrete stability}

As in the continuous analysis, we maintain the results from the original model, 
but update the time derivative terms \eqref{eq:stability:time-derivative-contraction-prep},
\begin{equation}
  \left\langle\mathcal JC_{t},\mu\right\rangle_{E,N} + \left\langle\mathcal 
  J\left(\sqrt{\rho}\svec{U}\right)_{t},\sqrt{\rho}\svec{U}\right\rangle_{E,N} - 
  \frac{1}{\rho_0c_0^2} \int_{\partial 
  e,N}P\svec{G}_{p,t}^{\star}\cdot\svec{n}\diff S +\frac{1}{\rho_0c_0^2} \left\langle\svec{\tilde 
  G}_{p,t}, \svec{\nabla}_{\xi}P\right\rangle_{E,N}.
  \label{eq:stability:second-model:time-derivative-contraction-prep}
\end{equation}
From the gradient definition \eqref{eq:dg:semi-discrete-approx:g}, with test 
function $\svec{\varphi}_{g}=\svec{G}_{p,t}$
\begin{equation}
  \left\langle \mathcal J\svec{G}_{p},\svec{G}_{p,t}\right\rangle_{E,N}= \int_{\partial e,N}\left(P^{\star}-P\right)\left(\svec{G}_{p,t}\cdot\svec{n}\right)\diff S + \left\langle 
  \svec{\tilde{G}}_{p,t},\svec{\nabla}_{\xi}P\right\rangle_{E,N},
\end{equation}
we get the time derivative terms,
\begin{equation}
\begin{split}
&  \left\langle \mathcal J\mathcal F_{t},1\right\rangle_{E,N}-\int_{\partial e,N}\frac{3}{2}\sigma\varepsilon \left(C_{t}\svec{{G}}_{c}^{\star}+C_{t}^{\star}\svec{G}_{c} - C_{t}\svec{G}_{c}\right)\cdot\svec{n}\diff 
  S+\left\langle\mathcal 
  J\left(\frac{1}{2}\rho V_{tot}^{2}\right)_{t},1\right\rangle_{E,N} \\
  &+ 
\left\langle \mathcal 
 J\left(\frac{1}{\rho_0c_0^2} 
 \svec{G}_{p}\right),\svec{G}_{p,t}\right\rangle_{E,N}-\frac{1}{\rho_0c_0^2}\int_{\partial 
 e,N}\left(P\svec{G}_{p,t}^{\star}+P^{\star}\svec{G}_{p,t}-P\svec{G}_{p,t}\right)\cdot\svec{n}\diff 
 S \\
=&  \left\langle \mathcal J\mathcal F_{t},1\right\rangle_{E,N}-\int_{\partial e,N}\frac{3}{2}\sigma\varepsilon \left(C_{t}\svec{{G}}_{c}^{\star}+C_{t}^{\star}\svec{G}_{c} - C_{t}\svec{G}_{c}\right)\cdot\svec{n}\diff 
  S+\left\langle\mathcal 
  J\mathcal K_{t},1\right\rangle_{E,N} \\
  &+ 
\left\langle \mathcal 
 J\left(\frac{1}{2\rho_0c_0^2} 
 |\svec{G}_{p}|_{t}^{2}\right),1\right\rangle_{E,N}-\frac{1}{\rho_0c_0^2}\int_{\partial 
 e,N}\left(P\svec{G}_{p,t}^{\star}+P^{\star}\svec{G}_{p,t}-P\svec{G}_{p,t}\right)\cdot\svec{n}\diff 
 S \\ 
 =&  \left\langle \mathcal J\mathcal E_{2,t},1\right\rangle_{E,N}-\int_{\partial e,N}\frac{3}{2}\sigma\varepsilon \left(C_{t}\svec{{G}}_{c}^{\star}+C_{t}^{\star}\svec{G}_{c} - C_{t}\svec{G}_{c}\right)\cdot\svec{n}\diff 
  S\\
  &-\frac{1}{\rho_0c_0^2}\int_{\partial 
 e,N}\left(P\svec{G}_{p,t}^{\star}+P^{\star}\svec{G}_{p,t}-P\svec{G}_{p,t}\right)\cdot\svec{n}\diff 
 S.
 \end{split}
\end{equation}
We have now obtained the time derivative of the alternative discrete entropy $\mathcal E_{2}$ \eqref{eq:continuous:second-model:entropy-def}. The first boundary integral belongs to the Cahn--Hilliard equation and was shown not to contribute at interior faces (Sec.~\ref{subsubsec:stability:boundary:interior-viscous-ch}) if $\beta=0$, and to modify the entropy if $\beta>0$. The surface free--energy is obtained in physical boundary faces (Sec.~\ref{subsubsec:stability:physical-boundary-wall}). The last integral was added by the discretization of the pressure time derivative Laplacian. The contribution from interior edges with the BR1 scheme \eqref{eq:ACM-2:BR1} is
\begin{equation}
\begin{split}
\mathrm{IBT}_{\nabla^{2}p_{t}} =& 
 \sum_{\interiorfaces}\int_{f,N}\left(\jump{P}\svec{G}_{p,t}^{\star}+\jump{\svec{G}_{p,t}}P^{\star}-\jump{P\svec{G}_{p,t}}\right)\cdot\svec{n}_{L}\mathrm{d}S\\
=&\sum_{\interiorfaces}\int_{f,N}\left(\jump{P}\aver{\svec{G}_{p,t}}+\jump{\svec{G}_{p,t}}\aver{P}-\jump{P\svec{G}_{p,t}}\right)\cdot\svec{n}_{L}\mathrm{d}S\\  
=&0.
\end{split}
\end{equation}
For physical boundaries, we apply the homogeneous Neumann boundary condition \eqref{eq:ACM-2:Neumann}, so that
\begin{equation}
\begin{split}
\mathrm{PBT}_{\nabla^{2}p_{t}}=&\sum_{\boundaryfaces}\int_{f,N}\left(P\svec{G}_{p,t}^\star+ P\svec{G}_{p,t}-P \svec{G}_{p,t}\right)\cdot\svec{n}\diff 
    S \\
    =&\sum_{\boundaryfaces}\int_{f,N}\left(0+P\svec{G}_{p,t}\svec{n}-P \svec{G}_{p,t}\right)\cdot\svec{n}\diff 
    S \\
    =&0.
\end{split}
\end{equation}

As a result, the alternative entropy \eqref{eq:continuous:second-model:entropy-def} satisfies a continuous bound \eqref{eq:continuous:second-model:entropy-balance-w-bcs}, and discrete entropy law,
\begin{equation}
    \frac{\diff}{\diff t}\left(\bar{\mathcal E}_{2}^{\beta}+\sum_{\boundaryfaces}\int_{f,N} f_{w}(C)\diff 
    S\right) \leqslant -\sum_{e}\left\langle \mathcal J\left(M_0|\svec{G}_{\mu}|^2+2\eta\mathcal 
  S:\mathcal S\right),1\right\rangle_{E,N}\le 0,
  \label{eq:ACM-2:entropy-eq-bc-final}
\end{equation}
where, similarly to \eqref{eq:stability:entropy-eq-bc-final}, the inequality depends on whether the ERS and or central fluxes are used, and the choice of the interface penalty parameter $\beta$.

\section{Entropy contraction of the inviscid fluxes}\label{app:inviscid-contraction}

We show the contraction of the inviscid fluxes \eqref{eq:continuous:inviscid-fluxes-contraction} 
for the iNS/CH system. To do so, we replace the entropy variables, fluxes, and 
non--conservative terms,
\begin{equation}
\begin{split}
\stvec{w}^{T}\left(\svec{\nabla}\cdot\ssvec{f}_{e}\left(\stvec{q}\right) + 
  \sum_{m=1}^{5}\ssvecg{\phi}_{m}\left(\stvec{q}\right)\cdot\svec{\nabla}w_{m}\right) 
 =&\mu\svec{\nabla}\cdot\left(c\svec{u}\right)+\frac{1}{2}u\svec{\nabla}\cdot\left(\rho u\svec{u}\right)+\frac{1}{2}v\svec{\nabla}\cdot\left(\rho v\svec{u}\right)
  +\frac{1}{2}w\svec{\nabla}\cdot\left(\rho w\svec{u}\right)+\svec{u}\cdot\svec{\nabla}p\\
  &+\frac{1}{2}\rho u\svec{u}\cdot\svec{\nabla}u+\frac{1}{2}\rho 
  v\svec{u}\cdot\svec{\nabla}v
  +\frac{1}{2}\rho 
  w\svec{u}\cdot\svec{\nabla}w+c\svec{u}\cdot\svec{\nabla}\mu+p\svec{\nabla}\cdot\svec{u}\\
=&  \svec{\nabla}\cdot\left(\left(\frac{1}{2}\rho v_{tot}^{2}+p+c\mu\right)\svec{u}\right)=
\svec{\nabla}\cdot\svec{f}^{\mathcal E}_{e}.
  \end{split}
\end{equation}
Thus, we see that the term $\svec{\nabla}\cdot\left(c\svec{u}\right)$ 
in the Cahn--Hilliard equation cancels the capillary pressure $c\svec{\nabla}\mu$ 
from the momentum non--conservative terms when multiplied by their 
respective entropy variables.
Similarly, it reveals that momentum terms from the conservative and non--conservative parts cancel each other, 
and that the pressure term in the momentum equation cancels the velocity divergence term by way of the 
artificial compressibility equation.

\section{Point--wise discretization}\label{app:point-wise}

In this section we list the steps to compute the solution time derivative with the DG approximation, \eqref{eq:dg:semi-discrete-approx}. To get the point--wise values, we replace the test function by the Lagrange polynomials $l_{i}\left(\xi\right)l_j\left(\eta\right)l_k\left(\zeta\right)$. Details on the extraction of point--wise values can be found in \cite{2009:Kopriva}.

\begin{enumerate}
	\item Compute the entropy variables, $\stvec{W}$. The first entropy variable is the chemical potential, $\mu$, which is computed from \eqref{eq:dg:semi-discrete-approx:c} and \eqref{eq:dg:semi-discrete-approx:mu},
	
\begin{equation}
\begin{split}
			\mathcal J_{ijk}\svec{G}_{c,ijk}=&\left(C^{\star}-C\right)_{ijk}\left(\frac{\delta_{im}}{w_{i}}\mathcal J\svec{a}^{1}_{ijk}+\frac{\delta_{jm}}{w_{j}}\mathcal J\svec{a}^{2}_{ijk}+\frac{\delta_{km}}{w_{k}}\mathcal J\svec{a}^{3}_{ijk}\right)\biggr|^{m=N}_{m=0}\\
			&+ \mathcal J\svec{a}^{1}_{ijk}\sum_{m=0}^{N}D_{im}C_{mjk}+ \mathcal J\svec{a}^{2}_{ijk}\sum_{m=0}^{N}D_{jm}C_{imk}+ \mathcal J\svec{a}^{3}_{ijk}\sum_{m=0}^{N}D_{km}C_{ijm}\\
			\mathcal J_{ijk}\mu_{ijk}=&\mathcal J_{ijk}F_0'(C_{ijk})-\frac{3}{2}\sigma\varepsilon\left(\frac{\delta_{im}}{w_{i}}{\tilde{G}}_{c,ijk}^{1,\star}+\frac{\delta_{jm}}{w_{j}}{\tilde{G}}_{c,ijk}^{2,\star}+\frac{\delta_{km}}{w_{k}}{\tilde{G}}_{c,ijk}^{3,\star}\right)\biggr|_{m=0}^{m=N}\\
			&+\frac{3}{2}\sigma\varepsilon\sum_{m=0}^{N}\left(\frac{w_{m}}{w_{i}}D_{mi}{\tilde{G}}_{c,mjk}^{1}+\frac{w_{m}}{w_{j}}D_{mj}{\tilde{G}}_{c,imk}^{2}+\frac{w_{m}}{w_{k}}D_{mk}{\tilde{G}}_{c,ijm}^{3}\right),		
\end{split}
\label{eq:pointwise:chemical-pot}\end{equation}
where $D_{ij}=l'_{j}(\xi_{i})$ is the derivative matrix, and recall the definition of contravariant components \eqref{eq:metrics:contravariant-flux}. We then get the entropy variables from the current values of the state vector, $\stvec Q$, and the chemical potential, $\mu$ computed by \eqref{eq:pointwise:chemical-pot},
\begin{equation}
\stvec{W}_{ijk} = \left(\mu,\frac{Q_{2}}{\sqrt{\rho}},\frac{Q_{3}}{\sqrt{\rho}},\frac{Q_{4}}{\sqrt{\rho}},\frac{Q_{5}}{\rho_0c_0^2}\right)_{ijk},~~\rho=\rho_{1}\hat{C}+\rho_{2}\left(1-\hat{C}\right),~~\hat{C}=\min\left(\max\left(C,0\right),1\right).
\end{equation}	

\item Compute the gradient of the entropy variables $\ssvec{G}$ from \eqref{eq:dg:semi-discrete-approx:g},
\begin{equation}
\begin{split}
 \mathcal J_{ijk}\ssvec{G}_{ijk}=&\left(\stvec{W}^{\star,T}-\stvec{W}^{T}\right)_{ijk}\left(\frac{\delta_{im}}{w_{i}}\mathcal J\svec{a}^{1}_{ijk}+\frac{\delta_{jm}}{w_{j}}\mathcal J\svec{a}^{2}_{ijk}+\frac{\delta_{km}}{w_{k}}\mathcal J\svec{a}^{3}_{ijk}\right)\biggr|^{m=N}_{m=0}\\
 &+ \mathcal J\svec{a}^{1}_{ijk}\sum_{m=0}^{N}D_{im}\stvec{W}_{mjk}+ \mathcal J\svec{a}^{2}_{ijk}\sum_{m=0}^{N}D_{jm}\stvec{W}_{imk}+ \mathcal J\svec{a}^{3}_{ijk}\sum_{m=0}^{N}D_{km}\stvec{W}_{ijm}.
\end{split}
\end{equation}

\item Compute the state vector time derivative,

\begin{equation}
\begin{split}
\mathcal J_{ijk}\smat{M}_{ijk}\frac{\diff\stvec{Q}_{ijk}}{\diff t}&+\left(\ssvec{F}_{e}^{\star}-\ssvec{F}_{e}+\sum_{m=1}^{5}\left(\ssvecg{\Phi}_{m}W_m\right)^{\diamondsuit}-\ssvecg{\Phi}_{m}W_{m}\right)_{ijk}\cdot\left(\frac{\delta_{im}}{w_{i}}\mathcal J\svec{a}^{1}_{ijk}+\frac{\delta_{jm}}{w_{j}}\mathcal J\svec{a}^{2}_{ijk}+\frac{\delta_{km}}{w_{k}}\mathcal J\svec{a}^{3}_{ijk}\right)\biggr|^{m=N}_{m=0}\\
&+\sum_{m=0}^{N}\left(D_{im}\tilde{\stvec{F}}_{e,mjk}^{1}+D_{jm}\tilde{\stvec{F}}_{e,imk}^{2}+D_{km}\tilde{\stvec{F}}_{e,ijm}^{3}\right)\\
&+\sum_{p=1}^{N}\sum_{m=1}^{5}\left(\tilde{\stvecg{\Phi}}_{m,ijk}^{1} D_{ip}W_{m,pjk}+\tilde{\stvecg{\Phi}}_{m,ijk}^{2} D_{jp}W_{m,ipk}+\tilde{\stvecg{\Phi}}_{m,ijk}^{3} D_{kp}W_{m,ijp}\right)\\
=&\left(\frac{\delta_{im}}{w_{i}}\tilde{\stvec{F}}_{v,ijk}^{1,\star}+\frac{\delta_{jm}}{w_{j}}\tilde{\stvec{F}}_{v,ijk}^{2,\star}+\frac{\delta_{km}}{w_{k}}\tilde{\stvec{F}}_{v,ijk}^{3,\star} \right)\biggr|_{m=0}^{m=N}\\
&+\sum_{m=0}^{N}\left(\frac{w_{m}}{w_{i}}D_{mi}\tilde{\stvec{F}}^{1}_{v,mjk}+\frac{w_{m}}{w_{j}}D_{mj}\tilde{\stvec{F}}^{2}_{v,imk}+\frac{w_{m}}{w_{k}}D_{mk}\tilde{\stvec{F}}^{3}_{v,ijm}\right) +\mathcal J_{ijk}\stvec{S}_{ijk}.
\end{split}
\end{equation}
\end{enumerate}

\bibliography{mybibfile}

\end{document}